\documentclass{amsart}
\usepackage{graphicx}
\pdfoutput=1
\allowdisplaybreaks[4]
\numberwithin{equation}{section}

\newfont{\Bigf}{cmr17}
\newfont{\bsft}{cmr10}
\newfont{\sft}{cmr7}
\newfont{\vsft}{cmr6}
\newfont{\vvsft}{cmr5}

\newcommand{\ts}{\textstyle}
\newcommand{\sss}{\scriptscriptstyle}

\def\tinysp{\mspace{1mu}}

\def\Wa{Wa\.{z}ewski}
\def\refHypoCr{\ref{hypo:C1}$^*$}
\def\refineqCr{\ref{hypo-ineqs}d}

\def\RR{{\mathbb{R}}}
\def\RRm{\RR^{m}}
\def\RRn{\RR^{n}}
\def\RRnm{\RR^{nm}}
\def\RRnmii{\RR^{nm^2}}
\def\RRnmj{\RR^{nm^j}}
\def\RRnmk{\RR^{nm^k}}

\def\RRmbym{\RR^{m\times m}}
\def\RRnbym{\RR^{n\times m}}
\def\RRnbyn{\RR^{n\times n}}
\def\RRnmbym{\RR^{nm\times m}}

\def\Torus{{\mathbb{T}}}
\def\BBoxd{{\mathbb{B}}_{d}}
\def\Wo{W^0}
\def\Wm{W^-}
\def\Ud{U_{\delta}}
\def\Udst{U_{\delta^*}}
\def\Udp{U_{\delta'}}
\def\Gaekd{\Gamma_{\epsilon, k, \delta}}
\def\Gaeked{\Gamma_{\epsilon, k_{\epsilon}, \delta}}
\def\Gaekedst{\Gamma_{\epsilon, k_{\epsilon}, \delta^*}}

\def\CC{{\mathcal{C}}}
\def\CF{{\mathcal{F}}}
\def\CG{{\mathcal{G}}}
\def\CI{{\mathcal{I}}}
\def\CK{{\mathcal{K}}}
\def\CL{{\mathcal{L}}}
\def\CM{{\mathcal{M}}}
\def\CO{{\mathcal{O}}}
\def\CP{{\mathcal{P}}}
\def\CQ{{\mathcal{Q}}}
\def\CR{{\mathcal{R}}}
\def\CV{{\mathcal{V}}}
\def\CW{{\mathcal{W}}}
\def\CX{{\mathcal{X}}}

\def\interior{\operatorname{int}}
\def\closure{\operatorname{cl}}
\def\textcap{\,{\ts \bigcap}\,}

\def\Mod{\operatorname{mod}}
\def\vect{\operatorname{vec}}
\def\vectnm{\vect_{\sss n\mspace{-1.5mu}, \mspace{-1.5mu}m}}
\def\vectmn{\vect_{\sss m\mspace{-1.5mu}, \mspace{-1.5mu}n}}
\def\vectnmm{\vect_{\sss n\mspace{-1.5mu} m\mspace{-1.5mu}, \mspace{-1.5mu}m}}
\def\vectnmjMm{\vect_{\sss n\mspace{-1.5mu} m^{\mspace{-1.5mu} j-1} \mspace{-3.5mu}, \mspace{-1.0mu}m}}
\def\vectnmkMm{\vect_{\sss n\mspace{-1.5mu} m^{\mspace{-1.5mu} k-1} \mspace{-3.5mu}, \mspace{-1.0mu}m}}

\def\Pibot{\Pi_{\sss{\bot}}}

\def\dFx{\dot{x} = F(x)}

\def\Ns{N^s}
\def\Nu{N^u}
\def\Me{M_{\epsilon}}
\def\Te{T_{\epsilon}}
\def\To{T_{0}}

\def\ddt{{\tfrac{d}{dt}}}

\def\au{\alpha_u}

\def\bp{\bar{p}}
\def\bq{\bar{q}}
\def\btheta{\bar{\theta}}

\def\Lfa{L_{f,a}}
\def\Lfz{L_{f,z}}
\def\Lga{L_{g,a}}
\def\Lgz{L_{g,z}}

\def\Lgp{L_{g,p}}

\def\Daf{D_{a}f}
\def\Dzf{D_{z}f}
\def\Dag{D_{a}\tinysp g}
\def\Dzg{D_{z}\tinysp g}

\def\DRf{D_{R}f}
\def\Dthetaf{D_{\theta}f}
\def\DRg{D_{R}\tinysp g}
\def\Dthetag{D_{\theta}\tinysp g}

\def\Drhof{D_{\rho}f}
\def\Drhog{D_{\rho}\tinysp g}

\def\xtil{\tilde{x}}
\def\atil{\tilde{a}}
\def\ztil{\tilde{z}}

\def\dx{\delta x}
\def\da{\delta a}
\def\dz{\delta z}

\def\xhat{\hat{x}}
\def\zhat{\hat{z}}

\def\bfx{{\boldsymbol x}}
\def\bfa{{\boldsymbol a}}
\def\bfz{{\boldsymbol z}}
\def\bfzhat{\hat{\bfz}}
\def\bfzero{{\boldsymbol 0}}

\def\Ibfz{\CI(\bfz)}
\def\Ibfzhat{\CI(\bfzhat)}

\def\Psix{\Psi_{x}}

\def\vi{v^{\sss 1}}
\def\Vi{V^{\sss 1}}
\def\zetai{\zeta^{\sss 1}}

\def\vii{v^{\sss 2}}
\def\Vii{V^{\sss 2}}
\def\zetaii{\zeta^{\sss 2}}

\def\vj{v^{\sss j}}
\def\zetaj{\zeta^{\sss j}}

\def\vjm{v^{\sss {j-1}}}
\def\zetajm{\zeta^{\sss {j-1}}}

\def\vk{v^{\sss k}}
\def\zetak{\zeta^{\sss k}}

\def\Dvifi{D_{\vi}\mspace{-2mu}f_1}
\def\Dzetaifi{D_{\zetai}\mspace{-2mu}f_1}

\def\Dviifii{D_{\vii}\mspace{-2mu}f_2}
\def\Dzetaiifii{D_{\zetaii}\mspace{-2mu}f_2}

\def\Dvjmfjm{D_{\vjm}\mspace{-2mu}f_{j-1}}
\def\Dzetajmfjm{D_{\zetajm}\mspace{-2mu}f_{j-1}}

\def\Dvkfk{D_{\vk}\mspace{-2mu}f_k}
\def\Dzetakfk{D_{\zetak}\mspace{-2mu}f_k}

\newtheorem{hypo}{Hypothesis}

\newtheorem{thm}{Theorem}[section]
\newtheorem{lem}[thm]{Lemma}
\newtheorem{prop}[thm]{Proposition}
\newtheorem{asp}[thm]{Assumption}

\newtheorem*{hypoCr}{Hypothesis \refHypoCr}
\newtheorem*{Waze}{Theorem ({\Wa})}

\theoremstyle{definition}

\newtheorem*{PosiInv}{Definition (Positive Invariance)}
\newtheorem*{defnWa}{Definition ({\Wa} Set)}

\begin{document}

\title{An Invariant Manifold Theory for ODEs and Its Applications}

\author{Dennis Guang Yang}
\address{Department of Mathematics, Cornell University, Ithaca, New York 14853}
\email{gy26@cornell.edu}


\date{September, 2009}



\begin{abstract}
For a system of ODEs defined on an open, convex domain $U$ containing a positively invariant set $\Gamma$, we prove that under appropriate hypotheses, $\Gamma$ is the graph of a $C^r$ function and thus a $C^r$ manifold. Because the hypotheses can be easily verified by inspecting the vector field of the system, this invariant manifold theory can be used to study the existence of invariant manifolds in systems involving a wide range of parameters and the persistence of invariant manifolds whose normal hyperbolicity vanishes when a small parameter goes to zero. We apply this invariant manifold theory to study three examples and in each case obtain results that are not attainable by classical normally hyperbolic invariant manifold theory.
\end{abstract}

\maketitle

\tableofcontents


\section{Introduction} \label{intro}

Consider the following ordinary differential equation
\begin{equation} \label{Fx}
\dot{x} = F(x) \,,
\end{equation}
where $\dot{}=\ddt$, $x \in \RRn$, and $F$ is $C^r$ for some $r \ge 1$. Let $\Phi(t,x)$ be the flow generated by (\ref{Fx}). Suppose $M$ is a $C^r$ invariant manifold of (\ref{Fx}). Following the work of Fenichel \cite{Fe71}, a simple version of normally hyperbolic invariant manifold theory can be roughly stated as follows:
{\it if there is a continuous splitting of the tangent bundle of\/ $\RRn$ restricted to\/ $M$: $T \RRn |_{M} = TM \oplus \Ns \oplus \Nu$ such that (1)\/ $TM \oplus \Ns$ and\/ $TM \oplus \Nu$ are invariant under\/ $D_x \Phi(t,M)$ (i.e., the linearization of\/ $\Phi(t,x)$ at\/ $M$) and (2)\/ $D_x \Phi(t,M)$ expands\/ $\Nu$ and contracts\/ $\Ns$ at rates at least\/ $r$ times of its expansion or contraction rate in\/ $TM$ ($r$-normal hyperbolicity), then\/ $M$ has a\/ $C^r$ stable manifold\/ $W^s(M)$ tangent to\/ $\Ns$ along\/ $M$ and a\/ $C^r$ unstable manifold\/ $W^u(M)$ tangent to\/ $\Nu$ along\/ $M$, and the manifolds\/ $M$, $W^s(M)$, and\/ $W^u(M)$ all persist with the same\/ $C^r$ smoothness under any sufficiently small (in\/ $C^1$ norm)\/ $C^r$ perturbation of the vector field\/ $F$.}
Here, we have omitted some technicalities such as the overflowing or inflowing invariance of $M$ (if it has a boundary) and the local invariance of $W^s(M)$ and $W^u(M)$. We refer readers to the work of Fenichel \cite{Fe71} and an extensive exposition by Wiggins \cite{Wi94} for the precise description of the theory and many other properties of normally hyperbolic invariant manifolds. Equivalent results can also be found in the work of Hirsch, Pugh, and Shub \cite{HiPuSh77}. A generalization to the case that $M$ is an invariant set of (\ref{Fx}) is given by Chow, Liu, and Yi in \cite{ChLiYi00}, where the authors proved the existence of center-stable, center-unstable, and center manifolds of $M$ and their persistence under small perturbations. In \cite{BaLuZe98}, Bates, Lu, and Zeng developed a normally hyperbolic invariant manifold theory for $C^1$ semiflows in general Banach spaces. They also provided a thorough review on the extensive history of invariant manifold theory.

The persistence of normally hyperbolic invariant manifolds can be applied to establish the existence of invariant manifolds in systems in the form of 
\begin{equation} \label{Fxe}
\dot{x} = F_0(x) + \epsilon F_1(x,\epsilon)
\end{equation}
provided that for the corresponding ``unperturbed'' system $\dot{x} = F_0(x)$, the existence of an ``unperturbed'' normally hyperbolic $C^r$ invariant manifold $M_0$ is already known. In particular, if $M_0$ is $r$-normally hyperbolic and the $C^1$ norm of the function $x \mapsto F_1(x,\epsilon)$ is $\CO(\epsilon)$ for $\epsilon \rightarrow 0^+$, then a $C^r$ invariant manifold $M_{\epsilon}$, which is the perturbed counterpart of $M_0$, exists for (\ref{Fxe}) with any $\epsilon \in (0, \epsilon_0]$ provided that $\epsilon_0$ is sufficiently small. However, it can be a difficult task to verify the normal hyperbolicity of $M_0$ since the precise knowledge about $M_0$ and the linearization of the unperturbed flow at $M_0$ is not attainable in many applications. Furthermore, because normally hyperbolic invariant manifold theory does not provide any further information about the size of $\epsilon_0$ other than being sufficiently small, we do not know to what extent the persistence result holds. 

In this paper, we present a new theory on the existence of $C^r$ invariant manifolds in systems of ordinary differential equations. Two special features of this new invariant manifold theory are: its hypotheses can be verified by inspecting the vector field instead of the flow of a specific system; and for systems depending on parameters, it provides a feasible way to compute the parameter ranges in which the desired invariant manifold results can be guaranteed. With these features, this new theory can work with systems that require delicate analysis of the perturbations involved as well as systems that cannot be treated as perturbation problems, and in both cases, it establishes results that are not attainable by classical invariant manifold theory.

\subsection{Main Results}

We consider the following system
\begin{equation} \label{sys}
\begin{split}
\dot a ={}& f(a,z) \,, \\
\dot z ={}& g(a,z) \,, 
\end{split}
\end{equation}
where $\dot{}=\ddt$ and $(a,z) \in \RRn \times \RRm$. For notational convenience, we define $x := (a,z)$ and $X := \RRn \times \RRm$. For any function defined on a subset of $X$, we use $x$ and $(a,z)$ interchangeably to denote its argument, e.g., $f(x) = f(a,z)$, and $g(x) = g(a,z)$. 

Suppose that $f$ and $g$ are at least $C^1$ on an open domain $U \subset \RRn \times \RRm$. Then (\ref{sys}) generates a flow $\Phi(t,x)$ on $U$ even though for some $x \in U$, $\Phi(t,x)$ may not be defined for all $t \in \RR$. For a subset of $U$, we define its positive invariance under the flow $\Phi$ as follows:

\begin{PosiInv}
$\Gamma \subseteq U$ is positively invariant under $\Phi$ if
\begin{enumerate}
   \item for every $x \in \Gamma$, $\Phi(t,x)$ is defined for all $t \ge 0$; and
   \item $\Phi(t,\Gamma) \subseteq \Gamma$ for all $t \ge 0$.
\end{enumerate}
\end{PosiInv}

First, we consider under what circumstances a positively invariant set $\Gamma$ is a $C^1$ manifold. Let $\langle \cdot\, , \cdot \rangle$ denote the usual Euclidean inner product, and let $\|\cdot\|$ denote the usual Euclidean norm as well as the induced operator norm. For each $x = (a,z) \in X$ and $d > 0$, let $\BBoxd(x)$ be a closed box neighborhood of $x$:
\begin{equation*}
\BBoxd(x) := \big\{ (a',z') \in X : \|a' - a\| \le d,\, \|z' - z\| \le d \big\} \,.
\end{equation*}
In addition, define $\CL : X \times X \rightarrow \RR$ as follows:
\begin{equation} \label{Lx}
\CL(x_1, x_2) := \| a_2 - a_1 \|^2 - \| z_2 - z_1 \|^2 \,,
\end{equation}
where $x_1=(a_1,z_1)$ and $x_2=(a_2,z_2)$. Then for each $x \in X$, we define the cone with vertex $x$ as follows:
\begin{equation} \label{cone}
\CC(x) := \big\{ x' \in X : \CL(x',x) \ge 0 \big\} \,.
\end{equation}

\begin{hypo} \label{hypo:U}
$U$ is an open, convex subset of\/ $\RRn \times \RRm$. In addition, there exists a\/ $d > 0$ such that\/ $\CC(x) \textcap U \subset \BBoxd(x)$ for all\/ $x \in U$.
\end{hypo}

\begin{hypo} \label{hypo:C1}
$f$ and\/ $g$ are\/ $C^1$ on\/ $U$. Furthermore, there exist a continuous, positive function\/ $\alpha : U \rightarrow \RR$, a continuous, nonnegative function\/ $\ell : U \rightarrow \RR$, and a constant\/ $c_1 > 0$ such that
\begin{subequations} \label{hypo-ineqs}
\begin{align}
& \langle a' , \Daf(x)\, a' \rangle \ge \alpha(x) \| a' \|^2 \text{ for any\/ $x \in U$ and any\/ $a' \in \RRn$,} \label{ineq:a-a} \\
& \langle z' , \Dzg(x)\, z' \rangle \le \ell(x) \| z' \|^2 \text{ for any\/ $x \in U$ and any\/ $z' \in \RRm$,} \label{ineq:z-z} \\ 
& \alpha(x) \ge \ell(x) + \|\Dzf(x)\| + \|\Dag(x)\| + c_1 \text{ for all\/ $x \in U$.} \label{ineq:C1}
\end{align}
\end{subequations}
\end{hypo}

\noindent Define the following projections:
\begin{align*}
\Pi :{}& (a,z) \mapsto a \,, \\
\Pibot :{}& (a,z) \mapsto z \,. 
\end{align*}

\begin{hypo} \label{hypo:Gamma}
$\Gamma \subset U$ is positively invariant under the flow of (\ref{sys}) and satisfy\/ $\Pibot(\Gamma) = \Pibot(U) \subseteq \RRm$. 
\end{hypo}

\begin{thm} \label{thm1}
Suppose Hypotheses \ref{hypo:U}, \ref{hypo:C1}, and \ref{hypo:Gamma} hold. Let\/ $K := \Pibot(\Gamma) = \Pibot(U)$. Then\/ $\Gamma$ contains all positively invariant subsets of\/ $U$, and there exists a\/ $C^1$ function\/ $h : K \rightarrow \RRn$ such that\/ $\Gamma = \big\{ (h(z), z) : z \in K \big\}$. Moreover, $\| h(z_2) - h(z_1) \| < \| z_2 - z_1 \|$ for any\/ $z_1$, $z_2 \in K$ with\/ $z_1 \neq z_2$, and\/ $\|Dh(z)\| < 1$ for all\/ $z \in K$.
\end{thm}

Next, we consider the further smoothness of a positively invariant set $\Gamma$ if we know {\it a priori} that $\Gamma$ is the graph of a $C^1$ function from an open subset of $\RRm$ to $\RRn$.

\begin{hypo} \label{hypo:K0}
$K_0$ is an open, convex subset of\/ $\RRm$. $h_0 : K_0 \rightarrow \RRn$ is\/ $C^1$, and there exists an\/ $\eta > 0$ such that\/ $\| Dh_0(z) \| < \eta$ for all\/ $z \in K_0$.
\end{hypo}

\begin{hypo} \label{hypo:Cr0}
There exists an open neighborhood\/ $U_0$ of\/ $\big\{ (h_0(z), z) : z \in K_0 \big\}$ such that\/ $f$ and\/ $g$ are\/ $C^r$ ($r \ge 2$) with their first to $r$-th derivatives all bounded on\/ $U_0$. Furthermore, there exist a continuous, positive function\/ $\alpha : U_0 \rightarrow \RR$, a continuous, nonnegative function\/ $\ell : U_0 \rightarrow \RR$, and a constant\/ $c_r > 0$ such that
\begin{subequations}
\begin{align}
& \langle a' , \Daf(x)\, a' \rangle \ge \alpha(x) \| a' \|^2 \text{ for any\/ $x \in U_0$ and any\/ $a' \in \RRn$,} \label{ineq:a-a:Cr0} \\
& \langle z' , \Dzg(x)\, z' \rangle \le \ell(x) \| z' \|^2 \text{ for any\/ $x \in U_0$ and any\/ $z' \in \RRm$,} \label{ineq:z-z:Cr0} \\ 
& \alpha(x) \ge r \ell(x) + ( r+1 ) \eta \| \Dag(x) \| + c_r \text{ for all\/ $x \in U_0$.} \label{ineq:Cr0}
\end{align}
\end{subequations}
\end{hypo}

\begin{thm} \label{thm2}
Suppose Hypotheses \ref{hypo:K0} and \ref{hypo:Cr0} hold. If\/ $\Gamma = \big\{ (h_0(z), z) : z \in K_0 \big\}$ is positively invariant under the flow of (\ref{sys}), then\/ $h_0 : K_0 \rightarrow \RRn$ is\/ $C^r$ with its first to $r$-th derivatives all bounded on\/ $K_0$. 
\end{thm}

It is evident that the combination of Theorem \ref{thm1} and Theorem \ref{thm2} establishes the existence of a $C^r$ positively invariant manifold for (\ref{sys}). Specifically, take $U_0 = U$, $h_0 = h$, and $\eta = 1$. Then by replacing Hypothesis \ref{hypo:C1} with a stronger one, we obtain a $C^r$ manifold theorem.

\begin{hypoCr}
$f$ and\/ $g$ are\/ $C^r$ ($r \ge 2$) with their first to $r$-th derivatives all bounded on\/ $U$. Furthermore, there exist a continuous, positive function\/ $\alpha : U \rightarrow \RR$, a continuous, nonnegative function\/ $\ell : U \rightarrow \RR$, and constants\/ $c_1$, $c_r > 0$ such that
\begin{align}
& \langle a' , \Daf(x)\, a' \rangle \ge \alpha(x) \| a' \|^2 \text{ for any\/ $x \in U$ and any\/ $a' \in \RRn$,} \notag \\
& \langle z' , \Dzg(x)\, z' \rangle \le \ell(x) \| z' \|^2 \text{ for any\/ $x \in U$ and any\/ $z' \in \RRm$,} \notag \\ 
& \alpha(x) \ge \ell(x) + \|\Dzf(x)\| + \|\Dag(x)\| + c_1 \text{ for all\/ $x \in U$,} \notag \\
& \alpha(x) \ge r \ell(x) + ( r+1 ) \| \Dag(x) \| + c_r \text{ for all\/ $x \in U$.} 
\tag{\ref{hypo-ineqs}d} 
\end{align}
\end{hypoCr}

\begin{thm} \label{thm3}
Suppose Hypotheses \ref{hypo:U}, \refHypoCr, and \ref{hypo:Gamma} hold. Let\/ $K := \Pibot(\Gamma) = \Pibot(U)$. Then\/ $\Gamma$ contains all positively invariant subsets of\/ $U$, and there exists a\/ $C^r$ function\/ $h : K \rightarrow \RRn$ such that\/ $\Gamma = \big\{ (h(z), z) : z \in K \big\}$. Moreover, $\| h(z_2) - h(z_1) \| < \| z_2 - z_1 \|$ for any\/ $z_1$, $z_2 \in K$ with\/ $z_1 \neq z_2$, $\|Dh(z)\| < 1$ for all\/ $z \in K$, and all higher derivatives of\/ $h$ are bounded on\/ $K$.
\end{thm}

To apply Theorem \ref{thm1} or Theorem \ref{thm3} to establish invariant manifolds for a given system of ordinary differential equations, we need to carry out the following five steps. The first step is to reformulate the system into the form of (\ref{sys}). It is often necessary to change variables and append auxiliary parameters. In particular, appropriate rescaling of variables can lead to sharper results from the inequalities (\ref{ineq:C1}) and (\refineqCr). The second step is to select a domain $U$ that satisfies Hypothesis \ref{hypo:U}. Sometimes, it is necessary to choose a family of candidate domains that depends on some parameters and then determine the parameters or their ranges in later analysis. The third step is to verify Hypothesis \ref{hypo:Gamma}. Very often we define $\Gamma \subset U$ to be the set of points whose images under the flow of the reformulated system stay in $U$ forever in forward time so that $\Gamma$ is positively invariant by definition. In addition, this definition of $\Gamma$ allows us to verify $\Pibot(\Gamma) = \Pibot(U)$ by simply inspecting the topological properties of the vector field at the boundary of $U$ and then applying an elementary topological argument---the {\Wa} principle (see Appendix \ref{appendix:Wa}). The fourth step is to verify Hypothesis \ref{hypo:C1} or \refHypoCr. Especially, using the inequalities (\ref{ineq:C1}) and (\refineqCr), we can estimate the set of admissible parameter values for the system or derive concrete, computable criteria for the smallness of perturbations. In the last step, we switch from the reformulated system back to the original system and identify the manifold that corresponds to $\Gamma$. Another advantage of the definition of $\Gamma$ introduced in the third step is that it allows us to easily establish important properties, which may include uniqueness, full invariance (both forward time and backward time), independence of the rescaling of variables in the reformulation of the system (in the first step), periodicity with respect to some variables, etc., for the manifold in the original system. In Section \ref{App}, we will illustrate all these five steps and the relevant technical arguments in detail with three examples. Specifically, the first example is not a perturbation problem and involves a wide range of parameter values, and the other two examples are related to the problem of weak hyperbolicity, which will be discussed in the next subsection.

Traditionally, proofs of invariant manifold theorems are based on either one of two complementary methods: Hadamard's graph transform method, and the Liapunov-Perron method. Both methods require the construction of a contraction in some Banach space such that its fixed point is a function whose graph is a Lipschitz invariant manifold. However, their approaches to constructing such a contraction are different. Specifically, the graph transform method obtains a contraction using the invariance of the manifold's graph representation under a time-$T$ map generated by the flow, whereas the Liapunov-Perron method achieves a contraction by setting up an integral equation using the variation of parameters formula. To prove the $C^1$ smoothness of the manifold, one constructs another contraction by formally differentiating the corresponding functional equations and then proves that its fixed point is in fact the desired derivative. If further smoothness is needed, one repeats this procedure inductively to show that the manifold is $C^r$ smooth. 

Alternatively, without referring to any contraction, one can establish invariant manifold results through a mixed use of invariant cones and topological arguments, an idea which, according to Jones \cite{Jo95}, was introduced by Charles Conley. In \cite{Jo95}, Jones used this method to construct a Lipschitz invariant manifold in a singularly perturbed slow-fast system. He also gave a very brief outline for showing the $C^1$ smoothness of the manifold. In an earlier reference \cite{McGehee73}, based on the same mixed use of invariant cones and topological arguments, McGehee supplied a proof of a $C^r$ local stable manifold theorem for a fixed point of a hyperbolic linear map plus a nonlinear term that is small in $C^r$ norm. In addition, Bates and Jones \cite{BaJo89} extended this method to an infinite-dimensional setting to construct Lipschitz invariant manifolds for a semilinear partial differential equation. 

In this paper, the proof of Theorem \ref{thm1}, which consists of two parts---the Lipschitz smoothness and the $C^1$ smoothness of $\Gamma$, is based on the mixed use of invariant cones and topological arguments. In fact, the Lipschitz smoothness part of the proof closely follows the relevant part of the proof given in \cite{Jo95}. However, in the $C^1$ smoothness part of the proof, we introduce a new strategy for combining invariant cones and topological arguments. Specifically, we first construct a vector bundle over $\Gamma$ such that it is invariant under the linearized dynamics along $\Gamma$ and also satisfies some desirable topological and dynamical properties, and then we show that this vector bundle is in fact the tangent bundle of $\Gamma$. This approach is completely different from what is described in the outline for showing $C^1$ smoothness in \cite{Jo95} and what is given in the relevant part of the proof in \cite{McGehee73}. The main technical issue is that the strategies used in \cite{Jo95,McGehee73} rely on the assumption of ``sufficient'' $C^r$-smallness of certain terms (or perturbations) and thus are not applicable in our case, where the corresponding boundedness conditions are given explicitly by the inequality (\ref{ineq:C1}). 

In the proof of Theorem \ref{thm2} for the $C^r$ smoothness of $\Gamma$, we also introduce a new inductive scheme, which utilizes the hypothesis that the manifold is known to be the graph of the $C^1$ function $h_0$. We show that in an appropriate coordinate system, along any solution trajectory of (\ref{sys}) on $\Gamma$ the dynamics of the derivative $D h_0$ are governed by a system in the form of (\ref{sys}). This allows us to apply Theorem \ref{thm1} to show that $D h_0$ is $C^1$ and thus $h_0$ is $C^2$. Once again, by choosing an appropriate coordinate system, we can express the dynamics (along any solution trajectory of (\ref{sys}) on $\Gamma$) of the second derivative $D^2 h_0$ by a system in the form of (\ref{sys}). Then we proceed inductively to establish the $C^r$ smoothness of $h_0$.

\subsection{The Problem of Weak Hyperbolicity} \label{intro:WH}

Although the invariant manifold theory presented in this paper has a wide range of applications, the initial motivation for this work is to study the persistence of invariant manifolds whose normal hyperbolicity is ``weak'' in the sense that it depends on a small parameter and vanishes as the parameter goes to zero. Specifically, consider the system
\begin{equation} \label{WH}
\begin{split}
\dot{w} ={}& \epsilon F_1(w) + \epsilon^2 F_2(w,\theta,\epsilon) \,, \\
\dot{\theta} ={}& G_0 (w) + \epsilon G_1(w) + \epsilon^2 G_2(w,\theta,\epsilon) \,,
\end{split}
\end{equation}
where $w \in \RR^n$, $\theta \in \Torus^m$, $0 < \epsilon \ll 1$, and all functions on the right-hand side are $C^r$ ($r \ge 1$) with respect to their arguments. Systems in the form of (\ref{WH}) often arise in situations involving averaging. In particular, in view of the absence of $\theta$ in the $\CO(1)$ and $\CO(\epsilon)$ terms, (\ref{WH}) can be regarded as the result of a first-order averaging procedure applied to a near-integrable system which is formulated in action-angle variables. In this case, the equation $\dot{w} = \epsilon F_1(w)$ is usually referred to as the averaged equation, and studying the dynamics of it is the first step to understand the dynamics of the full system (\ref{WH}). Suppose a $C^r$ invariant manifold $M \subset \RR^n$ is identified for the averaged equation $\dot{w} = \epsilon F_1(w)$. It follows that $M \times \Torus^m$ is a $C^r$ invariant manifold in the truncated system
\begin{equation} \label{WHtr}
\begin{split}
\dot{w} ={}& \epsilon F_1(w) \,, \\
\dot{\theta} ={}& G_0 (w) + \epsilon G_1(w) \,.
\end{split}
\end{equation}
Then the immediate question is the persistence of $M \times \Torus^m$ in the full system (\ref{WH}), which includes the $\CO(\epsilon^2)$ terms. However, this is a nontrivial problem even under the assumption that $M \times \Torus^m$ is normally hyperbolic with respect to (\ref{WHtr}) for any fixed $\epsilon \in (0, \epsilon_0]$ with $\epsilon_0$ being some positive constant. In particular, when we consider any fixed $\epsilon \in (0, \epsilon_0]$, it is unclear whether or not the $\CO(\epsilon^2)$ terms are ``sufficiently small'' in $C^1$ norm so that the aforementioned persistence of normally hyperbolic invariant manifolds applies. On the other hand, as we reduce the $C^1$ norm of the $\CO(\epsilon^2)$ terms by reducing $\epsilon$, the ``strength'' of the normal hyperbolicity of $M \times \Torus^m$ with respect to (\ref{WHtr}) is also weakened since the hyperbolicity is generated by the $\CO(\epsilon)$ term $\epsilon F_1(w)$. In the limit $\epsilon \rightarrow 0^+$, the normal hyperbolicity of $M \times \Torus^m$ with respect to (\ref{WHtr}) fails as $\epsilon F_1(w)$ vanishes. This situation is referred to as {\it weak hyperbolicity}\/ of the manifold $M \times \Torus^m$ \cite{Ha99,Wi94}. We also say that $M \times \Torus^m$ is only {\it weakly normally hyperbolic}\/ with respect to (\ref{WHtr}). Note that by rescaling $\theta$ and time, we can obtain a variant of (\ref{WH}) and correspondingly a variant of (\ref{WHtr}), with respect to which the normal hyperbolicity of $M \times \Torus^m$ no longer depends on $\epsilon$. However, the problem remains as other forms of singularities occur after rescaling variables. See \cite{ChLi00} for the related discussions.
  
In \cite{Wi94}, Wiggins adapted a continuation argument, which was originally proposed by Kopell \cite{Ko85} for a different class of systems, to show the persistence of $M \times \Torus^m$ in (\ref{WH}) for the case that $m \ge 1$ and $M$ is an attracting fixed point of $\dot{w} = \epsilon F_1(w)$. Specifically, fix $\epsilon > 0$, and rewrite (\ref{WH}) as a one-parameter family of systems
\begin{equation} \label{WH-aux}
\begin{split}
\dot{w} ={}& \epsilon F_1(w) + \alpha^2 F_2(w,\theta,\alpha) \,, \\
\dot{\theta} ={}& G_0 (w) + \epsilon G_1(w) + \alpha^2 G_2(w,\theta,\alpha) \,,
\end{split}
\end{equation}
where $\alpha \in [0, \epsilon]$ is an auxiliary parameter. When $\alpha = 0$, (\ref{WH-aux}) reduces to (\ref{WHtr}), under the flow of which $M \times \Torus^m$ is normally hyperbolic by the assumption that $M$ is an attracting fixed point of $\dot{w} = \epsilon F_1(w)$. Then the persistence of normally hyperbolic invariant manifolds implies that there exists an $\alpha_1\in(0, \epsilon]$ such that for any $\alpha \in [0, \alpha_1]$, (\ref{WH-aux}) has an invariant torus 
$\CM_{\alpha}$ which is a perturbation of $\CM_{0} := M \times \Torus^m$. Note that if $\CM_{\alpha_1}$ is normally hyperbolic with respect to (\ref{WH-aux}) with $\alpha = \alpha_1$, then $\alpha_1$ can be increased to some $\alpha_2 \in (\alpha_1, \epsilon]$ such that $\CM_{\alpha}$ exists for any $\alpha \in [0, \alpha_2]$. Now the question is whether or not this process can be repeated so that the family of invariant tori $\big\{ \CM_{\alpha} \big\}$ can be continued for all $\alpha \in [0, \epsilon]$. In \cite{Wi94} (pp. 168--170), Wiggins derived upper estimates of 
the generalized Lyapunov-type numbers\footnote{We refer readers to the references \cite{Fe71,Wi94} for the technical definitions of generalized Lyapunov-type numbers. Roughly speaking, under the linearized dynamics along a trajectory on the invariant manifold $M$, the generalized Lyapunov-type numbers measure the expansion and contraction rates in the bundles $\Nu$ and $\Ns$ and compare these rates with the expansion or contraction rate in the tangent bundle $TM$. Thus, the invariant manifold $M$ is normally hyperbolic if (1) $M$ is $C^r$ smooth for some $r \ge 1$, (2) it admits the splitting $T \RRn |_{M} = TM \oplus \Ns \oplus \Nu$ which satisfies the required invariance condition, and (3) the generalized Lyapunov-type numbers of $M$ are bounded below certain critical values. Note that (1) and (2) are required so that the generalized Lyapunov-type numbers of $M$ can be defined.} 
of $\CM_{\alpha}$. Then he argued that for sufficiently small $\epsilon > 0$, those upper estimates are bounded below the required critical values uniformly for all $\alpha \in [0, \epsilon]$ and thus $\CM_{\alpha}$ exists and remains normally hyperbolic for any $\alpha \in [0, \epsilon]$. 

This continuation argument was then used as a general strategy for establishing the persistence of weakly normally hyperbolic invariant manifolds in several different model systems (see, e.g., \cite{HaWi96,Ha99}). However, Chicone and Liu subsequently pointed out in \cite{ChLi00} that this argument contains a conceptual gap because it fails to address the issue that {\it the uniform boundedness of the generalized Lyapunov-type numbers does not guarantee uniform normal hyperbolicity for a continuous family of normally hyperbolic invariant manifolds}. Referring to the example discussed above, the technical issue can be described as follows. By extending the interval $[0, \alpha_1]$ to $[0, \alpha_2]$ and continuing in the same way, we can obtain an increasing sequence $\big\{ \alpha_1, \alpha_2, \alpha_3, ... \big\}$ such that $\CM_{\alpha}$ exists for any $\alpha \in [0, \alpha_i]$ with $i=1,2,3,...$. However, if $\big\{ \alpha_1, \alpha_2, \alpha_3, ... \big\}$ converges to a limit $\alpha^* < \epsilon$, we can only conclude that $\CM_{\alpha}$ exists for any $\alpha \in [0, \alpha^*)$. In order to continue the family of tori $\big\{ \CM_{\alpha} \big\}$ for $\alpha \ge \alpha^*$, we must establish the existence of a normally hyperbolic invariant torus $\CM_{\alpha^*}$ for (\ref{WH-aux}) with $\alpha = \alpha^*$. Unfortunately, even if the generalized Lyapunov-type numbers of $\CM_{\alpha}$ are bounded below their critical values uniformly for all $\alpha \in [0, \alpha^*)$, the existence of $\CM_{\alpha^*}$ cannot be guaranteed because $\CM_{\alpha}$ may still lose normal hyperbolicity in the limit $\alpha \rightarrow \alpha^{*-}$ by 
losing smoothness locally\footnote{In \cite{ChLi00}, Chicone and Liu presented a scenario that a family of normally hyperbolic limit cycles, with their generalized Lyapunov-type numbers uniformly bounded below the required critical values, cannot be continued further as the family converges to a nonsmooth homoclinic loop.} 
or losing the splitting $T (\RRn \times \Torus^m) |_{\CM_{\alpha}} = T\CM_{\alpha} \oplus \Ns \oplus \Nu$ locally\footnote{In \cite{HaLl07}, Haro and de la Llave studied continuation of invariant tori for quasi-periodic perturbations of the standard map. They numerically observed and analyzed a situation that a one-parameter family of normally hyperbolic invariant 1-tori, with their generalized Lyapunov-type numbers uniformly bounded below the required critical values, cannot be continued further as the subbundles $\Ns$ and $\Nu$ converge locally when the continuation parameter approaches a critical value. We remark that similar examples can also be constructed for flows.}. 
In particular, since the generalized Lyapunov-type numbers only reflect the global characteristics of the linearized dynamics along the invariant manifold, their uniform boundedness alone is not sufficient to rule out these local ``defects''. Note that although the expressions of the upper estimates used in the continuation argument may remain defined and bounded below the required critical values for $\alpha \ge \alpha^*$, they no longer have any significance since the corresponding generalized Lyapunov-type numbers are not even defined before we actually establish the existence of $\CM_{\alpha}$ for $\alpha \ge \alpha^*$.

Technically, the presence of the conceptual gap in the continuation argument does not imply that the result achieved by this argument is necessarily false. However, the assertion in \cite{Wi94} that an attracting invariant torus $\CM_{\alpha}$ exists for (\ref{WH-aux}) with any $\alpha \in [0, \epsilon]$ if $\epsilon$ is sufficiently small is in fact untrue, and so is the claim that for a sufficiently small $\epsilon > 0$, certain upper estimates of the generalized Lyapunov-type numbers of $\CM_{\alpha}$ are bounded below the required critical values uniformly for all $\alpha \in [0, \epsilon]$. Indeed, it is possible that the generalized Lyapunov-type numbers converge to their critical values at some $\alpha < \epsilon$ no matter how small $\epsilon$ is. To see this, consider the following simple example, which is a particular case of (\ref{WH}):
\begin{equation} \label{WH-CE}
\begin{split}
\dot{w} ={}& - \epsilon w + \alpha^2 \sin\theta \,, \\
\dot{\theta} ={}& w \,,
\end{split}
\end{equation}
where $w \in \RR$, $\theta \in \Torus^1$, $0 < \epsilon \ll 1$, and $\alpha \in [0, \epsilon]$. Clearly, for any fixed $\epsilon > 0$, $w = 0$ is an attracting fixed point of $\dot{w} = - \epsilon w$, and $\CM_0 = \big\{ (0,\theta) : \theta \in \Torus^1 \big\}$ is a normally hyperbolic invariant 1-torus for (\ref{WH-CE}) with $\alpha = 0$. However, no matter how small $\epsilon > 0$ is, (\ref{WH-CE}) with $\alpha = \epsilon$ has two fixed points, which are a saddle point at $(0,0)$ and a stable spiral at $(0,\pi)$, and thus possesses no invariant 1-torus. Furthermore, in this example the generalized Lyapunov-type numbers of $\CM_{\alpha}$ do reach their critical values as $\alpha$ converges to $\frac{1}{2}\epsilon$ from below, and the family of $C^1$ 1-tori $\big\{ \CM_{\alpha} \big\}$ only exists for $\alpha \in [0, \frac{1}{2}\epsilon)$.

The above counterexample also suggests that we should not expect the general persistence of weakly normally hyperbolic invariant manifolds in (\ref{WH}) and instead we should formulate additional hypotheses about the system in order to establish the desired invariant manifold results. For example, let $M$ be a hyperbolic fixed point or a hyperbolic periodic orbit of the averaged system $\dot{w} = \epsilon F_1(w)$. Suppose $\theta \in \Torus^1$ and $G_0 (w) \neq 0$ for any $w \in M$. Then $M \times \Torus^1$ persists in (\ref{WH}) as a hyperbolic periodic orbit or a hyperbolic torus (see, e.g., \cite{ArKoNe97,BoMi61,GuHo90}). In addition, Chow and Lu \cite{ChLu95} studied a particular form of (\ref{WH}):
\begin{equation} \label{sys:ChLu}
\begin{split}
\dot{w} ={}& \epsilon F_1(w) + \epsilon^2 F_2(w,\theta,\epsilon) \,, \\
\dot{\theta} ={}& \Theta_0 + \epsilon G_1(w) + \epsilon^2 G_2(w,\theta,\epsilon) \,,
\end{split}
\end{equation}
where $\theta \in \Torus^m$ and the frequency vector $\Theta_0$ is constant. For the case that $M$ is a fixed point of $\dot{w} = \epsilon F_1(w)$ with center-stable, center-unstable, and center manifolds $W^{\text{cs}}(M)$, $W^{\text{cu}}(M)$, and $W^{\text{c}}(M)$, respectively, they proved the persistence of $W^{\text{cs}}(M) \times \Torus^m$, $W^{\text{cu}}(M) \times \Torus^m$, and $W^{\text{c}}(M) \times \Torus^m$ in (\ref{sys:ChLu}) for sufficiently small $\epsilon > 0$. In \cite{ChLi00}, Chicone and Liu studied a 3-dimensional system (see (\ref{sys:ChLi})) which, by rescaling an angular variable and time, can be put into the following form:
\begin{equation} \label{sys:ChLi:res}
\begin{split}
\dot{w} ={}& \epsilon F_1(w) + \epsilon^{\frac{3}{2}} F_2(w,\theta,\epsilon^{\frac{1}{2}}) \,, \\
\dot{\theta} ={}& 1 \,,
\end{split}
\end{equation}
where $w \in \RR \times \Torus^1$, $\theta \in \Torus^1$, and $\epsilon$ is the square of the original small parameter used in \cite{ChLi00}. Chicone and Liu proved that if $M$ is an attracting or repelling limit cycle of $\dot{w} = \epsilon F_1(w)$ then $M \times \Torus^1$ persists in (\ref{sys:ChLi:res}) for sufficiently small $\epsilon > 0$. Note that the perturbation is $\CO(\epsilon^{\frac{3}{2}})$ in this case.

In Subsection \ref{ex3}, we will study the following system
\begin{equation} \label{GWH0}
\begin{split}
\dot{w} ={}& \epsilon F_1(w) + \epsilon^{1+\mu} F_2(w,\theta,\epsilon) \,, \\
\dot{\theta} ={}& \Theta_0 + \epsilon^{\nu} G_1(w,\epsilon) + \epsilon^{1+\gamma} G_2(w,\theta,\epsilon) \,,
\end{split} 
\end{equation}
where $w \in \RR^n$, $\theta \in \Torus^m$, $\mu > 0$, $\gamma > 0$, $1 \ge \nu \ge 0$, and $\mu + \nu > 1$. Note that both (\ref{sys:ChLu}) and (\ref{sys:ChLi:res}) are special cases of (\ref{GWH0}). Specifically, we have $\mu = 1$, $\gamma = 1$, and $\nu = 1$ for (\ref{sys:ChLu}) and $\mu = \frac{1}{2}$, $\gamma = 1$, and $\nu = 1$ for (\ref{sys:ChLi:res}) (considering $\dot{\theta} = 1 + \epsilon 0 + \epsilon^2 0$). We will apply Theorem \ref{thm1} and Theorem \ref{thm3} to show that if $M$ is a hyperbolic periodic orbit of $\dot{w} = \epsilon F_1(w)$ then $M \times \Torus^m$ persists in (\ref{GWH0}) for sufficiently small $\epsilon > 0$ (see Theorem \ref{thm4}). The crucial difference between (\ref{GWH0}) (or its special cases (\ref{sys:ChLu}) and (\ref{sys:ChLi:res})) and (\ref{WH-CE}) with $\alpha = \epsilon$ is that for the latter case, $\mu = 1$ and $\nu = 0$ and thus the inequality $\mu + \nu > 1$ is not satisfied. Finally, we mention that in Subsection \ref{ex2}, we will study (\ref{sys:ChLi:res}) in its original form in \cite{ChLi00} and solve an open problem posed by Chicone and Liu. In particular, by applying Theorem \ref{thm1} and Theorem \ref{thm3}, we will formulate a sufficient condition for the existence of a $C^r$ invariant torus without assuming the existence of any unperturbed invariant torus.

\subsection{Organization}

The balance of this paper is organized as follows. The complete proofs of Theorem \ref{thm1} and Theorem \ref{thm2} are given in Section \ref{Proof:thm1} and Section \ref{smoothness:Cr}, respectively. Then we illustrate the applications of Theorem \ref{thm1} and Theorem \ref{thm3} with three examples in Section \ref{App}. At the end of this paper, a concise statement of the {\Wa} principle is included in Appendix \ref{appendix:Wa}.

\section{Proof of Theorem \ref{thm1}} \label{Proof:thm1}

\subsection{The Existence and the Lipschitz Continuity of $h$} \label{smoothness:Lip}

In this subsection, we prove that there is a Lipschitz function $h : K \rightarrow \RRn$ such that the graph of $h$ is $\Gamma$. We achieve this goal in three steps. In the first step (Lemma \ref{lem:cone:inv}), we establish the ``invariance'' of the ``moving cones'' $\CC( \Phi(t,x) )$, which move in translation as their vertices move under the flow $\Phi$ in forward time.
In the second step (Lemma \ref{lem:cone:exp}), we show that trajectories in these cones drift away from the moving vertices at least exponentially fast in forward time. In the third step (Lemma \ref{lem:h:Lip}), we establish the existence and the Lipschitz continuity of $h$ using Lemmas \ref{lem:cone:inv} and \ref{lem:cone:exp}. Throughout this subsection, we will use the following notations for any $x_1$, $x_2 \in U$: 
\begin{align*}
(\atil_1(t), \ztil_1(t)) ={}& \xtil_1(t) := \Phi(t,x_1) \,, \\
(\atil_2(t), \ztil_2(t)) ={}& \xtil_2(t) := \Phi(t,x_2) \,, \\
(\da(t), \dz(t)) ={}& \dx(t) := \xtil_2(t) - \xtil_1(t) \,. 
\end{align*}
In addition, taking account of the convexity of $U$, we define the functions $\Lfa$, $\Lfz$, $\Lga$, and $\Lgz$, all mapping $U \times U$ into $\RR$, as follows:
\begin{equation} \label{Lfg}
\begin{split}
\Lfa(x_1,x_2) :={}& {\ts\int^1_0} \alpha( x_1 + s(x_2 - x_1) )\, ds \,, \\
\Lfz(x_1,x_2) :={}& {\ts\int^1_0} \big\| \Dzf( x_1 + s(x_2 - x_1) ) \big\|\, ds \,, \\
\Lga(x_1,x_2) :={}& {\ts\int^1_0} \big\| \Dag( x_1 + s(x_2 - x_1) ) \big\|\, ds \,, \\
\Lgz(x_1,x_2) :={}& {\ts\int^1_0} \ell( x_1 + s(x_2 - x_1) )\, ds \,. 
\end{split}
\end{equation}
Then it follows from (\ref{ineq:a-a}) and (\ref{ineq:z-z}) of Hypothesis \ref{hypo:C1} that 
\begin{equation}\label{est}
\begin{split}
\langle \da , \dot{\da} \rangle \ge{}& \Lfa(\xtil_1, \xtil_2) \| \da \|^2 - \Lfz(\xtil_1, \xtil_2) \| \da \| \| \dz \| \,, \\
\langle \dz , \dot{\dz} \rangle \le{}& \Lga(\xtil_1, \xtil_2) \| \da \| \| \dz \| + \Lgz(\xtil_1, \xtil_2) \| \dz \|^2 \,,
\end{split}
\end{equation}
whenever $\xtil_1(t) \in U$ and $\xtil_2(t) \in U$. 

\begin{lem} \label{lem:cone:inv}
For any\/ $x_1$, $x_2 \in U$ and\/ $t_0 > 0$, if\/ $x_2 \in \CC( x_1 )$, $\Phi([0,t_0],x_1) \subset U$, and\/ $\Phi([0,t_0],x_2) \subset U$, then\/ $\xtil_2(t) \in \CC( \xtil_1(t) )$ for all\/ $t \in [0,t_0]$.
\end{lem}

\begin{proof} The lemma is trivially true for $x_1 = x_2 \in U$. Thus, we only consider $x_1$, $x_2 \in U$ with $x_1 \neq x_2$. Then $\xtil_1(t) \neq \xtil_2(t)$ for all possible $t$ due to the uniqueness property of the solutions of (\ref{sys}) inside $U$. Suppose that $\Phi([0,t_0],x_1) \subset U$ and $\Phi([0,t_0],x_2) \subset U$ for some $t_0 > 0$.

Assume that $\xtil_2(\tau) \in \CC( \xtil_1(\tau) )$ for some $\tau \in [0, t_0]$. Then $\|\dz(\tau)\| \le \|\da(\tau)\|$ by (\ref{Lx}) and (\ref{cone}). Incorporate this relation into (\ref{est}). It follows that
\begin{subequations}
\begin{align}
\langle \da(\tau) , \dot{\da}(\tau) \rangle \ge{}& \big( \Lfa(\xtil_1(\tau), \xtil_2(\tau)) - \Lfz(\xtil_1(\tau), \xtil_2(\tau)) \big) \| \da(\tau) \|^2 \,, \label{est:a} \\
\langle \dz(\tau) , \dot{\dz}(\tau) \rangle \le{}& \big( \Lga(\xtil_1(\tau), \xtil_2(\tau)) + \Lgz(\xtil_1(\tau), \xtil_2(\tau)) \big) \| \da(\tau) \|^2 \,. 
\end{align}
\end{subequations}
Combining the above two inequalities and using (\ref{Lfg}) and (\ref{ineq:C1}), we obtain
\begin{equation*}
\langle \da(\tau) , \dot{\da}(\tau) \rangle - \langle \dz(\tau) , \dot{\dz}(\tau) \rangle \ge c_1 \| \da(\tau) \|^2 \,. 
\end{equation*}

Thus, for any $\tau \in [0, t_0]$ such that $\xtil_2(\tau) \in \CC( \xtil_1(\tau) )$ (i.e., $\CL ( \xtil_2(\tau), \xtil_1(\tau) ) \ge 0$), we have that
\begin{equation}
\big[ \ddt \CL ( \xtil_2(t), \xtil_1(t) ) \big]_{\!\sss t=\tau} \ge 2 c_1 \| \da(\tau) \|^2 > 0 \,, \label{dLdt}
\end{equation}
where the second (strict) inequality is assured by $\|\da(\tau)\| > 0$, a fact due to $\xtil_1(\tau) \neq \xtil_2(\tau)$ and $\|\da(\tau)\| \ge \| \dz(\tau) \|$. Now consider $x_2 \in \CC( x_1 )$. Then
\begin{equation}
\CL ( \xtil_2(0), \xtil_1(0) ) \ge 0 \,, \label{L0}
\end{equation}
and with $\tau = 0$, (\ref{dLdt}) becomes
\begin{equation}
\big[ \ddt \CL ( \xtil_2(t), \xtil_1(t) ) \big]_{\!\sss t=0} > 0 \,. \label{dLdt0}
\end{equation}
Altogether, (\ref{dLdt})--(\ref{dLdt0}) imply that $\CL ( \xtil_2(t), \xtil_1(t) ) \ge 0$ for all $t \in [0,t_0]$ if $x_2 \in \CC( x_1 )$.
\end{proof}

We now demonstrate that inside the moving cone $\CC( \Phi(t,x) )$, trajectories drift away from the moving vertex $\Phi(t,x)$ at least exponentially fast in forward time. 

\begin{lem} \label{lem:cone:exp}
For any\/ $x_1$, $x_2 \in U$ and\/ $t_0 > 0$, if\/ $x_2 \in \CC( x_1 )$, $\Phi([0,t_0],x_1) \subset U$, and\/ $\Phi([0,t_0],x_2) \subset U$, then\/ $\| \atil_2(t) - \atil_1(t) \| \ge \| a_2 - a_1 \|\, e^{c_1 t}$ for all\/ $t \in [0, t_0]$.
\end{lem}

\begin{proof}
Take $x_1$, $x_2 \in U$ and $t_0>0$ such that $\Phi([0,t_0],x_1) \subset U$ and $\Phi([0,t_0],x_2) \subset U$. Assume that $\xtil_2(\tau) \in \CC( \xtil_1(\tau) )$ for some $\tau \in [0, t_0]$. Recall (\ref{est:a}), and note that $\Lfa(x,x') - \Lfz(x,x') \ge c_1$ for any $(x,x') \in U \times U$ by (\ref{ineq:C1}). Then we have
\begin{equation}
\big[ {\ddt} \| \atil_2(t)-\atil_1(t) \|^2 \big]_{\!\sss t=\tau} \ge 2 c_1 \| \atil_2(\tau)-\atil_1(\tau) \|^2 \,. \label{ddadt}
\end{equation}

Now consider $x_2 \in \CC( x_1 )$. By Lemma \ref{lem:cone:inv}, $\xtil_2(\tau) \in \CC( \xtil_1(\tau) )$ for all $\tau \in [0,t_0]$. Thus (\ref{ddadt}) holds for all $\tau \in [0,t_0]$. It follows that $\| \atil_2(t)-\atil_1(t) \| \ge \| a_2 - a_1 \|\, e^{c_1 t}$ for all $t \in [0,t_0]$.
\end{proof}

Having established Lemmas \ref{lem:cone:inv} and \ref{lem:cone:exp}, we now show that $\Gamma$ is the largest positively invariant subset of $U$ and it is the graph of a Lipschitz function $h : K \rightarrow \RRn$.

\begin{lem} \label{lem:h:Lip}
$\Gamma$ contains all positively invariant subsets of\/ $U$, and there exists a Lipschitz function\/ $h : K \rightarrow \RRn$ such that\/ $\Gamma = \big\{ (h(z),z) : z \in K \big\}$. Moreover, $\| h(z_2) - h(z_1) \| < \| z_2 - z_1 \|$ for any\/ $z_1$, $z_2 \in K$ with\/ $z_1 \neq z_2$.
\end{lem}

\begin{proof}
Let $\Gamma'$ be the union of all positively invariant subsets of $U$. Consider $x_1$, $x_2 \in \Gamma'$ with $x_1 \neq x_2$. First, we prove $x_2 \not\in \CC(x_1)$ by contradiction. Suppose $x_2 \in \CC(x_1)$. Since $\Gamma'$ is positively invariant, $\Phi([0,t_0],x_1) \subset \Gamma' \subseteq U$ and $\Phi([0,t_0],x_2) \subset \Gamma' \subseteq U$ for all $t_0 > 0$. Applying Lemma \ref{lem:cone:exp}, we obtain that for all $t \ge 0$,
\begin{equation*}
\| \atil_2(t)-\atil_1(t) \| \ge \| a_2 - a_1 \| \, e^{c_1 t} \,.
\end{equation*}
In addition, $\| a_2 - a_1 \| > 0$ since $x_1 \neq x_2$ and $\| a_2 - a_1 \| \ge \| z_2 - z_1 \|$ for $x_2 \in \CC(x_1)$. Thus, $\| a_2 - a_1 \| \, e^{c_1 t}$ increases exponentially (as opposed to being constant $0$) for $t \ge 0$. On the other hand, Lemma \ref{lem:cone:inv} implies that $\xtil_2(t) \in \CC( \xtil_1(t) )$ for all $t \ge 0$. Then it follows from Hypothesis \ref{hypo:U} that $\| \atil_2(t)-\atil_1(t) \| < d$ for all $t \ge 0$, which constitutes a contradiction.

It follows that for each $x \in \Gamma'$, $\CC(x) \textcap \Gamma' = \big\{ x \big\}$. Since $K = \Pibot(\Gamma) = \Pibot(U)$ by Hypothesis \ref{hypo:Gamma} and $\Pibot(\Gamma) \subseteq \Pibot(\Gamma') \subseteq \Pibot(U)$, we have that $\Pibot(\Gamma') = K$. Then for each $z \in K$ there is a unique $x_z = (a_z,z) \in \Gamma'$ such that $\Pibot(x_z)=z$. Thus, we can define a function $h : K \rightarrow \RRn$ by letting $h(z) := a_z$ for each $z \in K$. Obviously, $\Gamma' = \big\{ (h(z),z) : z \in K \big\}$. In addition, we have $\Gamma = \Gamma'$ since $\Pibot(\Gamma) = \Pibot(\Gamma') = K$. 

Next, we consider $x_1 = ( h(z_1), z_1 )$ and $x_2 = ( h(z_2), z_2 )$ for any $z_1$, $z_2 \in K$ with $z_1 \neq z_2$. Since $x_1$, $x_2 \in \Gamma$ and $x_1 \neq x_2$, we have $x_2 \not\in \CC(x_1)$. Then it follows from (\ref{Lx}) and (\ref{cone}) that $\| h(z_2) - h(z_1) \| < \| z_2 - z_1 \|$.
\end{proof}

\subsection{The $C^1$ Smoothness of $h$}\label{smoothness:C1}

We need to work with the variational equation of (\ref{sys}) along trajectories in $\Gamma$. Rewrite (\ref{sys}) in a compact form: $\dFx$, and let $\bfa \in \RRn$, $\bfz \in \RRm$, and $\bfx=(\bfa,\bfz) \in X$ be the variations of $a$, $z$, and $x$, respectively. Then the variational equation of (\ref{sys}) along a trajectory $\Phi(t,x)$ for any $x \in \Gamma$ is given by
\begin{equation}
\dot{\bfx} = DF(\Phi(t,x))\, \bfx \,, \label{sys:vari}
\end{equation}
where $t \ge 0$ and $x \in \Gamma$ now serves as a parameter. 

For the linear system (\ref{sys:vari}) with any parameter $x \in \Gamma$, the solution that originates at $\bfx$ at $t=0$ can be represented as $Q(t,x)\,\bfx$, where $Q(t,x)$ is a linear transformation of $X$ to itself for each fixed $t \ge 0$ with $Q(0,x) = I$ being the identity transformation of $X$. Note that for $f$ and $g$ that are $C^1$ on $U$, $DF(\Phi(t,x))$ depends on $(t,x)$ continuously on $[0, \infty) \times \Gamma$. 
Then for each fixed $x \in \Gamma$, $Q(t,x)$ is defined for all $t \ge 0$, and furthermore, the map $(t,x) \mapsto Q(t,x)$ is continuous on $[0, \infty) \times \Gamma$. 

It is important to note that $Q(t,x)\,\bfx$ is {\it not}\/ a flow on $X$ since the system (\ref{sys:vari}) is {\it nonautonomous}. However, for $(t,x) \in [0,\infty) \times \Gamma$, the family of linear transformations $Q(t,x)$ forms a {\it cocycle}\/ over the flow $\Phi(t,x)$, i.e., for any $\tau$, $t \in [0,\infty)$ and $x \in \Gamma$, 
\begin{align} 
Q(0,x) ={}& I \,, \notag \\
Q(\tau + t,x) ={}& Q(\tau,\Phi(t,x)) \circ Q(t,x) \label{t+tau} \,.
\end{align}

Note that our construction of cones in $X$ still applies when using the variational variable $\bfx=(\bfa,\bfz)$, i.e., 
\begin{align} 
\CL(\bfx_1,\bfx_2) ={}& \| \bfa_2 - \bfa_1 \|^2 - \| \bfz_2 - \bfz_1 \|^2 \,, \label{Lx:vari} \\
\CC(\bfx) ={}& \big\{ \bfx' \in X : \CL(\bfx',\bfx) \ge 0 \big\} \,.  \label{cone:vari}
\end{align}
Let $\bfzero$ be the zero vector in $X$. For each $x \in \Gamma$, define the set $T(x)$ as follows:
\begin{equation} \label{LSTx}
\begin{split}
T(x) :={}& \big\{ \bfx \in X : Q(t,x)\,\bfx \not\in \interior(\CC( \bfzero )) \text{ for all } t \ge 0 \big\} \\
={}& \big\{ \bfx \in X : \CL( Q(t,x)\,\bfx, \bfzero ) \le 0 \text{ for all } t \ge 0 \big\} \,.
\end{split} 
\end{equation}
Obviously, $\bfzero \in T(x)$, and the image of $T(x)$ under the linear transformation $Q(t,x)$ does not intersect $\interior(\CC( \bfzero ))$ for any $t \ge 0$. In addition, it follows from (\ref{t+tau}) that for any $\bfx \in T(x)$,
\begin{equation} \label{LSTx:inv}
Q(t,x)\,\bfx \in T( \Phi(t,x) ) \text{ for all } t \ge 0 \,.
\end{equation}

The outline of the proof of the $C^1$ smoothness of $h$ is the following. In \ref{LinSpaTx}, we will first show that for each $x \in \Gamma$, $T(x)$ is in fact a linear subspace of $X$ and it can be represented as the graph of a linear operator $H(x) : \RRm \rightarrow \RRn$. Then we will demonstrate in \ref{H:continuous} that the map $x \mapsto H(x)$ is continuous for all $x \in \Gamma$. Finally, in \ref{H:derivative}, we will show that $H(h(z),z)$ is indeed the derivative of $h$ at $z$ for any $z \in K$.


\subsubsection{$T(x)$ is a linear subspace of\/ $X$}\label{LinSpaTx}

Consider (\ref{sys:vari}) in terms of $(\bfa,\bfz)$, i.e.,
\begin{equation}\label{sys:vari:a-z}
\begin{split}
\dot{\bfa} ={}& \Daf(\Phi(t,x))\, \bfa + \Dzf(\Phi(t,x))\, \bfz \,, \\
\dot{\bfz} ={}& \Dag(\Phi(t,x))\, \bfa + \Dzg(\Phi(t,x))\, \bfz \,.
\end{split}
\end{equation}
Applying (\ref{ineq:a-a}) and (\ref{ineq:z-z}), we obtain the following inequalities:
\begin{subequations} \label{est:vari}
\begin{align}
\langle \bfa , \dot{\bfa} \rangle \ge{}& \alpha(\Phi(t,x)) \|\bfa\|^2 - \big\| \Dzf(\Phi(t,x)) \big\| \|\bfa\| \|\bfz\| \,, \label{est:vari:bfa} \\
\langle \bfz , \dot{\bfz} \rangle \le{}& \big\| \Dag(\Phi(t,x)) \big\| \|\bfa\| \|\bfz\| + \ell(\Phi(t,x)) \|\bfz\|^2 \,, \label{est:vari:bfz}
\end{align}
\end{subequations}
which are similar to (\ref{est}). In addition, {\it for any\/ $x \in \Gamma$, both inequalities of (\ref{est:vari}) hold for all\/ $t \ge 0$ since\/ $\Phi(t,x) \in \Gamma \subset U$ for all\/ $t \ge 0$.}

The next three lemmas are simple consequences of (\ref{est:vari}).

\begin{lem} \label{lem:C1:cone:inv}
For any\/ $x \in \Gamma$ and any\/ $\bfx \in \CC( \bfzero )$, $Q(t,x)\,\bfx \in \CC( \bfzero )$ for all\/ $t \ge 0$.
\end{lem}

\begin{lem} \label{lem:C1:est:a}
For any\/ $x \in \Gamma$ and any\/ $\bfx \in \CC( \bfzero )$,
\begin{gather*}
\big\| \Pi( Q(t,x)\,\bfx ) \big\| \ge \| \Pi( \bfx ) \|\, e^{\int^t_0 ( \alpha(\Phi(s,x)) - \| \Dzf(\Phi(s,x)) \| )\, ds} \text{ for all\/ $t \ge 0$.}
\end{gather*}
\end{lem}

We omit the proofs of Lemma \ref{lem:C1:cone:inv} and Lemma \ref{lem:C1:est:a} since they are essentially the same as the proofs of Lemma \ref{lem:cone:inv} and Lemma \ref{lem:cone:exp}, respectively. A difference here is that we only consider $x \in \Gamma$. Thus the statements of Lemma \ref{lem:C1:cone:inv} and Lemma \ref{lem:C1:est:a} can be shown true for all $t \ge 0$ (as opposed to only for $t \in [0,t_0]$). 

By (\ref{ineq:C1}), we have that for any $(t,x) \in [0, \infty) \times \Gamma$,
\begin{equation} \label{ineq:exp-rates}
\alpha(\Phi(t,x)) - \big\| \Dzf(\Phi(t,x)) \big\| \ge c_1 + \big\| \Dag(\Phi(t,x)) \big\| + \ell(\Phi(t,x)) \,.
\end{equation}
Thus Lemma \ref{lem:C1:est:a} implies that for any $x \in \Gamma$ and $\bfx \in \CC( \bfzero )$ with $\bfx \neq \bfzero$, $\big\| \Pi( Q(t,x)\, \bfx ) \big\|$ and $\| Q(t,x)\, \bfx \|$ both grow at least exponentially fast as $t$ increases.

Next, for $\bfx \in T(x)$, we have the following growth estimate of $\| Q(t,x)\, \bfx \|$.

\begin{lem} \label{lem:C1:est:LSTx}
For any\/ $x \in \Gamma$ and any\/ $\bfx \in T(x)$, 
\begin{gather*}
\| Q(t,x)\, \bfx \| \le 2 \| \bfx \|\, e^{\int^t_0 ( \| \Dag(\Phi(s,x)) \| + \ell(\Phi(s,x)) )\, ds} \text{ for all\/ $t \ge 0$.}
\end{gather*}
\end{lem}

\begin{proof}
Take an arbitrary $x \in \Gamma$ and then an arbitrary $\bfx \in T(x)$. 
Recall (\ref{Lx:vari}) and (\ref{cone:vari}) for $\CC( \bfzero )$. By the definition of $T(x)$ (see (\ref{LSTx})), we have that 
\begin{align}
\big\| \Pi( Q(t,x)\, \bfx ) \big\| \le{}& \big\| \Pibot( Q(t,x)\, \bfx ) \big\| \text{ and} \label{bfa-vs-bfz} \\
\| Q(t,x)\, \bfx \| \le{}& 2 \big\| \Pibot( Q(t,x)\, \bfx ) \big\|  \label{bfx-vs-bfz}
\end{align}
for all $t \ge 0$. Then it follows from (\ref{est:vari:bfz}) and (\ref{bfa-vs-bfz}) that for all $t \ge 0$,
\begin{gather*}
{\ddt} \big\| \Pibot( Q(t,x)\, \bfx ) \big\|^2 \le 2 \Big( \big\| \Dag(\Phi(t,x)) \big\| + \ell(\Phi(t,x)) \Big) \big\| \Pibot( Q(t,x)\, \bfx ) \big\|^2 \,. 
\end{gather*}
Thus, for all $t \ge 0$,
\begin{equation*}
\big\| \Pibot( Q(t,x)\, \bfx ) \big\| \le \| \Pibot( \bfx ) \|\, e^{\int^t_0 ( \| \Dag(\Phi(s,x)) \| + \ell(\Phi(s,x)) )\, ds} \,.
\end{equation*}
Then using (\ref{bfx-vs-bfz}), we obtain that for all $t \ge 0$,
\begin{align*}
\| Q(t,x)\, \bfx \| \le{}& 2 \| \Pibot( \bfx ) \|\, e^{\int^t_0 ( \| \Dag(\Phi(s,x)) \| + \ell(\Phi(s,x)) )\, ds} \\
\le{}& 2 \| \bfx \|\, e^{\int^t_0 ( \| \Dag(\Phi(s,x)) \| + \ell(\Phi(s,x)) )\, ds} \,. \qedhere
\end{align*}
\end{proof}

We are now ready to demonstrate that $T(x)$ is a linear subspace of $X$.

\begin{lem} \label{lem:C1:subspace}
For each\/ $x \in \Gamma$, $T(x)$ is a linear subspace of\/ $X$.
\end{lem}

\begin{proof}
Take an arbitrary $x \in \Gamma$. Consider any $\bfx_1$, $\bfx_2 \in T(x)$ and any $\lambda_1$, $\lambda_2 \in \RR$. We prove that $\bfx_3 = \lambda_1 \bfx_1 + \lambda_2 \bfx_2 \in T(x)$ by contradiction.

Suppose $\bfx_3 \not\in T(x)$. Then there exists a $t_0 \ge 0$ such that $Q(t_0,x)\,\bfx_3 \in \interior(\CC( \bfzero ))$. Applying Lemma \ref{lem:C1:est:a} to $x' := \Phi(t_0,x) \in \Gamma$ and $\bfx' := Q(t_0,x)\,\bfx_3$, we obtain that for all $\tau \ge 0$,
\begin{equation*}
\big\| \Pi( Q(\tau, x')\,\bfx' ) \big\| \ge \| \Pi( \bfx' ) \|\, e^{\int^{\tau}_0 ( \alpha(\Phi(s,x')) - \| \Dzf(\Phi(s,x')) \| )\, ds} \,,
\end{equation*}
which, by (\ref{t+tau}), can be rewritten in terms of $x$ and $\bfx_3$ as follows: 
\begin{gather*}
\big\| \Pi( Q( \tau + t_0, x )\,\bfx_3 ) \big\| \ge \big\| \Pi( Q(t_0,x)\,\bfx_3 ) \big\|\, e^{\int^{\tau}_0 ( \alpha(\Phi(s+t_0,x)) - \| \Dzf(\Phi(s+t_0,x)) \| )\, ds} \,.
\end{gather*}
Let $t = \tau + t_0$. Then the above estimate becomes
\begin{equation} \label{bfx3:exp}
\big\| \Pi( Q(t,x)\,\bfx_3 ) \big\| \ge \big\| \Pi( Q(t_0,x)\,\bfx_3 ) \big\|\, e^{\int^t_{t_0} ( \alpha(\Phi(s,x)) - \| \Dzf(\Phi(s,x)) \| )\, ds} 
\end{equation}
for all $t \ge t_0$. Note that $\big\| \Pi( Q(t_0,x)\,\bfx_3 ) \big\| > 0$ since $Q(t_0,x)\,\bfx_3 \in \interior(\CC( \bfzero ))$. Thus, the right-hand side of (\ref{bfx3:exp}) increases at least exponentially (as opposed to being constant $0$) for $t \ge t_0$.

On the other hand, it follows from Lemma \ref{lem:C1:est:LSTx} that for all $t \ge t_0$,
\begin{align*}
\| Q(t,x)\,\bfx_3 \| \le{}& |\lambda_1| \| Q(t,x)\,\bfx_1 \| + |\lambda_2| \| Q(t,x)\,\bfx_2 \|  \\
\le{}& 2 ( |\lambda_1| \| \bfx_1 \| + |\lambda_2| \| \bfx_2 \| ) \, e^{\int^t_0 ( \| \Dag(\Phi(s,x)) \| + \ell(\Phi(s,x)) )\, ds} \\
\le{}& \mu \, e^{\int^t_{t_0} ( \| \Dag(\Phi(s,x)) \| + \ell(\Phi(s,x)) )\, ds} \,,
\end{align*}
where $\mu$ in the last inequality is a sufficiently large constant. By (\ref{ineq:exp-rates}), this estimate contradicts (\ref{bfx3:exp}).
\end{proof}

Let $\Ibfz$ denote the cross-section of $X$ at $\bfz$ along the $\bfa$-direction, i.e.,
\begin{equation*}
\Ibfz := \big\{ (\bfa,\bfz) : \bfa \in \RRn \big\} \,. 
\end{equation*}
A consequence of Lemma \ref{lem:C1:subspace} is that $T(x) \textcap \Ibfz$ contains at most one point.

\begin{lem} \label{lem:C1:at-most-1}
For any\/ $x \in \Gamma$, $T(x) \textcap \Ibfz$ contains at most one point for each\/ $\bfz \in \RRm$.
\end{lem}

\begin{proof}
Take an arbitrary $x \in \Gamma$. Suppose there exist $\bfx_1 \neq \bfx_2$ such that $\bfx_1$, $\bfx_2 \in T(x) \textcap \Ibfzhat$ for some $\bfzhat \in \RRm$. Then $\bfx_3 = \bfx_2 - \bfx_1 \neq \bfzero$, and, by Lemma \ref{lem:C1:subspace}, $\bfx_3 \in T(x)$. On the other hand,
\begin{equation*}
\| \Pibot(\bfx_3) \| = \| \Pibot(\bfx_2) - \Pibot(\bfx_1) \| = \|\bfzhat - \bfzhat\| =0 \,.
\end{equation*}
Since $\bfx_3 \neq \bfzero$, we have $\| \Pi(\bfx_3) \| > \| \Pibot(\bfx_3) \| = 0$. Then $\bfx_3 \in \interior(\CC( \bfzero ))$. This constitutes a contradiction.
\end{proof}

In fact, $T(x) \textcap \Ibfz$ contains exactly one point for every $\bfz \in \RRm$.

\begin{lem} \label{lem:C1:exactly-1}
For any\/ $x \in \Gamma$, $T(x) \textcap \Ibfz$ contains exactly one point for each\/ $\bfz \in \RRm$.
\end{lem}

\begin{proof}
Take an arbitrary $x \in \Gamma$. By Lemma \ref{lem:C1:at-most-1}, we only need to prove that $T(x) \textcap \Ibfz \neq \emptyset$ for any $\bfz \in \RRm$. We will use the {\Wa} theorem (see Appendix \ref{appendix:Wa}), which requires us to work with a {\it flow}. Thus, we append $\dot \tau =1$ to (\ref{sys:vari}) with the argument $t$ in $\Phi(t,x)$ replaced by $\tau$ to form an autonomous system
\begin{align*}
\dot \tau ={}& 1 \,, \\
\dot{\bfx} ={}& DF(\Phi(\tau,x))\, \bfx \,,
\end{align*}
where $(\tau, \bfx) \in (-\epsilon_x, \infty) \times X$ with $\epsilon_x > 0$ being a constant that depends on the chosen $x \in \Gamma$. Let $\Psix(t,(\tau,\bfx))$ be the flow generated by the above system:
\begin{equation*} 
\Psix(t,(\tau,\bfx)) = \big( t + \tau,\, Q(t,\Phi(\tau,x))\,\bfx \big) \,.
\end{equation*}

We define a set $\CW \subset (-\epsilon_x, \infty) \times X$ as follows:
\begin{align}
\CW :={}& \big\{ (\tau, \bfx) : \tau \ge 0,\, \bfx \not\in \interior(\CC(\bfzero)) \big\} \notag\\
={}& \big\{ (\tau, \bfx) : \tau \ge 0,\, \CL( \bfx, \bfzero ) \le 0 \big\} \,. \label{CW}
\end{align}
Since $\CW$ is a closed subset of $(-\epsilon_x, \infty) \times X$, it automatically satisfies the condition (W1) in the definition of a {\Wa} set (see Appendix \ref{appendix:Wa}), i.e., if $(\tau, \bfx) \in \CW$ and $\Psix([0,t],(\tau,\bfx)) \subset \closure(\CW)$ then $\Psix([0,t],(\tau,\bfx)) \subset \CW$. Next, we define the sets $\CW^0$ and $\CW^-$ as follows:
\begin{align*}
\CW^0 :={}& \big\{ (\tau, \bfx) \in \CW : \exists\, t >0 \text{ such that } \Psix(t, (\tau,\bfx)) \not\in \CW \big\} \,, \\
\CW^- :={}& \big\{ (\tau, \bfx) \in \CW : \Psix([0,t), (\tau,\bfx)) \not\subseteq \CW \text{ for all } t > 0 \big\} \,. 
\end{align*}
Note that $\CW^0$ is the set of points that do not stay in $\CW$ forever under the flow $\Psix$ in forward time, and $\CW^-$ is the set of points that immediately leave $\CW$ in forward time. Clearly, $\CW^- \subseteq \CW^0 \subseteq \CW$. In order to verify that $\CW$ is a {\Wa} set, we need to check the condition (W2), that is, $\CW^-$ is closed relative to $\CW^0$. 

It is obvious that $\CW^- \subseteq \partial\CW$, which is the boundary of $\CW$ and is the union of two disjoint sets $S_1 := \big\{ (\tau, \bfx) :  \tau \ge 0,\, \CL( \bfx, \bfzero ) = 0 \big\}$ and $S_2 := \big\{ (0, \bfx) : \CL( \bfx, \bfzero ) < 0 \big\}$. For any $(\tau,\bfx) \in S_1$ with $\bfx \neq \bfzero$, by the fact that $\Phi(\tau,x) \in \Gamma \subset U$, (\ref{est:vari}) holds, and it leads to
\begin{equation*}
\Big[ {\ddt} \CL \big( Q(t,\Phi(\tau,x))\,\bfx , \bfzero \big) \Big]_{\!\sss t=0} > 0 \,.
\end{equation*}
Thus, $(\tau,\bfx) \in \CW^-$ if $(\tau, \bfx) \in S_1$ and $\bfx \neq \bfzero$. For each $(0,\bfx) \in S_2$, by the continuity of $\CL$, there exists a $t_0>0$ such that $\CL \big( Q(t,\Phi(0,x))\,\bfx , \bfzero \big) < 0$ for all $t \in [0, t_0]$. It follows that $\Psix(t,(0,\bfx)) = \big( t,\, Q(t,\Phi(0,x))\,\bfx \big) \in \CW$ for all $t \in [0, t_0]$. Thus, $\CW^- \textcap S_2 = \emptyset$. Furthermore, note that for any $(\tau, \bfzero) \in S_1$, $\Psix(t, (\tau,\bfzero)) = (t + \tau, \bfzero) \in S_1 \subset \CW$ for all $t \ge 0$. Altogether, we have simultaneously
\begin{gather}
\CW^0 \textcap \big\{ (\tau, \bfzero) : \tau \ge 0 \big\} = \emptyset \,, \label{empty-CW0} \\
\CW^- = S_1 \setminus \big\{ (\tau, \bfzero) : \tau \ge 0 \big\} 
= \big\{ (\tau, \bfx) :  \tau \ge 0,\, \bfx \neq \bfzero,\, \CL( \bfx, \bfzero ) = 0 \big\} \,. \label{CW-}
\end{gather}
In addition, by (\ref{empty-CW0}) and (\ref{CW-}), we have
\begin{equation*}
\CW^0 \setminus \CW^- \subseteq \big\{ (\tau, \bfx) :  \tau \ge 0,\, \CL( \bfx, \bfzero ) < 0 \big\} \,.
\end{equation*}
Thus, for every $(\tau,\bfx) \in \CW^0 \setminus \CW^-$, there is a neighborhood $\CV$, which contains $(\tau,\bfx)$ and is open in $\RR \times X$, such that $\CV \textcap \CW^- = \emptyset$. Then, $\CW^0 \setminus \CW^-$ is open relative to $\CW^0$, and $\CW^-$ is closed relative to $\CW^0$. Therefore, $\CW$ is a {\Wa} set. By the {\Wa} theorem, there exists a continuous function $\CR : \CW^0 \times [0,1] \rightarrow \CW^0$ such that $\CR$ is a strong deformation retraction of $\CW^0$ onto $\CW^-$.

We now prove that $T(x) \textcap \Ibfz \neq \emptyset$ for any $\bfz \in \RRm$. Since $\bfzero \in T(x)$, we only need to show that $T(x) \textcap \Ibfz \neq \emptyset$ for any $\bfz \in \RRm$ with $\|\bfz\| > 0$. We prove this by contradiction.

Assume that $T(x) \textcap \Ibfzhat = \emptyset$ for some $\bfzhat \in \RRm$ with $\|\bfzhat\| > 0$. Let $\hat{R} := \|\bfzhat\|$. Then by the definition of $T(x)$ (see (\ref{LSTx})) and (\ref{CW}), we have
\begin{equation*}
\big\{ (0, (\bfa,\bfzhat)) : \|\bfa\| \le \hat{R} \big\} = \big\{ (0, \bfx) : \bfx \in \Ibfzhat \big\} \textcap \CW \subset \CW^0 \,.
\end{equation*}
Thus the set $\big\{ (0, (\bfa,\bfzhat)) : \|\bfa\| \le \hat{R} \big\}$ is contained inside the domain of the continuous function $\CR(\,\cdot\,,1) : \CW^0 \rightarrow \CW^-$.

Moreover, for all $(\tau, \bfx) \in \CW^-$, $\|\Pi(\bfx)\| = \|\Pibot(\bfx)\| > 0$ according to (\ref{CW-}). Thus we can define a projection $\CP : \CW^- \rightarrow \big\{ \bfa \in \RRn : \|\bfa\|=\hat{R} \big\}$ as follows:
\begin{equation*}
\CP (\tau,\bfx) := \hat{R} \tfrac{\Pi(\bfx)}{\|\Pi(\bfx)\|} \,.
\end{equation*}
By taking the composition of $\CP$, $\CR(\,\cdot\,,1)$, and the map $\bfa \mapsto (0, (\bfa,\bfzhat))$, we obtain a continuous map $\CG : \big\{ \bfa \in \RRn : \|\bfa\| \le \hat{R} \big\} \rightarrow \big\{ \bfa \in \RRn : \|\bfa\|=\hat{R} \big\}$ as follows:
\begin{equation*}
\CG(\bfa) := \CP \circ \CR \big( (0, (\bfa,\bfzhat) ), 1 \big) \,.
\end{equation*}
Since $\CR$ is a strong deformation retraction of $\CW^0$ onto $\CW^-$, $\CR ( (\tau, \bfx), 1 ) = (\tau, \bfx)$ for all $(\tau, \bfx) \in \CW^-$. Note that for all $\bfa$ with $\|\bfa\|=\hat{R}$, $(0, (\bfa,\bfzhat) ) \in \CW^-$ according to (\ref{CW-}). Thus, for all $\bfa$ with $\|\bfa\|=\hat{R}$,
\begin{align*}
\CG(\bfa) ={}& \CP \circ \CR \big( ( 0, (\bfa,\bfzhat) ), 1 \big) \\
={}& \CP ( 0, (\bfa,\bfzhat) ) \\
={}& \bfa \,.
\end{align*}
The existence of such a $\CG$ contradicts the fact that there is no retraction that maps a closed $n$-ball onto its boundary (i.e., an $(n\!-\!1)$-sphere).
\end{proof}

Based on Lemma \ref{lem:C1:subspace} and Lemma \ref{lem:C1:exactly-1}, we can associate each $x \in \Gamma$ with a unique linear operator $H(x) : \RRm \rightarrow \RRn$ such that $T(x)$ is the graph of $H(x)$, i.e.,
\begin{equation} \label{LSTx:graph}
T(x) = \big\{ (H(x)\,\bfz, \bfz) : \bfz \in \RRm \big\} \,.
\end{equation}
In addition, it follows from (\ref{CW-}) that for each fixed $x \in \Gamma$, 
\begin{equation*}
\CL\big( (H(x)\,\bfz, \bfz), \bfzero \big) < 0 \text{ for all $\bfz \in \RRm$ with } \|\bfz\| = 1 \,.
\end{equation*}
Then the compactness of the set $\big\{ \bfz \in \RRm : \|\bfz\| = 1 \big\}$ implies that for any $x \in \Gamma$,
\begin{equation} \label{Hx:norm}
\|H(x)\| = \max_{\|\bfz\| = 1} \| H(x)\,\bfz \| < 1 \,.
\end{equation}


\subsubsection{The map\/ $x \mapsto H(x)$ is continuous.} \label{H:continuous}

We first prove the following lemma.

\begin{lem}\label{lem:C1:point-continuous}
For any fixed\/ $\bfz \in \RRm$, the map from\/ $\Gamma$ into\/ $\RRn$: $x \mapsto H(x)\,\bfz$ is continuous.
\end{lem}

\begin{proof}
Take an arbitrary $\xhat \in \Gamma$ and any $\bfz \in \RRm$. We prove that $\lim_{x \rightarrow \xhat} H(x)\,\bfz = H(\xhat)\,\bfz$ by contradiction. Suppose that there exists a constant $\varepsilon > 0$ and a sequence $\big\{ x_1,x_2,x_3, ... \big\} \subset \Gamma$ such that $x_i \rightarrow \xhat$ and at the same time
\begin{equation} \label{Hxi-Hxhat}
\| H(x_i)\,\bfz - H(\xhat)\,\bfz \| > \varepsilon 
\end{equation}
for all $i=1,2,3,...$. Since $\|H(x_i)\,\bfz \| < \| \bfz \|$ for $i=1,2,3,...$, we can take a convergent subsequence $\big\{ H(x_{i_1})\,\bfz, H(x_{i_2})\,\bfz, H(x_{i_3})\,\bfz, ... \big\}$ and denote its limit by $\hat{\bfa}$, i.e., $H(x_{i_k})\,\bfz \rightarrow \hat{\bfa}$. Then (\ref{Hxi-Hxhat}) implies that $\| \hat{\bfa} - H(\xhat)\,\bfz \| \ge \varepsilon$. Thus $\hat{\bfx} := (\hat{\bfa}, \bfz) \not\in T( \hat{x} )$, and there exists a $t_0 \ge 0$ such that $Q(t_0,\xhat)\,\hat{\bfx} \in \interior(\CC( \bfzero ))$. Let $q$ be the radius of a closed ball that is centered at $Q(t_0,\xhat)\,\hat{\bfx}$ and contained in $\interior(\CC( \bfzero ))$.

Let $\bfx_{i_k} := (H(x_{i_k})\,\bfz, \bfz) \in T(x_{i_k})$ for $k=1,2,3,...$. Consider the convergent sequence $\big\{ (x_{i_1}, \bfx_{i_1}), (x_{i_2}, \bfx_{i_2}), (x_{i_3}, \bfx_{i_3}), ... \big\} \subset \Gamma \times X$ with limit $(\xhat, \hat\bfx)$. Since $Q(t_0,x)\,\bfx$ depends on $(x,\bfx)$ continuously for all $(x,\bfx) \in \Gamma \times X$, we can choose an $i_{k^*}$ large enough such that
\begin{equation*}
\| Q(t_0,x_{i_{k^*}})\,\bfx_{i_{k^*}} - Q(t_0,\xhat)\,\hat{\bfx} \| \le q \,.
\end{equation*}
Then $Q(t_0,x_{i_{k^*}})\,\bfx_{i_{k^*}} \in \interior(\CC( \bfzero ))$. However, this is impossible since $\bfx_{i_{k^*}} \in T( x_{i_{k^*}} )$ by our construction.
\end{proof}

Since $H(x)$ is a linear operator from $\RRm$ to $\RRn$ for each $x \in \Gamma$, Lemma \ref{lem:C1:point-continuous} implies the continuity of the map $x \mapsto H(x)$ for all $x \in \Gamma$.



\subsubsection{$H(h(z),z)$ is the derivative of\/ $h$ at\/ $z$.}\label{H:derivative}

The final step of establishing the $C^1$ smoothness of $h$ is to prove the next lemma.

\begin{lem}\label{lem:C1:h-differentiable}
$h : K \rightarrow \RRn$ is differentiable at all\/ $z \in K$. In addition, $Dh(z) = H(h(z),z)$ for all\/ $z \in K$.
\end{lem}

\begin{proof}
We will show that $H(x)$ satisfies the definition of the derivative of $h$, i.e., for any $z \in K$ and correspondingly $x = ( h(z), z ) \in \Gamma$,
\begin{equation*}
\lim_{\|\bfz\| \rightarrow 0} \tfrac{1}{\| \bfz \|} {\| h(z + \bfz) - h(z) - H(x)\,\bfz \|} = 0 \,.
\end{equation*}
We prove this by contradiction. Assume that for some $\zhat \in K$ and correspondingly $\xhat = ( h(\zhat), \zhat ) \in \Gamma$, there exist $\big\{ \bfz_1, \bfz_2, \bfz_3, ... \big\}$, a sequence of {\it non-zero} vectors in $\RRm$ with $\|\bfz_i\| \rightarrow 0$, and a constant $\varepsilon > 0$ such that $\zhat + \bfz_i \in K$ and
\begin{equation} \label{bigger-than-epsilon}
\tfrac{1}{\| \bfz_i \|} {\| h(\zhat + \bfz_i) - h(\zhat) - H(\xhat)\,\bfz_i \|} > \varepsilon
\end{equation}
for all $i=1,2,3,...$. Let $\bfa_i = ( h(\zhat + \bfz_i) - h(\zhat) ) / \| \bfz_i \|$ for $i=1,2,3,...$. By Lemma \ref{lem:h:Lip}, we have $\|\bfa_i\| < 1$ for any $i$. Let $\bfx_i = (\bfa_i, \bfz_i / \| \bfz_i \| )$ for $i=1,2,3,...$. Then $\big\{ \bfx_1, \bfx_2, \bfx_3, ... \big\}$ is a bounded sequence in $X$. Take a convergent subsequence $\big\{ \bfx_{i_1}, \bfx_{i_2}, \bfx_{i_3}, ... \big\}$, and denote its limit by $\hat{\bfx} = (\hat{\bfa}, \hat{\bfz})$. Note that $\bfa_{i_k} \! \rightarrow \hat{\bfa}$ and $\bfz_{i_k}/\| \bfz_{i_k} \| \rightarrow \hat{\bfz}$. Then (\ref{bigger-than-epsilon}) implies that $\| \hat{\bfa} - H(\xhat)\,\hat{\bfz} \| \ge \varepsilon$. Thus $\hat{\bfx} \not\in T( \xhat )$, and there exists a $t_0 \ge 0$ such that $Q(t_0,\xhat)\,\hat{\bfx} \in \interior(\CC( \bfzero ))$. Let $q$ be the radius of a closed ball that is centered at $Q(t_0,\xhat)\,\hat{\bfx}$ and contained in $\interior(\CC( \bfzero ))$.

Furthermore, by the constructions of $\bfa_i$ and $\bfx_i$, we have that for any $i_k$,
\begin{align*}
\xhat + \|\bfz_{i_k}\| \bfx_{i_k} ={}& ( h(\zhat) + \|\bfz_{i_k}\| \bfa_{i_k} ,\, \zhat + \bfz_{i_k} ) \\
={}& ( h(\zhat + \bfz_{i_k}) ,\, \zhat + \bfz_{i_k} ) \in \Gamma \,.
\end{align*}
Then $\Phi(t_0, \xhat) \in \Gamma$ and $\Phi( t_0, \xhat + \|\bfz_{i_k}\| \bfx_{i_k} ) \in \Gamma$ according to the positive invariance of $\Gamma$. Note that $\|\bfz_{i_k}\| \bfx_{i_k} \neq \bfzero$ for any $i_k$. Thus $\Phi( t_0, \xhat + \|\bfz_{i_k}\| \bfx_{i_k} ) \not\in \CC( \Phi(t_0, \xhat) )$ (see the proof of Lemma \ref{lem:h:Lip}) for any $i_k$.

On the other hand, note that $Q(t_0,\xhat)$ is just the derivative of $\Phi(t_0,x)$ with respect to $x$ evaluated at $x = \xhat$. Recall that by our construction, $\big\{ \bfx_{i_1}, \bfx_{i_2}, \bfx_{i_3}, ... \big\}$ is bounded and $\|\bfz_{i_k}\| \rightarrow 0$. Then there exists an $N$ such that for any $i_k \ge N$,
\begin{equation} \label{Q102}
\begin{split}
\tfrac{q}{3} \ge{}& \tfrac{1}{\|\bfz_{i_k}\|} \big\| \Phi( t_0, \xhat + \|\bfz_{i_k}\| \bfx_{i_k} ) - \Phi(t_0, \xhat) - Q(t_0,\xhat) (\|\bfz_{i_k}\| \bfx_{i_k}) \big\| \\
={}& \Big\| \tfrac{1}{\|\bfz_{i_k}\|} \big( \Phi( t_0, \xhat + \|\bfz_{i_k}\| \bfx_{i_k} ) - \Phi(t_0, \xhat) \big) - Q(t_0,\xhat)\,\bfx_{i_k} \Big\|  \,.  
\end{split}
\end{equation}
In addition, we can choose an $i_{k^*} \ge N$ such that
\begin{equation} \label{Q1976}
\| Q(t_0,\xhat)\,\bfx_{i_{k^*}} - Q(t_0,\xhat)\,\hat{\bfx} \| \le \tfrac{q}{3} \,.
\end{equation}
Combining (\ref{Q102}) and (\ref{Q1976}), we obtain
\begin{equation*}
\Big\| \tfrac{1}{\|\bfz_{i_{k^*}}\|} \big( \Phi( t_0, \xhat + \|\bfz_{i_{k^*}}\| \bfx_{i_{k^*}} ) - \Phi(t_0, \xhat ) \big) - Q(t_0,\xhat)\,\hat{\bfx} \Big\| \le \tfrac{q}{3} + \tfrac{q}{3} < q \,,
\end{equation*}
which implies that 
\begin{equation*}
\tfrac{1}{\|\bfz_{i_{k^*}}\|} \big( \Phi( t_0, \xhat + \|\bfz_{i_{k^*}}\| \bfx_{i_{k^*}} ) - \Phi(t_0, \xhat) \big) \in \interior(\CC( \bfzero )) \,.
\end{equation*}
It follows that $\Phi( t_0, \xhat + \|\bfz_{i_{k^*}}\| \bfx_{i_{k^*}} ) \in \interior \big( \CC( \Phi(t_0, \xhat) ) \big)$. This constitutes a contradiction.
\end{proof}

Note that the map from $K$ onto $\Gamma$: $z \mapsto (h(z),z)$ is continuous. Thus $Dh(z) = H(h(z),z)$ is continuous with respect to $z$ on $K$. Finally, it follows from (\ref{Hx:norm}) that $\|Dh(z)\| < 1$ for all $z \in K$. This concludes the proof of the $C^1$ smoothness of $h$.

\section{Proof of Theorem \ref{thm2}}\label{smoothness:Cr}

We begin with the assumption that the positively invariant set $\Gamma$ is a $C^1$ manifold embedded in $\RRn \times \RRm$ as the graph of the $C^1$ function $h_0$ from an open, convex $K_0 \subseteq \RRm$ to $\RRn$. The proof of the $C^r$ smoothness of $h_0$ proceeds as follows. In Subsection \ref{dyn:Dh0}, we derive a system that describes the dynamics of $Dh_0$ along trajectories in the $C^1$ manifold $\Gamma$. In Subsection \ref{case:C2}, we show that this system can also be put into the form of (\ref{sys}). Then by applying Theorem \ref{thm1}, we prove that the set $\big\{ (Dh_0(z), z) : z \in K_0 \big\}$ can be smoothly embedded into an appropriate Euclidean space as a $C^1$ manifold, which implies the $C^2$ smoothness of $h_0$. Finally, in Subsection \ref{case:Cr}, we demonstrate how to apply this argument inductively to establish the $C^r$ smoothness of $h_0$.

\subsection{The Dynamics of $Dh_0$ along Trajectories in $\Gamma$} \label{dyn:Dh0}

We shall continue to use some of the notations that have been introduced in Subsection \ref{smoothness:C1}. However, from now on we represent all linear operators as matrices. In particular, for each $x \in \Gamma$, the $(n+m) \times (n+m)$ matrix function $t \mapsto Q(t,x)$ is just the {\it fundamental matrix solution} of (\ref{sys:vari}) with parameter $x$, satisfying the matrix differential equation
\begin{equation*}
\dot{Q}(t,x) = DF(\Phi(t,x))\, Q(t,x) 
\end{equation*}
with $Q(0,x) = I_{n+m}$, which is the $(n+m) \times (n+m)$ identity matrix.

Restricting the $z$-component of (\ref{sys}) to $\Gamma$, we obtain
\begin{equation*} 
\dot{z} = g(h_0(z),z) \,,
\end{equation*}
which generates a flow $\phi(t,z)$ on $K_0$. Due to the positive invariance of $\Gamma$, $\phi(t,z)$ is related to $\Phi(t,x)$ in the following obvious way:
\begin{gather*} 
\big( h_0(\phi(t,z)), \phi(t,z) \big) = \Phi \big( t, (h_0(z),z) \big) \text{ for any $z \in K_0$ and any $t \ge 0$.} 
\end{gather*}
Then the dynamics of $Dh_0$ along trajectories in $\Gamma$ are governed by a matrix Riccati differential equation as stated in the next lemma.

\begin{lem} \label{lem:MRDE}
Consider any system in the form of (\ref{sys}) with\/ $f$ and\/ $g$ at least\/ $C^1$. If\/ $h_0: K_0\subseteq \RRm \rightarrow \RRn$ is\/ $C^1$ and\/ $\Gamma = \big\{ (h_0(z), z) : z \in K_0 \big\}$ is positively invariant, then for each fixed\/ $z \in K_0$, the matrix function\/ $Dh_0( \phi(\,\cdot\, , z) ) : [0,\infty) \rightarrow \RRnbym$ (the space of\/ $n \times m$ real matrices) is a solution of the following matrix Riccati differential equation:
\begin{equation} \label{MRDE}
\dot \Vi = \Daf \, \Vi - \Vi \Dzg - \Vi \Dag \, \Vi + \Dzf \,,
\end{equation}
where\/ $\Vi \in \RRnbym$ and\/ $\Daf$, $\Dzf$, $\Dag$, and\/ $\Dzg$ are all evaluated along the trajectory\/ $\big( h_0(\phi(t,z)), \phi(t,z) \big)$ for\/ $t \ge 0$.
\end{lem}

\begin{proof}
Take an arbitrary $z \in K_0$ and correspondingly $x = (h_0(z),z) \in \Gamma$. First, we show that the matrix function $Dh_0( \phi(\,\cdot\, , z) ) : [0,\infty) \rightarrow \RRnbym$ is continuously differentiable. 

Note that for any $t \ge 0$, $Q(t,x)$ maps the tangent space of $\Gamma$ at $x$ to the tangent space of $\Gamma$ at $\Phi(t,x)$. It follows that for any $\bfz \in \RRm$, 
\begin{equation} \label{LSTx:inv:matrix}
Q(t,x) 
\begin{pmatrix} Dh_0(z) \\ I_{m} \end{pmatrix} \bfz =
\begin{pmatrix} Dh_0( \phi(t,z) ) \\ I_m \end{pmatrix}
\begin{pmatrix} 0_{m \times n} \!&\!\! I_m \end{pmatrix} 
Q(t,x)
\begin{pmatrix} Dh_0(z) \\ I_{m} \end{pmatrix} \bfz \,.
\end{equation}
Note that for any $\bfz \in \RRm$ with $\|\bfz\|>0$, the left hand side of (\ref{LSTx:inv:matrix}) represents a solution (of (\ref{sys:vari}) with the parameter $x = (h_0(z),z) \in \Gamma$) that cannot reach the origin for any finite $t \ge 0$. Thus it must be true that for any $t \ge 0$, the null space of the product of the second, third, and forth matrices on the right-hand side of (\ref{LSTx:inv:matrix}) is the trivial subspace of $\RRm$. Therefore, this product produces an invertible matrix for any $t \ge 0$. Solving (\ref{LSTx:inv:matrix}) for $Dh_0( \phi(t,z) )$, we obtain
$$ 
Dh_0( \phi(t,z) ) = 
\begin{pmatrix} I_n \!&\!\! 0_{n\times m} \end{pmatrix} 
Q(t,x) \begin{pmatrix} Dh_0(z) \\ I_{m} \end{pmatrix} \!\left(
\begin{pmatrix} 0_{m \times n} \!&\!\! I_m \end{pmatrix} Q(t,x)
\begin{pmatrix} Dh_0(z) \\ I_{m} \end{pmatrix} \right)^{-1} .
$$ 
In view of the differentiability of $Q(t,x)$ with respect to $t$, $Dh_0( \phi(t,z) )$ is continuously differentiable with respect to $t$ for all $t \ge 0$.

Next, we derive a differential equation for $Dh_0( \phi(t,z) )$. Choose any $t \ge 0$, and consider (\ref{sys:vari:a-z}) at $\big( Dh_0(\phi(t,z))\,\bfz, \bfz \big) \in T_{\Phi(t,x)}\Gamma$ for an arbitrary $\bfz \in \RRm$. Then we have
\begin{align*}
\big( \ddt Dh_0(\phi(t,z)) \big)\, \bfz + Dh_0(\phi(t,z))\, \dot{\bfz} \\
={}& \Daf(\Phi(t,x))\, Dh_0(\phi(t,z))\, \bfz + \Dzf(\Phi(t,x))\, \bfz \,, \\
\dot{\bfz} ={}& \Dag(\Phi(t,x))\, Dh_0(\phi(t,z))\, \bfz + \Dzg(\Phi(t,x))\, \bfz \,.
\end{align*}
Combining these two equations and reorganizing terms, we obtain that $Dh_0(\phi(t,z))$ satisfies the matrix Riccati differential equation (\ref{MRDE}) with the coefficient matrices $\Daf$, $\Dzf$, $\Dag$, and $\Dzg$ all evaluated along the trajectory $\Phi(t,x) = \big( h_0(\phi(t,z))$, $\phi(t,z) \big)$ for $t \ge 0$. 
\end{proof}

\subsection{The $C^2$ Smoothness of $h_0$} \label{case:C2}

For each $z \in K_0$, (\ref{MRDE}) defines a time varying system in $\RRnbym$ along the trajectory $\big( h_0(\phi(t,z)), \phi(t,z) \big)$ for $t \ge 0$. 
Thus, we can couple (\ref{MRDE}) with the underlying equation $\dot{z} = g(h_0(z),z)$, which generates the flow $\phi(t,z)$, to form a system in $\RRnbym \times K_0$ as follows:
\begin{subequations} \label{sys1M}
\begin{align}
\dot \Vi ={}& \Daf(h_0(z),z)\, \Vi - \Vi \Dzg(h_0(z),z) - \Vi \Dag(h_0(z),z)\, \Vi + \Dzf(h_0(z),z) \,, \label{sys1M:V1} \\
\dot{z} ={}& g(h_0(z),z) \,, \label{sys1M:z}
\end{align}
\end{subequations}
where $\Vi \in \RRnbym$ and $z \in K_0$. 

In the subsequent analysis, it is more convenient to consider the ``vectorization'' of the matrix differential equation (\ref{sys1M:V1}). Recall that the {\it vectorization} of a matrix is a transformation that converts the matrix into a column vector by stacking the columns of the matrix. Let $\vectnm : \RRnbym \rightarrow \RRnm$ be the vectorization of $n \times m$ matrices, and let $\vectnm^{-1} : \RRnm \rightarrow \RRnbym$ be the inverse of $\vectnm$. Notice that for any $A \in \RRnbyn$, $B \in \RRmbym$, and $P \in \RRnbym$, 
\begin{equation*}
\vectnm\!(APB) = (B^T \otimes A) \vectnm\!(P) \,,
\end{equation*}
where ``$\otimes$'' denotes the Kronecker (or tensor) product of two matrices. Then (\ref{sys1M:V1}) can be converted into
\begin{equation*} 
\begin{split} 
\dot\vi ={}& \Big( I_m \otimes \Daf(h_0(z),z) - \big( \Dzg(h_0(z),z) \big)^{T} \otimes I_n \Big)\, \vi \\
& - \Big( I_m \otimes \big( \vectnm^{-1}\!(\vi)\, \Dag(h_0(z),z) \big) \Big)\, \vi + \vectnm\! \big( \Dzf(h_0(z),z) \big) \,,
\end{split}
\end{equation*}
where $\vi := \vectnm\!(\Vi) \in \RRnm$.

Furthermore, we rescale $z$ to $\zetai := \sigma_1^{-1} z$ with the scaling factor $\sigma_1 > 0$ to be determined later. Let $K_1 := \big\{ \zetai \in \RRm : \sigma_1 \zetai \in K_0 \big\}$. Then we transform (\ref{sys1M}) into a system defined on $\RRnm \times K_1$ as follows:
\begin{equation} \label{sys1}
\begin{split}
\dot\vi ={}& f_1(\vi,\zetai) \,, \\
\dot\zetai ={}& \gamma_1(\zetai) \,, 
\end{split}
\end{equation}
where $\vi \in \RRnm$, $\zetai \in K_1$, and the functions $\gamma_1 : K_1 \rightarrow \RRm$ and $f_1 : \RRnm \times K_1 \rightarrow \RRnm$ are defined as follows:
\begin{gather}
\gamma_1(\zetai) := \tfrac{1}{\sigma_1} g( h_0(\sigma_1\zetai), \sigma_1\zetai ) \,, \label{gamma1:def} \\
\begin{split} \label{f1:def}
f_1(\vi,\zetai) :={}& \Big( I_m \otimes \Daf( h_0(\sigma_1\zetai), \sigma_1\zetai ) - \big( \Dzg( h_0(\sigma_1\zetai), \sigma_1\zetai ) \big)^{T} \otimes I_n \Big)\, \vi \\
& - \Big( I_m \otimes \big( \vectnm^{-1}\!(\vi)\, \Dag( h_0(\sigma_1\zetai), \sigma_1\zetai ) \big) \Big)\, \vi \\
& + \vectnm\! \big( \Dzf( h_0(\sigma_1\zetai), \sigma_1\zetai ) \big) \,.
\end{split}
\end{gather}

Choose $J_1$, a bounded, open subset of $\RRnm$, as follows:
\begin{equation*} 
J_1 := \big\{ \vi \in \RRnm : \| \vectnm^{-1}\!(\vi) \| < \eta \big\} \,,
\end{equation*}
where, as postulated in Hypothesis \ref{hypo:K0}, $\eta$ is a positive constant such that $\|Dh_0(z)\| < \eta$ for all $z \in K_0$. Obviously $J_1 \times K_1$ satisfies Hypothesis \ref{hypo:U}. In what follows, we will verify that the restriction of (\ref{sys1}) to $J_1 \times K_1$ satisfies Hypothesis \ref{hypo:C1}.

\begin{lem} \label{lem:f1:zeta1}
$f_1$ is\/ $C^1$ on\/ $J_1 \times K_1$, and\/ $\gamma_1$ is\/ $C^1$ on\/ $K_1$. Furthermore, there exists a constant\/ $L_1 > 0$ independent of the choice of\/ $\sigma_1$ such that\/ $\| \Dzetaifi(\vi,\zetai) \| \le \sigma_1 L_1$ for any\/ $(\vi,\zetai) \in J_1 \times K_1$.
\end{lem}

\begin{proof}
By Hypothesis \ref{hypo:Cr0}, $f$ and $g$ are $C^r$ ($r \ge 2$) with their first to $r$-th derivatives all bounded on $U_0$. In addition, $h_0$ is $C^1$ with $\|Dh_0(z)\| < \eta$ for all $z \in K_0$. Then by inspecting (\ref{gamma1:def}) and (\ref{f1:def}), we obtain the abovementioned properties of $f_1$ and $\gamma_1$.
\end{proof}

Define a continuous function $\alpha_1 : K_1 \rightarrow \RR$ as follows:
\begin{equation} \label{alpha1:def}
\alpha_1(\zetai) := \alpha( h_0(\sigma_1\zetai), \sigma_1\zetai ) - \ell( h_0(\sigma_1\zetai), \sigma_1\zetai ) - 2 \eta \big\| \Dag( h_0(\sigma_1\zetai), \sigma_1\zetai ) \big\| \,.
\end{equation}
By (\ref{ineq:Cr0}) of Hypothesis \ref{hypo:Cr0}, we have that $\alpha_1(\zetai) > 0$ for all $\zetai \in K_1$.

\begin{lem} \label{lem:f1:v1}
For any\/ $(\vi_0,\zetai_0) \in J_1 \times K_1$ and any\/ $\vi \in \RRnm$, 
\begin{equation} \label{ineq:f1:v1}
\langle \vi, \Dvifi(\vi_0,\zetai_0)\, \vi \rangle \ge \alpha_1(\zetai_0) \|\vi\|^2 \,.
\end{equation}
\end{lem}

\begin{proof}
Take an arbitrary $(\vi_0,\zetai_0) \in J_1 \times K_1$. It follows from (\ref{f1:def}) that 
\begin{equation*} 
\begin{split} 
\Dvifi(\vi_0,\zetai_0) ={}& I_m \otimes \Daf( h_0(\sigma_1\zetai_0), \sigma_1\zetai_0 ) - \big( \Dzg( h_0(\sigma_1\zetai_0), \sigma_1\zetai_0 ) \big)^{T} \otimes I_n \\
& - I_m \otimes \big( \vectnm^{-1}\!(\vi_0)\, \Dag( h_0(\sigma_1\zetai_0), \sigma_1\zetai_0 ) \big) \\
& - \big( \Dag( h_0(\sigma_1\zetai_0), \sigma_1\zetai_0 )\, \vectnm^{-1}\!(\vi_0) \big)^T \otimes I_n \,.
\end{split}
\end{equation*}

Taking account of (\ref{ineq:a-a:Cr0}), it is straightforward to verify that for all $\vi \in \RRnm$,
\begin{equation*}
\langle \vi, ( I_m \otimes \Daf )\, \vi \rangle \ge \alpha \|\vi\|^2 \,,
\end{equation*}
where both $\Daf$ and $\alpha$ are evaluated at $( h_0(\sigma_1\zetai_0), \sigma_1\zetai_0 )$.

Next, take any $\vi \in \RRnm$ and correspondingly $\Vi = \vectnm^{-1}\!(\vi)$. We have that
\begin{align*}
\big\langle \vi, ( ( \Dzg )^{T} \otimes I_n )\, \vi \big\rangle ={}& \langle \vectnm\!(\Vi), \vectnm\!( \Vi \Dzg ) \rangle \\
={}& \big\langle \vectmn\!({\Vi}^{T}), \vectmn\!( (\Dzg)^{T} {\Vi}^{T} ) \big\rangle \\
={}& \big\langle \vectmn\!({\Vi}^{T}), ( I_n \otimes (\Dzg)^{T} ) \vectmn\!( {\Vi}^{T} ) \big\rangle \\
={}& \big\langle \vectmn\!({\Vi}^{T}), ( I_n \otimes (\Dzg)^{T} )^{T} \vectmn\!( {\Vi}^{T} ) \big\rangle \\
={}& \big\langle \vectmn\!({\Vi}^{T}), ( I_n \otimes \Dzg ) \vectmn\!( {\Vi}^{T} ) \big\rangle \\
\le{}& \ell \| \vectmn\!({\Vi}^{T}) \|^2 = \ell \| \vi \|^2 \,,
\end{align*}
where both $\Dzg$ and $\ell$ are evaluated at $( h_0(\sigma_1\zetai_0), \sigma_1\zetai_0 )$ throughout and the inequality follows from (\ref{ineq:z-z:Cr0}).

Similarly, we obtain that for any $\vi \in \RRnm$,
\begin{align*}
\Big| \big\langle \vi, \big( I_m \otimes ( \vectnm^{-1}\!(\vi_0)\, \Dag ) \big)\, \vi \big\rangle \Big| \le{}& \| \vectnm^{-1}\!(\vi_0)\, \Dag \| \|\vi\|^2 < \eta \| \Dag \| \|\vi\|^2 \,, \\
\Big| \big\langle \vi, \big( ( \Dag\, \vectnm^{-1}\!(\vi_0) )^T \otimes I_n \big)\, \vi \big\rangle \Big| \le{}& \| \Dag\, \vectnm^{-1}\!(\vi_0) \| \|\vi\|^2 < \eta \| \Dag \| \|\vi\|^2 \,,
\end{align*}
where $\Dag$ is evaluated at $( h_0(\sigma_1\zetai_0), \sigma_1\zetai_0 )$ throughout and $\| \vectnm^{-1}\!(\vi_0) \| < \eta$ by our choice of $J_1$.

Combining all the estimates above, we attain the inequality (\ref{ineq:f1:v1}).
\end{proof}

Define a nonnegative, continuous function $\ell_1 : K_1 \rightarrow \RR$ as follows:
\begin{equation} \label{ell1:def}
\ell_1(\zetai) := \eta \big\| \Dag( h_0(\sigma_1\zetai), \sigma_1\zetai ) \big\| + \ell( h_0(\sigma_1\zetai), \sigma_1\zetai ) \,.
\end{equation}
Then in view of (\ref{gamma1:def}), the next lemma is a trivial consequence of (\ref{ineq:z-z:Cr0}).

\begin{lem} \label{lem:gamma1}
For any\/ $\zetai_0 \in K_1$ and any\/ $\zetai \in \RRm$,
\begin{equation} \label{ineq:gamma1:zeta1}
\langle \zetai, D \gamma_1(\zetai_0)\, \zetai \rangle \le \ell_1(\zetai_0) \|\zetai\|^2 \,.
\end{equation}
\end{lem}

In addition, we have the estimate described in the following lemma.

\begin{lem} \label{lem:hyper:C2}
There exists a sufficiently small\/ $\sigma_1$ such that for any\/ $\zetai \in K_1$, 
\begin{equation} \label{ineq:hyper:sys1}
\alpha_1(\zetai) \ge \ell_1(\zetai) + \sigma_1 L_1 + \tfrac{1}{2} c_r \,.
\end{equation}
\end{lem}

\begin{proof}
It follows from (\ref{ineq:Cr0}) of Hypothesis \ref{hypo:Cr0} that for any $z \in K_0$,
\begin{equation} \label{ineq:hyper:Cr:2}
\alpha( h_0(z), z ) \ge 2 \ell( h_0(z), z ) + 3 \eta \big\| \Dag( h_0(z), z ) \big\| + c_r \,.
\end{equation}
In addition, since $L_1 > 0$ is independent of the choice of $\sigma_1$ (see Lemma \ref{lem:f1:zeta1}), we can choose $\sigma_1$ sufficiently small so that $\sigma_1 L_1 \le \tfrac{1}{2} c_r$. Then by combining this with (\ref{ineq:hyper:Cr:2}), we obtain the inequality (\ref{ineq:hyper:sys1}) after replacing $z$ with $\sigma_1\zetai$ and rearranging terms according to (\ref{alpha1:def}) and (\ref{ell1:def}).
\end{proof}

Note that the $\zetai$-component of (\ref{sys1}) is independent of $\vi$. Thus Lemmas \ref{lem:f1:zeta1}--\ref{lem:hyper:C2} verify that the restriction of (\ref{sys1}) to $J_1 \times K_1$ satisfies Hypotheses \ref{hypo:C1}.

Recall that $\|Dh_0(z)\| < \eta$ for all $z \in K_0$. Define $h_1 : K_1 \rightarrow J_1$ as follows:
\begin{equation} \label{h1:def}
h_1(\zetai) := \vectnm\! ( Dh_0(\sigma_1\zetai) ) \,.
\end{equation}
Now consider the set 
\begin{equation*}
\Gamma_1 := \big\{ (h_1(\zetai), \zetai) : \zetai \in K_1 \big\} \subset J_1 \times K_1 \,.
\end{equation*} 
Clearly, $\Gamma_1$ is positively invariant under the flow of (\ref{sys1}). Then, by applying Theorem \ref{thm1}, we obtain the following result.

\begin{lem} \label{lem:h1}
$\Gamma_1$ is a\/ $C^1$ manifold. In particular, $h_1 : K_1 \rightarrow J_1$ is\/ $C^1$ with\/ $\|Dh_1(\zetai)\| < 1$ for all\/ $\zetai \in K_1$.
\end{lem}

Note that $Dh_0(z) = \vectnm^{-1}\!( h_1(\sigma_1^{-1} z) )$ for any $z \in K_0$. Therefore, we have proven that $h_0 : K_0 \rightarrow \RRn$ is $C^2$ with its first and second derivatives bounded on $K_0$. Furthermore, 
\begin{equation*}
D^2 h_0(z)(\bfz_1,\bfz_2) = \sigma_1^{-1} \vectnm^{-1}\!( Dh_1(\sigma_1^{-1} z)\, \bfz_1 )\, \bfz_2
\end{equation*}
for any $z \in K_0$ and any $\bfz_1$, $\bfz_2 \in \RRm$.

\subsection{The $C^r$ Smoothness of $h_0$} \label{case:Cr}

In this subsection, we apply the argument used in the previous subsection inductively to establish the $C^r$ smoothness of $h_0$ for the case with $r \ge 3$. 

We illustrate how the argument works again in the proof of the $C^3$ smoothness of $h_0$. First, let $\phi_1(t,\zetai)$ be the flow on $K_1$ generated by $\dot\zetai = \gamma_1(\zetai)$. Then by Lemma \ref{lem:MRDE}, we have that for each fixed $\zetai \in K_1$, the matrix function $Dh_1( \phi_1(\,\cdot\, , \zetai) ) : [0,\infty) \rightarrow \RRnmbym$ is a solution of the matrix Riccati differential equation:
\begin{equation} \label{MRDE:2}
\dot \Vii = \Dvifi \, \Vii - \Vii D\gamma_1 + \Dzetaifi \,,
\end{equation}
where $\Vii \in \RRnmbym$ and $\Dvifi$, $\Dzetaifi$, and $D\gamma_1$ are all evaluated along the trajectory $\big( h_1( \phi_1(t,\zetai) ), \phi_1(t,\zetai) \big)$ for $t \ge 0$. Compared to (\ref{MRDE}), the right-hand side of (\ref{MRDE:2}) does not include the term quadratic in $\Vii$ since $D_{\vi}\mspace{-2mu} \gamma_1 = 0$. Next, we couple (\ref{MRDE:2}) with the underlying equation $\dot\zetai = \gamma_1(\zetai)$ to form a system in $\RRnmbym \times K_1$ as follows:
\begin{subequations} \label{sys2M}
\begin{align}
\dot \Vii ={}& \Dvifi(h_1(\zetai),\zetai) \, \Vii - \Vii D\gamma_1(\zetai) + \Dzetaifi(h_1(\zetai),\zetai) \,, \label{sys2M:V2} \\
\dot\zetai ={}& \gamma_1(\zetai) \,. 
\end{align}
\end{subequations}

We shall again consider the vectorization of the matrix differential equation (\ref{sys2M:V2}). Let $\vectnmm : \RRnmbym \rightarrow \RRnmii$ be the vectorization of $nm \times m$ matrices, and let $\vectnmm^{-1}$ be the inverse of $\vectnmm$. Furthermore, we rescale $\zetai$ to $\zetaii := \sigma_2^{-1} \zetai$ with the scaling factor $\sigma_2 > 0$ to be determined. Take $K_2 := \big\{ \zetaii \in \RRm : \sigma_2 \zetaii \in K_1 \big\}$. Then we transform (\ref{sys2M}) into a system defined on $\RRnmii \times K_2$ as follows:
\begin{equation} \label{sys2}
\begin{split}
\dot\vii ={}& f_2(\vii,\zetaii) \,, \\
\dot\zetaii ={}& \gamma_2(\zetaii) \,, 
\end{split}
\end{equation}
where $\vii := \vectnmm\!(\Vii) \in \RRnmii$, $\zetaii \in K_2$, and the functions $\gamma_2 : K_2 \rightarrow \RRm$ and $f_2 : \RRnmii \times K_2 \rightarrow \RRnmii$ are defined as follows:
\begin{gather}
\gamma_2(\zetaii) := \tfrac{1}{\sigma_2} \gamma_1( \sigma_2\zetaii ) \,, \label{gamma2:def} \\
\begin{split} \label{f2:def}
f_2(\vii,\zetaii) :={}& \big( I_m \otimes \Dvifi( h_1(\sigma_2\zetaii), \sigma_2\zetaii ) - ( D\gamma_1( \sigma_2\zetaii ) )^{T} \otimes I_{nm} \big)\, \vii \\
& + \vectnmm\! \big( \Dzetaifi( h_1(\sigma_2\zetaii), \sigma_2\zetaii ) \big) \,.
\end{split}
\end{gather}

Choose $J_2$, a bounded, open subset of $\RRnmii$, as follows:
\begin{equation*} 
J_2 := \big\{ \vii \in \RRnmii : \| \vectnmm^{-1}\!(\vii) \| < 1 \big\} \,.
\end{equation*}
Obviously $J_2 \times K_2$ satisfies Hypothesis \ref{hypo:U}. In addition, the restriction of (\ref{sys2}) to $J_2 \times K_2$ again satisfies Hypothesis \ref{hypo:C1} as verified below.

\begin{lem} \label{lem:f2:zeta2}
$f_2$ is\/ $C^1$ on\/ $J_2 \times K_2$, and\/ $\gamma_2$ is\/ $C^2$ on\/ $K_2$. Furthermore, there exists a constant\/ $L_2 > 0$ independent of the choice of\/ $\sigma_2$ such that\/ $\| \Dzetaiifii(\vii,\zetaii) \| \le \sigma_2 L_2$ for any\/ $(\vii,\zetaii) \in J_2 \times K_2$.
\end{lem}

\begin{proof}
So far, we have proven that 1) $h_0$ is $C^2$ with $Dh_0$ and $D^2h_0$ both bounded on $K_0$ and 2) $h_1$ defined by (\ref{h1:def}) is $C^1$ with $Dh_1$ bounded on $K_1$. Thus, by our assumption that $f$ and $g$ are $C^r$ ($r \ge 3$ in this case) with their first to $r$-th derivatives all bounded on $U_0$, now $f_1$ and $\gamma_1$ defined by (\ref{f1:def}) and (\ref{gamma1:def}) are both $C^2$ with their first and second derivatives bounded on $J_1 \times K_1$ and $K_1$, respectively. In turn, $f_2$ is $C^1$ with $\| \Dzetaiifii \| \le \sigma_2 L_2$ on $J_2 \times K_2$ for some $L_2 > 0$, which is independent of the choice of $\sigma_2$. In addition, $\gamma_2$ is $C^2$ since it has the same smoothness as $\gamma_1$.
\end{proof}

Define continuous functions $\alpha_2 : K_2 \rightarrow \RR$ and $\ell_2 : K_2 \rightarrow \RR$ as follows:
\begin{align} 
\alpha_2(\zetaii) :={}& \alpha_1( \sigma_2\zetaii ) - \ell_1( \sigma_2\zetaii ) \,, \label{alpha2:def} \\
\ell_2(\zetaii) :={}& \ell_1( \sigma_2\zetaii ) \,. \label{ell2:def}
\end{align}
Clearly, $\alpha_2$ is positive because of (\ref{ineq:hyper:sys1}), and $\ell_2$ is nonnegative as $\ell_1$ is. In addition, we have a set of estimates for $\alpha_2$ and $\ell_2$ similar to what we have for $\alpha_1$ and $\ell_1$ in Lemmas \ref{lem:f1:v1}--\ref{lem:hyper:C2}. 

\begin{lem} \label{lem:f2:v2}
For any\/ $(\vii_0,\zetaii_0) \in J_2 \times K_2$ and any\/ $\vii \in \RRnmii$, 
\begin{equation} \label{ineq:f2:v2}
\langle \vii, \Dviifii(\vii_0,\zetaii_0)\, \vii \rangle \ge \alpha_2(\zetaii_0) \|\vii\|^2 \,.
\end{equation}
\end{lem}

\begin{proof}
Consider an arbitrary $(\vii_0,\zetaii_0) \in J_2 \times K_2$. Notice that $\Dviifii$ is simply the coefficient matrix left-multiplied to $\vii$ in (\ref{f2:def}). By the same arguments as what we have used in the proof of Lemma \ref{lem:f1:v1} and incorporating the estimates (\ref{ineq:f1:v1}) and (\ref{ineq:gamma1:zeta1}), we have that for all $\vii \in \RRnmii$,
\begin{align*}
\langle \vii, ( I_m \otimes \Dvifi )\, \vii \rangle \ge{}& \alpha_1 \|\vii\|^2 \,, \\
\big\langle \vii, ( ( D\gamma_1 )^{T} \otimes I_{nm} )\, \vii \big\rangle \le{}& \ell_1 \|\vii\|^2 \,,
\end{align*}
where $\Dvifi$ is evaluated at $( h_1(\sigma_2\zetaii_0), \sigma_2\zetaii_0 )$ and $D\gamma_1$, $\alpha_1$, and $\ell_1$ are all evaluated at $\sigma_2\zetaii_0$. Then (\ref{ineq:f2:v2}) follows.\end{proof}

\begin{lem} \label{lem:gamma2}
For any\/ $\zetaii_0 \in K_2$ and any\/ $\zetaii \in \RRm$,
\begin{equation} \label{ineq:gamma2:zeta2}
\langle \zetaii, D \gamma_2(\zetaii_0)\, \zetaii \rangle \le \ell_2(\zetaii_0) \|\zetaii\|^2 \,.
\end{equation}
\end{lem}

\begin{lem} \label{lem:hyper:C3}
There exists a sufficiently small\/ $\sigma_2$ such that for any\/ $\zetaii \in K_2$, 
\begin{equation} \label{ineq:hyper:sys2}
\alpha_2(\zetaii) \ge \ell_2(\zetaii) + \sigma_2 L_2 + \tfrac{1}{2} c_r \,.
\end{equation}
\end{lem}

\begin{proof}
Since (\ref{ineq:Cr0}) holds for some $r \ge 3$, we have that for any $z \in K_0$,
\begin{equation} \label{ineq:hyper:Cr:3}
\alpha( h_0(z), z ) \ge 3 \ell( h_0(z), z ) + 4 \eta \big\| \Dag( h_0(z), z ) \big\| + c_r \,.
\end{equation}
In addition, since $L_2 > 0$ is independent of the choice of $\sigma_2$ (see Lemma \ref{lem:f2:zeta2}), we can choose $\sigma_2$ sufficiently small so that $\sigma_2 L_2 \le \tfrac{1}{2} c_r$. Combine this with (\ref{ineq:hyper:Cr:3}). Recall (\ref{alpha1:def}), (\ref{ell1:def}), (\ref{alpha2:def}), and (\ref{ell2:def}), and rearrange terms accordingly. Then we obtain (\ref{ineq:hyper:sys2}) once replacing $z$ with $\sigma_1\sigma_2\zetaii$.
\end{proof}

Again, since the $\zetaii$-component of (\ref{sys2}) is independent of $\vii$, Lemmas \ref{lem:f2:zeta2}--\ref{lem:hyper:C3} verify that the restriction of (\ref{sys2}) to $J_2 \times K_2$ satisfies Hypotheses \ref{hypo:C1}.

Recall that $\|Dh_1(\zetai)\| < 1$ for all $\zetai \in K_1$ (see Lemma \ref{lem:h1}). Define $h_2 : K_2 \rightarrow J_2$ as follows:
\begin{equation} \label{h2:def}
h_2(\zetaii) := \vectnmm\! ( Dh_1(\sigma_2\zetaii) ) \,.
\end{equation}
Now consider the set 
\begin{equation*}
\Gamma_2 := \big\{ (h_2(\zetaii), \zetaii) : \zetaii \in K_2 \big\} \subset J_2 \times K_2 \,.
\end{equation*} 
Clearly, $\Gamma_2$ is positively invariant under the flow of (\ref{sys2}). Then, by applying Theorem \ref{thm1} again, we obtain the $C^1$ smoothness of $h_2$.

\begin{lem} \label{lem:h2}
$\Gamma_2$ is a\/ $C^1$ manifold. In particular, $h_2 : K_2 \rightarrow J_2$ is\/ $C^1$ with\/ $\|Dh_2(\zetaii)\| < 1$ for all\/ $\zetaii \in K_2$.
\end{lem}

Since $Dh_1(\zetai) = \vectnmm^{-1}\!( h_2(\sigma_2^{-1} \zetai) )$ for any $\zetai \in K_1$, we have proven that $h_1 : K_1 \rightarrow J_1$ is $C^2$ with its first and second derivatives bounded on $K_1$. Therefore, $h_0 : K_0 \rightarrow \RRn$ is now $C^3$ with its first to third derivatives bounded on $K_0$. Furthermore, 
$$ 
D^3 h_0(z)(\bfz_1,\bfz_2,\bfz_3) = \sigma_2^{-1} \sigma_1^{-2} \vectnm^{-1}\!\big( \vectnmm^{-1}\!( Dh_2(\sigma_2^{-1}\sigma_1^{-1} z)\, \bfz_1 )\, \bfz_2 \big)\, \bfz_3
$$ 
for any $z \in K_0$ and any $\bfz_1$, $\bfz_2$, and $\bfz_3 \in \RRm$.

\subsubsection{Induction} 

Up to this point, the inductive scheme of the proof has become obvious. Consider $r \ge 4$. Suppose that we have proven the $C^k$ smoothness of $h_0$ for some $3 \le k \le r-1$ and are about to prove the $C^{k+1}$ smoothness of $h_0$. Then we have already established the following prerequisites.  
\begin{enumerate} 
\item 
   For $j = 2, ..., k$, we have defined $\zetaj$, $K_j$, and $h_j : K_j \rightarrow \RRnmj$ recursively as follows:
\begin{align}
\zetaj :={}& \sigma_j^{-1} \zeta^{\sss j-1} \,, \notag \\
K_j :={}& \big\{ \zetaj \in \RRm : \sigma_j \zetaj \in K_{j-1} \big\} \,, \notag \\
h_j (\zetaj) :={}& \vectnmjMm\! ( Dh_{j-1}(\sigma_j\zetaj) ) \,, \label{hj:hjM}
\end{align}
with each $\sigma_j$ chosen sufficiently small. In addition, we have already proven that for each $j = 1, ..., k-1$, $h_j$ is $C^{k-j}$ with $\| Dh_j (\zetaj) \| < 1$ for all $\zetaj \in K_j$ and all applicable higher derivatives bounded on $K_j$.
\item
   For $j = 2, ..., k$, we have defined $\gamma_j : K_j \rightarrow \RRm$ and $f_j : \RRnmj \times K_j \rightarrow \RRnmj$ recursively as follows:
\begin{gather*}
\gamma_j(\zetaj) := \tfrac{1}{\sigma_j} \gamma_{j-1}( \sigma_j\zetaj ) \,, \\
\begin{split} 
f_j(\vj,\zetaj) :={}& \big( I_m \otimes \Dvjmfjm( h_{j-1}(\sigma_j\zetaj), \sigma_j\zetaj ) - ( D\gamma_{j-1}( \sigma_j\zetaj ) )^{T} \otimes I_{nm^{j-1}} \big)\, \vj \\
& + \vectnmjMm\! \big( \Dzetajmfjm( h_{j-1}(\sigma_j\zetaj), \sigma_j\zetaj ) \big) \,,
\end{split}
\end{gather*}
where $\vj \in \RRnmj$ and $\zetaj \in K_j$.
\item
   For $j = 2, ..., k$, we have defined continuous functions $\alpha_j : K_j \rightarrow \RR$ and $\ell_j : K_j \rightarrow \RR$ recursively as follows:
\begin{align*} 
\alpha_j(\zetaj) :={}& \alpha_{j-1}( \sigma_j\zetaj ) - \ell_{j-1}( \sigma_j\zetaj ) \,, \\
\ell_j(\zetaj) :={}& \ell_{j-1}( \sigma_j\zetaj ) \,. 
\end{align*}
\end{enumerate}

Then we shall consider the system
\begin{equation} \label{sysk}
\begin{split}
\dot\vk ={}& f_k(\vk,\zetak) \,, \\
\dot\zetak ={}& \gamma_k(\zetak) \,, 
\end{split}
\end{equation}
which is defined on $\RRnmk \times K_k$. Next, choose $J_k$, a bounded, open subset of $\RRnmk$, as follows:
\begin{equation*}
J_k := \big\{ \vk \in \RRnmk : \| \vectnmkMm^{-1}\!(\vk) \| < 1 \big\} \,.
\end{equation*}
Obviously $J_k \times K_k$ satisfies Hypothesis \ref{hypo:U}. In addition, the restriction of (\ref{sysk}) to $J_k \times K_k$ again satisfies Hypothesis \ref{hypo:C1}, and the verification is now routine. In particular, following the same arguments as used in the proofs of Lemmas \ref{lem:f2:zeta2}--\ref{lem:hyper:C3}, it is straightforward to verify the following statements about $f_k$, $\gamma_k$, $\alpha_k$, and $\ell_k$.
\begin{enumerate} 
\item 
   $f_k$ is $C^1$ on $J_k \times K_k$, and $\gamma_k$ is $C^k$ on $K_k$. Furthermore, there exists a constant $L_k > 0$ independent of the choice of $\sigma_k$ such that $\| \Dzetakfk(\vk,\zetak) \| \le \sigma_k L_k$ for any $(\vk,\zetak) \in J_k \times K_k$.
\item
   For any $(\vk_0,\zetak_0) \in J_k \times K_k$ and any $\vk \in \RRnmk$, 
\begin{equation*} 
\langle \vk, \Dvkfk(\vk_0,\zetak_0)\, \vk \rangle \ge \alpha_k(\zetak_0) \|\vk\|^2 \,.
\end{equation*}
\item
   For any $\zetak_0 \in K_k$ and any $\zetak \in \RRm$,
\begin{equation*}
\langle \zetak, D \gamma_k(\zetak_0)\, \zetak \rangle \le \ell_k(\zetak_0) \|\zetak\|^2 \,.
\end{equation*}
\item
   There exists a sufficiently small $\sigma_k$ such that for any $\zetak \in K_k$, 
\begin{equation} \label{ineq:hyper:sysk}
\alpha_k(\zetak) \ge \ell_k(\zetak) + \sigma_k L_k + \tfrac{1}{2} c_r \,.
\end{equation}
\end{enumerate}

Now we consider the set 
\begin{equation*}
\Gamma_k := \big\{ (h_k(\zetak), \zetak) : \zetak \in K_k \big\} \subset J_k \times K_k \,,
\end{equation*} 
which, by our construction, is positively invariant under the flow of (\ref{sysk}). By Theorem \ref{thm1}, $\Gamma_k$ is a $C^1$ manifold. In particular, $h_k : K_k \rightarrow J_k$ is $C^1$ with $\|Dh_k(\zetak)\| < 1$ for all $\zetak \in K_k$. Then, in view of the inverse of the recurrence relation (\ref{hj:hjM}), we have proven that for each $j = 1, ..., k$, $h_j$ is $C^{k+1-j}$ with $\| Dh_j (\zetaj) \| < 1$ for all $\zetaj \in K_j$ and all applicable higher derivatives bounded on $K_j$. Therefore, $h_0 : K_0 \rightarrow \RRn$ is now $C^{k+1}$ with its first to $(k+1)$-th derivatives bounded on $K_0$. Furthermore, 
\begin{multline*}
D^{k+1} h_0(z)(\bfz_1, ..., \bfz_{k+1}) \\
= \sigma_k^{-1} \cdots \sigma_1^{-k} \vectnm^{-1}\!\big( \cdots \vectnmkMm^{-1}\!( Dh_k(\sigma_k^{-1} \cdots \sigma_1^{-1} z)\, \bfz_1 )\, \bfz_2 \cdots \big)\, \bfz_{k+1}
\end{multline*}
for any $z \in K_0$ and any $\bfz_1, ..., \bfz_{k+1} \in \RRm$.

Finally, we remark that the inequality (\ref{ineq:hyper:sysk}) is guaranteed by (\ref{ineq:Cr0}) provided that $k+1 \le r$. Therefore, the induction has to cease at $k=r-1$ with the outcome of the last iteration being that $h_0$ is $C^{k+1} = C^r$. This concludes the proof of Theorem \ref{thm2}.

\section{Applications} \label{App}

\subsection{Invariant Tori in a $2$-Parameter Family of Systems} \label{ex1}

Consider the following $2$-parameter family of systems:
\begin{equation} \label{example1}
\begin{split}
\dot R ={}& R(8\beta-R) + \sin\phi_1 + \sin\phi_2 \,, \\
\dot \phi_1 ={}& \beta R + \beta^2 \sin\phi_1 \sin\phi_2 \,, \\
\dot \phi_2 ={}& \omega \,, 
\end{split}
\end{equation}
where $R \in \RR$, $\phi_1 \in [0, 2\pi]$ and $\phi_2 \in [0, 2\pi]$ are two angular variables with the end points $0$ and $2\pi$ identified, and $\beta \in (0, \infty)$ and $\omega \in \RR$ are the two parameters. We are interested in determining for what parameter values (\ref{example1}) has $C^r$ ($r \ge 1$) invariant tori. However, even with the help of numerical methods, it is a challenging task to identify the exact set of parameters for the existence of $C^r$ invariant tori. Here we apply Theorem \ref{thm1} and Theorem \ref{thm3} to derive sufficient conditions on the parameters for the existence of $C^r$ invariant tori for (\ref{example1}).

We need to modify (\ref{example1}) slightly to fit the form of (\ref{sys}). This is done by lifting $\phi_1$ and $\phi_2$ to $\RR$. Furthermore, we rescale $\phi_1$ to $\theta_1 := k^{-1} \phi_1$ and $\phi_2$ to $\theta_2 := (k \gamma)^{-1} \phi_2$ with the scaling factors $k > 0$ and $\gamma > 0$ to be determined later. The purpose of this rescaling of variables is to obtain a larger set of parameters for which the existence of $C^r$ invariant tori can be guaranteed by Theorem \ref{thm1} and Theorem \ref{thm3}. Write $\theta = (\theta_1, \theta_2)$. Then (\ref{example1}) becomes
\begin{equation} \label{example1:res}
\begin{split}
\dot R ={}& f(R,\theta) := R(8\beta-R) + \sin k\theta_1 + \sin k \gamma \theta_2 \,, \\
\dot \theta ={}& g(R,\theta) := \begin{pmatrix} \tfrac{\beta}{k} R + \tfrac{\beta^2}{k} \sin k\theta_1 \sin k \gamma \theta_2 \\ \tfrac{\omega}{k \gamma} \end{pmatrix} , 
\end{split}
\end{equation}
where we have treated $R \in \RR$ and $\theta \in \RR^2$ as ``$a$'' and ``$z$'' of (\ref{sys}), respectively. In addition, for any $\delta > 0$, define the set $\Ud$ as follows:
\begin{equation*}
\Ud := \big\{ (R, \theta) : |R| < \delta, \theta \in \RR^2 \big\} \,,
\end{equation*}
which clearly satisfies Hypothesis \ref{hypo:U}.

In order to apply Theorem \ref{thm1}, we have to first establish the existence of a positively invariant set $\Gamma$ that is contained in $\Ud$ and satisfies 
\begin{equation} \label{proj:Gamma1}
\Pibot(\Gamma) = \Pibot(\Ud) = \RR^2 \,,
\end{equation} 
where $\Pibot$ is the projection onto the $\theta$-coordinate. We define $\Gamma \subset \Ud$ to be the set of points that stay in $\Ud$ forever along solution trajectories of (\ref{example1:res}) in forward time. Obviously, $\Gamma$ is the largest positively invariant subset of $\Ud$ by this definition. Next, we show that $\Gamma$ satisfies (\ref{proj:Gamma1}) using the {\Wa} theorem. 

Note that at the boundary of $\Ud$,
\begin{align*}
\dot R \ge{}& \delta (8\beta - \delta) - 2 && \hspace{-0.4in} \text{for any $(R, \theta)$ with $R = \delta$,} \\
\dot R \le{}& -\delta (8\beta + \delta) + 2 && \hspace{-0.4in} \text{for any $(R, \theta)$ with $R = -\delta$.}
\end{align*}
Thus, for any $\beta > 0$ and $\delta > 0$ satisfying 
\begin{equation} \label{ineq:example1:Wa}
\delta (8\beta - \delta) - 2 > 0 \,,
\end{equation}
$\dot R > 0$ for any $(R, \theta)$ with $R = \delta$, and $\dot R < 0$ for any $(R, \theta)$ with $R = -\delta$. Then it can be easily verified that $\closure(\Ud)$ is a {\Wa} set with its boundary $\partial \Ud$ being the set of points which leave $\closure(\Ud)$ immediately along solution trajectories of (\ref{example1:res}) in forward time. Let $\Ud^0 \subseteq \closure(\Ud)$ be the set of points that do not stay in $\closure(\Ud)$ forever in forward time. Then by the {\Wa} theorem, there exists a continuous function $\CR : \Ud^0 \times [0,1] \rightarrow \Ud^0$ such that $\CR$ is a strong deformation retraction of $\Ud^0$ onto $\partial \Ud$. Now suppose that for some $\theta^*$,
\begin{equation*}
\Gamma \textcap \big\{ (R, \theta^*): |R| \le \delta \big\} = \emptyset \,.
\end{equation*}
Then the set $\big\{ (R, \theta^*): |R| \le \delta \big\}$ is contained inside the domain of the continuous function $\CR(\,\cdot\,,1) : \Ud^0 \rightarrow \partial \Ud$. Immediately, we obtain a continuous function defined on $[-\delta, \delta]$ as follows:
\begin{equation*}
R \mapsto \Pi \circ \CR((R, \theta^*), 1) \,,
\end{equation*}
where $\Pi$ is the projection onto the $R$-coordinate. However, the above continuous function has the impossible property that 
\begin{align*}
& \Pi \circ \CR((R, \theta^*), 1) \in \big\{ -\delta, \delta \big\} \,, \\
& \Pi \circ \CR((\delta, \theta^*), 1) = \delta \,, \\
& \Pi \circ \CR((-\delta, \theta^*), 1) = -\delta \,.
\end{align*}
Thus, $\Gamma \textcap \big\{ (R, \theta): |R| \le \delta \big\} \neq \emptyset$ for any $\theta \in \RR^2$. This confirms (\ref{proj:Gamma1}).

Next, we verify Hypothesis \ref{hypo:C1} and Hypothesis \refHypoCr for (\ref{example1:res}) on $\Ud$. Straightforward calculation shows that for any $(R,\theta) \in \Ud$ and any $(R',\theta') \in \RR \times \RR^2$,
\begin{align*}
\langle R' , \DRf(R,\theta)\, R' \rangle \ge{}& (8\beta-2\delta) | R' |^2 \,, \\
\langle \theta' , \Dthetag(R,\theta)\, \theta' \rangle \le{}& \beta^2\sqrt{1+\gamma^2}\, \| \theta' \|^2 \,, \\ 
\| \Dthetaf(R,\theta) \| \le{}& k\sqrt{1+\gamma^2} \,, \\
\| \DRg(R,\theta) \| ={}& \tfrac{\beta}{k} \,.
\end{align*}
It follows that for any $\beta > 0$, $\delta > 0$, and $k > 0$ satisfying
\begin{equation} \label{ineq:example1:C1}
8\beta - 2\delta - \beta^2 - k - \tfrac{\beta}{k} > 0 \,,
\end{equation}
there exist sufficiently small $\gamma > 0$ and $c_1 > 0$, both depending on $\beta$, $\delta$, and $k$, such that 
\begin{equation*}
8\beta - 2\delta \ge \beta^2 \sqrt{1+\gamma^2} + k \sqrt{1+\gamma^2} + \tfrac{\beta}{k} + c_1 \,,
\end{equation*}
which implies (\ref{ineq:C1}). Therefore, by Theorem \ref{thm1}, (\ref{ineq:example1:Wa}) and (\ref{ineq:example1:C1}) together guarantee that $\Gamma$ is a $C^1$ positively invariant manifold. Furthermore, if there exists an integer $r \ge 2$ such that $\beta$, $\delta$, and $k$ also satisfy
\begin{equation} \label{ineq:example1:Cr}
8\beta - 2\delta - r\beta^2 - (r+1) \tfrac{\beta}{k} > 0 \,,
\end{equation}
then, taking a sufficiently small $c_r > 0$ and reducing $\gamma$ further if necessary, we have 
\begin{equation*} 
8\beta - 2\delta \ge r \beta^2 \sqrt{1+\gamma^2} + (r+1) \tfrac{\beta}{k} + c_r \,,
\end{equation*}
which implies (\refineqCr). Consequently, if $\beta$, $\delta$, and $k$ satisfy (\ref{ineq:example1:Wa})--(\ref{ineq:example1:Cr}) for some $r \ge 2$, then $\Gamma$ is in fact a $C^r$ manifold by Theorem \ref{thm3}.

Note that the inequalities (\ref{ineq:example1:Wa})--(\ref{ineq:example1:Cr}) involve the auxiliary parameters $\delta$ and $k$, which do not appear in the original system (\ref{example1}). In order to eliminate these auxiliary parameters to obtain conditions on $\beta$ only, we consider the $\beta$-projection of the solution set of (\ref{ineq:example1:Wa})--(\ref{ineq:example1:Cr}). Let $Q_1$ denote the set of triples $(\beta, \delta, k)$ that satisfy both (\ref{ineq:example1:Wa}) and (\ref{ineq:example1:C1}). In addition, for $r = 2, 3, 4, ...$, let $Q_r$ denote the set of triples $(\beta, \delta, k)$ that satisfy (\ref{ineq:example1:Wa}), (\ref{ineq:example1:C1}), and (\ref{ineq:example1:Cr}) for the corresponding $r$. Elementary computation reveals that $Q_r = \emptyset$ for all $r \ge 8$. Furthermore, we find that for each $r = 1, ..., 7$, the $\beta$-projection of $Q_r$ is a connected interval $(\beta_r^-, \beta_r^+)$ with the approximate values of $\beta_r^-$ and $\beta_r^+$ given in the table below. Note that $(\beta_1^-, \beta_1^+) \supset (\beta_2^-, \beta_2^+) \supset \cdots \supset (\beta_7^-, \beta_7^+)$.

\begin{table}[h]
\begin{tabular}{*{8}{| c} |}
\hline 
$r$ & 1 & 2 & 3 & 4 & 5 & 6 & 7 \\
\hline
$\beta_r^-$ & 0.395 & 0.404 & 0.420 & 0.441 & 0.472 & 0.518 & 0.634 \\
\hline
$\beta_r^+$ & 7.248 & 3.887 & 2.542 & 1.849 & 1.412 & 1.093 & 0.781 \\
\hline
\end{tabular}
\medskip
\caption{$\beta_r^-$ and $\beta_r^+$.}
\label{tab}
\vspace{-0.2in}
\end{table}

Now consider (\ref{example1}) with an arbitrary $\omega \in \RR$ and a $\beta$ drawn from one of the intervals $(\beta_r^-, \beta_r^+)$, $r = 1, ..., 7$. Since $(\beta_r^-, \beta_r^+)$ is the $\beta$-projection of $Q_r$, we can always find appropriate $\delta$, $k$, and correspondingly a sufficiently small $\gamma$ such that Theorem \ref{thm1} (if $r = 1$) or Theorem \ref{thm3} (if $2 \le r \le 7$) is applicable to the modified system (\ref{example1:res}) with $(\beta, \omega, k, \gamma)$ on the domain $\Ud$. Then we establish the existence of the $C^r$ positively invariant manifold $\Gamma$, which, however, depends on $\beta$ and $\omega$ as well as the choice of the auxiliary parameters $k$, $\gamma$, and $\delta$. Specifically, by Theorem \ref{thm1} or Theorem \ref{thm3}, $\Gamma$ is the graph of a $C^r$ function $h( \cdot , \cdot ; \beta, \omega, k, \gamma) : \RR^2 \rightarrow (-\delta, \delta)$, i.e., 
\begin{equation*}
\Gamma = \big\{ (h(\theta_1,\theta_2; \beta, \omega, k, \gamma), \theta_1, \theta_2) : (\theta_1, \theta_2) \in \RR^2 \big\} \subset \Ud \,.
\end{equation*}
Consequently, for (\ref{example1}) with $(\phi_1, \phi_2)$ lifted to $\RR^2$, we obtain a $C^r$ positively invariant manifold
\begin{equation} \label{example1:tori}
\big\{ (h( \tfrac{\phi_1}{k},\tfrac{\phi_2}{k\gamma}; \beta, \omega, k, \gamma), \phi_1, \phi_2) : (\phi_1, \phi_2) \in \RR^2 \big\} \,.
\end{equation}
On the other hand, suppose we can choose $k'$, $\gamma'$, and $\delta'$ such that for (\ref{example1:res}) with $(\beta, \omega, k', \gamma')$ there is also a $C^r$ positively invariant manifold $\Gamma'$, which is the graph of a $C^r$ function $h( \cdot , \cdot ; \beta, \omega, k', \gamma') : \RR^2 \rightarrow (-\delta', \delta')$, i.e., 
\begin{equation*}
\Gamma' = \big\{ (h(\theta'_1,\theta'_2; \beta, \omega, k', \gamma'), \theta'_1, \theta'_2) : (\theta'_1, \theta'_2) \in \RR^2 \big\} \subset \Udp \,.
\end{equation*}
Then, for (\ref{example1}) with $(\phi_1, \phi_2)$ lifted to $\RR^2$, we obtain ``another'' $C^r$ positively invariant manifold
\begin{equation*}
\big\{ (h( \tfrac{\phi_1}{k'},\tfrac{\phi_2}{k'\gamma'}; \beta, \omega, k', \gamma'), \phi_1, \phi_2) : (\phi_1, \phi_2) \in \RR^2 \big\} \,.
\end{equation*}
We show that in fact
\begin{equation} \label{hkg-vs-hkpgp}
h( \tfrac{\phi_1}{k},\tfrac{\phi_2}{k\gamma}; \beta, \omega, k, \gamma) = h(\tfrac{\phi_1}{k'},\tfrac{\phi_2}{k'\gamma'}; \beta, \omega, k', \gamma') \text{ for all } (\phi_1, \phi_2) \in \RR^2 \,.
\end{equation}
Without loss of generality, we assume $\delta \le \delta'$. It follows that for all $(\theta_1, \theta_2) \in \RR^2$,
\begin{equation*}
| h(\theta_1, \theta_2; \beta, \omega, k, \gamma) | < \delta' \,.
\end{equation*}
In addition, by the relation $(\phi_1, \phi_2) = (k\theta_1, k\gamma \theta_2) =(k'\theta'_1, k'\gamma'\theta'_2)$, the positive invariance of $\Gamma$ under the flow of (\ref{example1:res}) with $(\beta, \omega, k, \gamma)$ implies that 
\begin{equation*}
\tilde{\Gamma} := \big\{ (h( \tfrac{k'\theta'_1}{k}, \tfrac{k'\gamma'\theta'_2}{k\gamma}; \beta, \omega, k, \gamma), \theta'_1, \theta'_2) : (\theta'_1, \theta'_2) \in \RR^2 \big\} \subset \Udp 
\end{equation*}
is positively invariant under the flow of (\ref{example1:res}) with $(\beta, \omega, k', \gamma')$. Recall that by definition $\Gamma'$ is the largest positively invariant subset of $\Udp$ for (\ref{example1:res}) with $(\beta, \omega, k', \gamma')$. Thus, $\tilde{\Gamma} \subseteq \Gamma'$, which implies that for all $(\theta'_1, \theta'_2) \in \RR^2$,
\begin{equation*}
h( \tfrac{k'\theta'_1}{k},\tfrac{k'\gamma'\theta'_2}{k\gamma}; \beta, \omega, k, \gamma) = h(\theta'_1, \theta'_2; \beta, \omega, k', \gamma') \,.
\end{equation*}
Then we obtain (\ref{hkg-vs-hkpgp}) by substituting $\theta'_1 = \tfrac{\phi_1}{k'}$ and $\theta'_2 = \tfrac{\phi_2}{k'\gamma'}$ into the above identity. Therefore, for each $(\beta, \omega) \in (\beta_r^-, \beta_r^+) \times \RR$,  (\ref{example1}) with $(\phi_1, \phi_2)$ lifted to $\RR^2$ has a unique $C^r$ positively invariant manifold that is contained in $\Ud$ for a certain $\delta > 0$ and that is the graph of a $C^r$ function from $\RR^2$ to $\RR$ as given by (\ref{example1:tori}).

We still need to show that (\ref{example1:tori}) corresponds to a $C^r$ invariant torus. By the definition of $\Gamma$, the trajectory that passes through $(h(\theta_1,\theta_2; \beta, \omega, k, \gamma), \theta_1, \theta_2)$ is contained in $\Ud$ forever in forward time. Then by the periodicity of (\ref{example1:res}) in $\theta_1$ and $\theta_2$, the trajectory that passes through $(h(\theta_1,\theta_2; \beta, \omega, k, \gamma), \theta_1 \!+\! \frac{2n\pi}{k}, \theta_2 \!+\! \frac{2m\pi}{k\gamma})$ for any integers $n$ and $m$ is also contained in $\Ud$ forever in forward time. Thus $(h(\theta_1,\theta_2; \beta, \omega, k, \gamma), \theta_1 \!+\! \frac{2n\pi}{k}, \theta_2 \!+\! \frac{2m\pi}{k\gamma}) \in \Gamma$, which implies that
\begin{equation} \label{h:periodicity}
h(\theta_1,\theta_2; \beta, \omega, k, \gamma) = h(\theta_1 \!+\! \tfrac{2n\pi}{k}, \theta_2 \!+\! \tfrac{2m\pi}{k\gamma}; \beta, \omega, k, \gamma) \,.
\end{equation}
For $(\phi_1,\phi_2) \in \RR^2$, let $\rho^u(\phi_1,\phi_2; \beta, \omega) := h( \tfrac{\phi_1}{k},\tfrac{\phi_2}{k\gamma}; \beta, \omega, k, \gamma)$, which is well defined in view of (\ref{hkg-vs-hkpgp}). Clearly, $\rho^u(\phi_1 , \phi_2 ; \beta, \omega)$ is $C^r$ in $\phi_1$ and $\phi_2$, and it is $2\pi$-periodic in $\phi_1$ and $\phi_2$ by (\ref{h:periodicity}). Therefore, (\ref{example1:tori}) corresponds to the $C^r$ torus
\begin{gather*}
T^u_{\beta, \omega} := \big\{ ( \rho^u(\phi_1, \phi_2 ; \beta, \omega), \phi_1, \phi_2) : \phi_1, \phi_2 \in \RR (\Mod{2\pi}) \big\} \,.
\end{gather*}
It is easy to verify that $T^u_{\beta, \omega}$ is invariant under the flow of (\ref{example1}) (with the corresponding $\beta$ and $\omega$) in both forward time and backward time. Note that the superscript ``$u$'' reflects the obvious fact that $T^u_{\beta, \omega}$ is unstable.

After establishing $T^u_{\beta, \omega}$, we also consider the time reversal of (\ref{example1:res}) on the domain $\big\{ (R, \theta) : |R - 8\beta| < \delta,\, \theta \in \RR^2 \big\}$. By repeating all the analysis presented so far in this subsection, we find that for any $(\beta, \omega) \in (\beta_r^-, \beta_r^+) \times \RR$ with $1 \le r \le 7$, (\ref{example1}) has another $C^r$ invariant torus 
\begin{gather*}
T^s_{\beta, \omega} := \big\{ ( \rho^s(\phi_1, \phi_2 ; \beta, \omega), \phi_1, \phi_2) : \phi_1, \phi_2 \in \RR (\Mod{2\pi}) \big\} \,,
\end{gather*}
where $\rho^s(\phi_1 , \phi_2 ; \beta, \omega)$ is $C^r$ and $2\pi$-periodic in $\phi_1$ and $\phi_2$ for $(\phi_1, \phi_2) \in \RR^2$. $T^s_{\beta, \omega}$ is stable and contained in the domain $\big\{ (R, \theta) : |R - 8\beta| < \delta,\, \theta \in \RR^2 \big\}$ for a certain $\delta > 0$.  

Furthermore, for any $(\beta_0, \omega_0) \in (\beta_r^-, \beta_r^+) \times \RR$ with $1 \le r \le 7$, let $\beta = \beta_0 + \epsilon \xi_1$ and $\omega = \omega_0 + \epsilon \xi_2$ for a sufficiently small $\epsilon > 0$ such that $\beta = \beta_0 + \epsilon \xi_1 \in (\beta_r^-, \beta_r^+)$ for all $|\xi_1| < 1$. Next, we append $\dot \xi_1 = 0$ and $\dot \xi_2 = 0$ to (\ref{example1}) to form an enlarged system with the additional components $\xi_1$ and $\xi_2$. Then by adapting our previous analysis for the enlarged system, reducing $\epsilon$ further if necessary, we obtain that the enlarged system has $C^r$ invariant manifolds $M^u$ and $M^s$ which are the graphs of the $C^r$ functions, respectively,  
\begin{align*}
(\phi_1, \phi_2, \xi_1, \xi_2) \mapsto{}& \rho^u(\phi_1, \phi_2 ; \beta_0 + \epsilon \xi_1, \omega_0 + \epsilon \xi_2) \text{ and} \\
(\phi_1, \phi_2, \xi_1, \xi_2) \mapsto{}& \rho^s(\phi_1, \phi_2 ; \beta_0 + \epsilon \xi_1, \omega_0 + \epsilon \xi_2)
\end{align*}
on the domain $\RR (\Mod{2\pi}) \times \RR (\Mod{2\pi}) \times (-1, 1) \times (-1, 1)$. It follows that $\rho^u(\phi_1 , \phi_2 ; \beta, \omega)$ and $\rho^s(\phi_1 , \phi_2 ; \beta, \omega)$ are also $C^r$ in $\beta$ and $\omega$ for $(\beta, \omega) \in (\beta_r^-, \beta_r^+) \times \RR$. Therefore, for each $r = 1, ..., 7$, $\big\{ T^u_{\beta, \omega} : (\beta, \omega) \in (\beta_r^-, \beta_r^+) \times \RR \big\}$ and $\big\{ T^s_{\beta, \omega} : (\beta, \omega) \in (\beta_r^-, \beta_r^+) \times \RR \big\}$ are $C^r$ families of $C^r$ invariant tori.

\begin{figure}[h]
\centering
\includegraphics[scale=1]{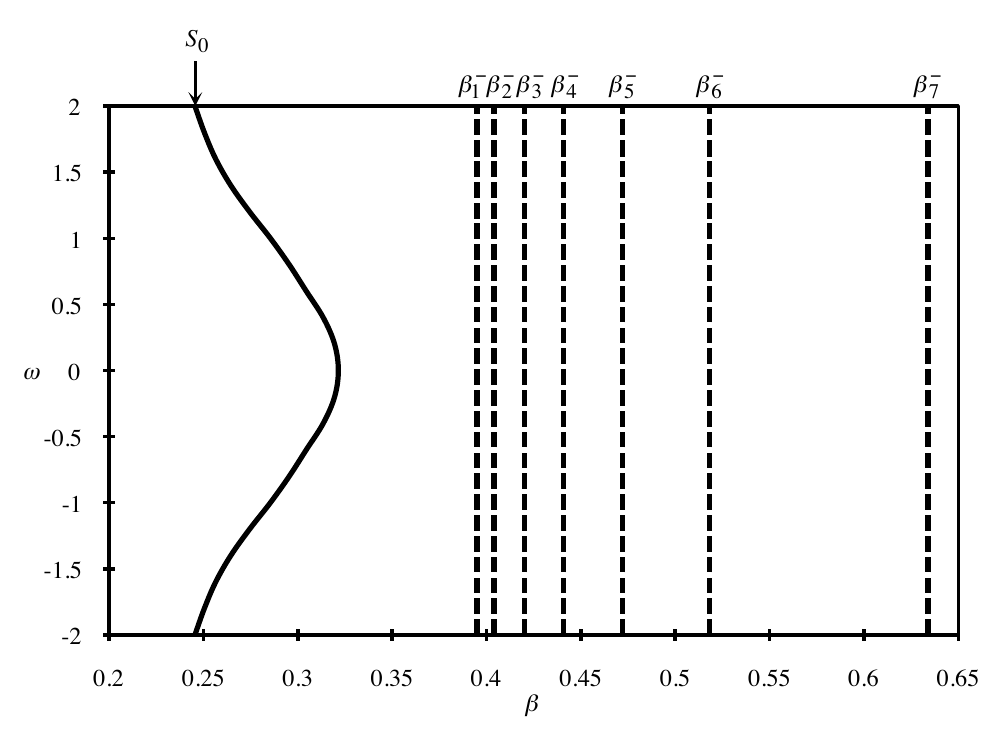}
\caption{$\beta_r^-$ and the curve $S_0$ on the $(\beta, \omega)$-plane. No invariant torus exists for (\ref{example1}) with any $(\beta, \omega)$ in the region bounded between the line $\beta = 0$ and the curve $S_0$.}
\label{fig1}
\end{figure}

\begin{figure}[h]
\centering
\includegraphics[scale=1]{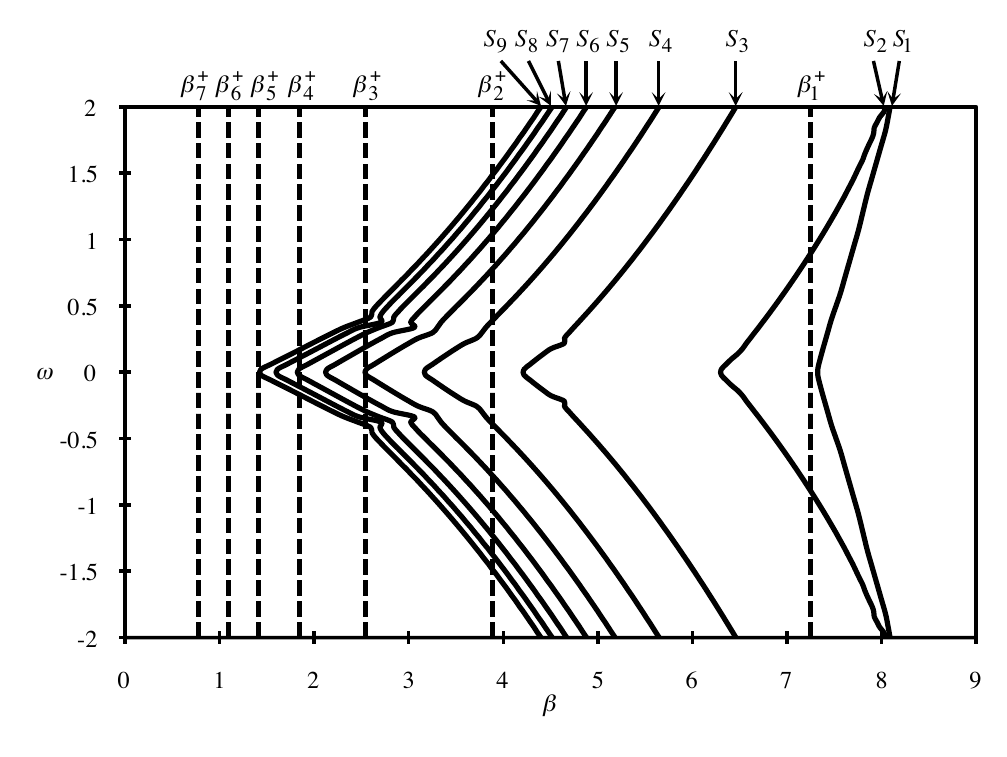}
\caption{$\beta_r^+$ and the curves $S_i$ on the $(\beta, \omega)$-plane. $T^u_{\beta, \omega}$ does not exist for any $(\beta, \omega)$ on the right-hand side of $S_1$, and $T^u_{\beta, \omega}$ is exactly $C^r$ (i.e., not $C^{r+1}$) for any $(\beta, \omega)$ in the region bounded between $S_{r+1}$ and $S_r$ ($r = 1, ..., 8$). Note that our theoretical estimates predict that $T^u_{\beta, \omega}$ is (at least) $C^r$ for any $(\beta, \omega) \in (\beta_r^-, \beta_r^+) \times \RR$ with $1 \le r \le 7$.}
\label{fig2}
\end{figure}

Finally, we compare our estimates (i.e., $\beta_r^-$ and $\beta_r^+$ in Table \ref{tab}) to some preliminary numerical results:

\begin{enumerate}

\item
By numerically integrating (\ref{example1}), we identify a region on the $(\beta, \omega)$-plane bounded between the line $\beta = 0$ and the curve $S_0$ (see Figure \ref{fig1}) such that for (\ref{example1}) with any $(\beta, \omega)$ in this region, solution trajectories that start on the plane $\big\{ (R, \phi_1, \phi_2) : R = 5 \big\}$ all reach the plane $\big\{ (R, \phi_1, \phi_2) : R = -5 \big\}$. Note that for (\ref{example1}) with $\beta < \beta_1^-$, $\dot R < 0$ if $R \le -5$ or $R \ge 5$. Thus, we conclude that no invariant torus exists for (\ref{example1}) with any $(\beta, \omega)$ in this region. In addition, when we continue the tori $T^u_{\beta, \omega}$ and $T^s_{\beta, \omega}$ along a path with constant $\omega$ starting at a $\beta$ between $\beta_1^-$ and $\beta_2^-$ and moving towards the curve $S_0$, we observe that $T^u_{\beta, \omega}$ and $T^s_{\beta, \omega}$ continue to exist even for $\beta < \beta_1^-$, but as $(\beta, \omega)$ gets close to $S_0$, both tori go through a series of bifurcations and then disappear. We have not yet studied these bifurcations in detail since it requires a much more delicate numerical analysis, which is out of the scope of this work.

\item
By sweeping the region $\beta_7^+ \le \beta \le 9$ and $-2 \le \omega \le 2$, we find that for some $(\beta, \omega)$, the $\alpha$-limit set of $T^u_{\beta, \omega}$ consists of a saddle-type periodic orbit $P_1$ and a non-saddle-type unstable periodic orbit $P_2$. In addition, the tangent bundle of $\RR^3$ restricted to $P_2$ has the splitting: $T \RR^3 |_{P_2} = N_1 \oplus N_2 \oplus TP_2$ such that both $N_1$ and $N_2$ are invariant under the linearized flow of (\ref{example1}) along $P_2$ and $T^u_{\beta, \omega}$ is tangent to $N_2$ along $P_2$. Since the linearized flow along $P_2$ expands both $N_1$ and $N_2$, the ratio between the expansion rate in $N_1$ and the expansion rate in $N_2$ determines the smoothness of $T^u_{\beta, \omega}$. By continuing $T^u_{\beta, \omega}$ and $P_2$ and monitoring the expansion rates in $N_1$ and $N_2$, we obtain the level curves $S_2$, ..., $S_9$ (see Figure \ref{fig2}), along which the ratio between the expansion rates in $N_1$ and $N_2$ are integers $2$, ..., $9$, respectively. In addition, we find that for any $(\beta, \omega)$ on the curve $S_1$, the bundles $N_1$ and $N_2$ merge, and then for $(\beta, \omega)$ on the right-hand side of $S_1$, the two bundles emerge again, however, both rotating $n \pi$ ($n \in \mathbb{Z} \setminus \big\{ 0 \big\}$) for one period of $P_2$. Thus, $T^u_{\beta, \omega}$ does not exist for any $(\beta, \omega)$ on the right-hand side of $S_1$, and $T^u_{\beta, \omega}$ is exactly $C^r$ (i.e., not $C^{r+1}$) for any $(\beta, \omega)$ in the region bounded between $S_{r+1}$ and $S_r$ ($r = 1, ..., 8$). We remark that $T^u_{\beta, \omega}$ can be $C^r$ with $r \ge 10$ for some $(\beta, \omega)$ on the left-hand side of $S_9$.
\end{enumerate}

\subsection{Invariant Torus in a System with Rapid Oscillations} \label{ex2}

In \cite{ChLi00}, Chicone and Liu studied the existence of an invariant torus in the following system:
\begin{equation} \label{sys:ChLi}
\begin{split}
\dot\rho ={}& \Delta(\sigma) \rho + \Lambda(\sigma) + \mu R(\rho,\sigma,\tau/\mu^2,\mu) \,, \\
\dot\sigma ={}& \rho + \mu S(\rho,\sigma,\tau/\mu^2,\mu) \,, \\
\dot\tau ={}& 1 \,,
\end{split}
\end{equation}
where $\dot{} = \frac{d}{ds}$ for a slow time $s$, $\mu > 0$ is a real parameter, $\rho \in \RR$, and $\sigma \in [0, 2\pi]$ and $\tau \in [0, 2\pi\mu^2/\omega]$ for a fixed $\omega > 0$ are angular variables with the corresponding end points identified. In addition, (\ref{sys:ChLi}) is assumed to be class $C^{\infty}$. Note that $\mu R$ and $\mu S$ are fast oscillatory in $\tau$ for $\mu$ small.

After the truncation of the $\CO(\mu)$ terms $\mu R$ and $\mu S$, (\ref{sys:ChLi}) becomes 
\begin{equation} \label{sys:ChLi0}
\begin{split}
\dot\rho ={}& \Delta(\sigma) \rho + \Lambda(\sigma) \,, \\
\dot\sigma ={}& \rho \,, \\
\dot\tau ={}& 1 \,.
\end{split}
\end{equation}
It is assumed in \cite{ChLi00} that $\Lambda(\sigma) \neq 0$ and $\Delta(\sigma)<0$ for all $\sigma \in [0, 2\pi]$. Then the $\rho$-$\sigma$ subsystem of (\ref{sys:ChLi0}) has an attracting limit cycle, whose existence can be proved by an easy application of the Poincar\'e-Bendixson theorem, and the suspension of this limit cycle in (\ref{sys:ChLi0}) forms an invariant torus $\CM_0$. An immediate question is whether or not $\CM_0$ persists in the system (\ref{sys:ChLi}) for $\mu>0$ but small. Note that the attracting torus $\CM_0$ is $r$-normally hyperbolic with respect to (\ref{sys:ChLi0}) for any integer $r \ge 1$ and its strength of normal hyperbolicity is independent of any parameter. However, the persistence theory of $r$-normally hyperbolic invariant manifolds by Fenichel \cite{Fe71} and Hirsch, Pugh, and Shub \cite{HiPuSh77} is not applicable in this case. The reason is that (\ref{sys:ChLi}) may not be $C^1$ close to (\ref{sys:ChLi0}) even for small $\mu$ since the partial derivatives of $\mu R$ and $\mu S$ with respect to $\tau$ can be very large for small $\mu$. On the other hand, by rescaling time to $t := s/\mu^2$ and taking $\theta := \tau/\mu^2$, we transform (\ref{sys:ChLi}) into (\ref{sys:ChLi:res}) with $w = (\rho, \sigma)$ and $\epsilon = \mu^2$. Then, in the $\CO(\epsilon)$ truncation of (\ref{sys:ChLi:res}) (i.e., $\dot{w} = \epsilon F_1(w)$, $\dot{\theta} = 1$), the normal hyperbolicity of the corresponding unperturbed invariant manifold depends on the small parameter $\epsilon$, which raises the problem of weak hyperbolicity as discussed in Subsection \ref{intro:WH}.

The main result of Chicone and Liu is that {\it for any integer\/ $r \ge 2$, (\ref{sys:ChLi}) with\/ $\mu > 0$ but sufficiently small has an\/ $r$-normally hyperbolic invariant torus\/ $C^r$ diffeomorphic to\/ $\CM_0$.} In addition, by time reversal, the same conclusion holds under the assumption that $\Lambda(\sigma) \neq 0$ and $\Delta(\sigma)>0$ for all $\sigma \in [0, 2\pi]$. Chicone and Liu also remarked that this result is not valid if $\Lambda$ is allowed to have zeros and the formulation of correct hypotheses needed to prove an analogous result in this case is ``an interesting open problem''. One of the difficulties is that even the existence of a $C^r$ invariant torus for (\ref{sys:ChLi0}) is not readily available if $\Lambda$ has zeros.

We now apply Theorem \ref{thm1} and Theorem \ref{thm3} to establish a $C^r$ ($r \ge 1$) invariant torus for (\ref{sys:ChLi}) even if $\Lambda$ has zeros. In preparation for our analysis, we lift $\sigma$ and $\tau$ to $\RR$ and then rescale them to $\theta_1 := k^{-1} \sigma$ and $\theta_2 := \tau/\mu^2$ with the scaling factor $k > 0$ to be determined later. In addition, we introduce an auxiliary variable $\theta_3 \in (-2, 2)$ to (\ref{sys:ChLi}) to form an enlarged system, which is written in the form of (\ref{sys}) as follows:
\begin{equation} \label{sys:ChLi:ext}
\begin{split}
\dot \rho ={}& f(\rho,\theta) := \Delta(k\theta_1) \rho + \Lambda(k\theta_1) + \theta_3 \mu R(\rho, k\theta_1, \theta_2, \mu) \,, \\
\dot \theta ={}& g(\rho,\theta) := \begin{pmatrix} \tfrac{1}{k} \rho + \tfrac{1}{k} \theta_3 \mu S(\rho, k\theta_1, \theta_2, \mu) \\ \tfrac{1}{\mu^2} \\ 0 \end{pmatrix} .
\end{split}
\end{equation}
where $\theta := (\theta_1, \theta_2, \theta_3)$, $\Delta$ and $\Lambda$ are $2\pi/k$-periodic in $\theta_1$, and $R$ and $S$ are $2\pi/k$-periodic in $\theta_1$ and $2\pi/\omega$-periodic in $\theta_2$. Note that we have treated $\rho$ and $\theta$ as ``$a$'' and ``$z$'' of (\ref{sys}), respectively.

Assume $\Delta(\sigma) > 0$ for all $\sigma \in [0, 2\pi]$. Define constants $E_1$, $E_2$, and $E$ as follows:
\begin{align}
E_1 :={}& \min_{\sigma \in [0, 2\pi]} \bigg(\! -\frac{\Lambda(\sigma)}{\Delta(\sigma)} \bigg) \,, \notag \\
E_2 :={}& \max_{\sigma \in [0, 2\pi]} \bigg(\! -\frac{\Lambda(\sigma)}{\Delta(\sigma)} \bigg) \,, \notag \\
E :={}& \max \big\{ |E_1|, |E_2| \big\} \,. \label{def:const:E}
\end{align}
We shall consider (\ref{sys:ChLi:ext}) on the domain 
\begin{equation*}
\Ud := (E_1-\delta, E_2+\delta) \times \RR^2 \times (-2, 2) \,,
\end{equation*}
where $\delta \in (0, 1]$ is another fixed constant to be determined. Obviously, $\Ud$ as defined above satisfies Hypothesis \ref{hypo:U}. Now, we need to establish the existence of a positively invariant set $\Gamma$ that is contained in $\Ud$ and satisfies 
\begin{equation} \label{proj:Gamma:ex2}
\Pibot(\Gamma) = \Pibot(\Ud) = \RR^2 \times (-2, 2) \,,
\end{equation} 
where $\Pibot$ is the projection onto the $\theta$-coordinate. In the same way as we have done for the example in Subsection \ref{ex1}, we define $\Gamma \subset \Ud$ to be the set of points that stay in $\Ud$ forever along solution trajectories of (\ref{sys:ChLi:ext}) in forward time so that $\Gamma$ is the largest positively invariant subset of $\Ud$ by definition. Then we show that $\Gamma$ satisfies (\ref{proj:Gamma:ex2}), again, using the {\Wa} theorem. 

By the definitions of $E_1$ and $E_2$ and the assumption that $\Delta(\sigma) > 0$ for all $\sigma \in [0, 2\pi]$, it is straightforward to verify that for $M := \min_{\sigma \in [0, 2\pi]} \Delta(\sigma) > 0$,
\begin{align*}
\dot\rho \ge{}& M \delta - 2 \mu \|R\|_0 & \text{for any $(\rho, \theta)$ with $\rho \ge E_2+\delta$,} \\
\dot\rho \le{}& -M \delta + 2 \mu \|R\|_0 & \text{for any $(\rho, \theta)$ with $\rho \le E_1-\delta$,}
\end{align*}
where $\|R\|_0$ is the uniform norm of $R(\rho, k\theta_1, \theta_2, \mu)$ on the domain $(E_1-1, E_2+1) \times \RR^2 \times (0, 1]$. Then for any $\delta, \mu \in (0,1]$ satisfying 
\begin{equation} \label{ineq:ex2:Wa}
M \delta - 2 \mu \|R\|_0 > 0 \,,
\end{equation}
$\dot\rho > 0$ for all $(\rho, \theta)$ with $\rho \ge E_2+\delta$, and $\dot\rho < 0$ for all $(\rho, \theta)$ with $\rho \le E_1-\delta$. It follows that $\big\{ E_1-\delta, E_2+\delta \big\}  \times \RR^2 \times (-2, 2)$ is the set of points which leave the set $[E_1-\delta, E_2+\delta] \times \RR^2 \times (-2, 2)$ immediately along solution trajectories of (\ref{sys:ChLi:ext}) in forward time. Then, the latter is a {\Wa} set, and by the same arguments as what we have used in Subsection \ref{ex1}, we can show that $\Gamma$ does satisfy (\ref{proj:Gamma:ex2}).

Next, we verify Hypothesis \ref{hypo:C1} and Hypothesis \refHypoCr for (\ref{sys:ChLi:ext}) on $\Ud$. Straightforward calculation shows that for any $(\rho,\theta) \in \Ud$ and any $(\rho',\theta') \in \RR \times \RR^3$,
\begin{align*}
\langle \rho' , \Drhof(\rho,\theta)\, \rho' \rangle \ge{}& (\Delta(k\theta_1) - 2 \mu \| D_1 R \|_0) | \rho' |^2 \,, \\
\langle \theta' , \Dthetag(\rho,\theta)\, \theta' \rangle \le{}& \mu ( \tfrac{1}{k^2} \|S\|_0^2 + 4 \|D_2 S\|_0^2 + \tfrac{4}{k^2} \|D_3 S\|_0^2 )^{\frac{1}{2}} \| \theta' \|^2 \,, \\ 
\| \Dthetaf(\rho,\theta) \| \le{}& \Big( k^2 \big( \|D\Delta(k\theta_1)\| (E+\delta) + \|D\Lambda(k\theta_1)\| + 2 \mu \|D_2R\|_0 \big)^2 \\
& + 4 \mu^2 \|D_3R\|_0^2 + \mu^2 \|R\|_0^2 \Big)^{\frac{1}{2}} \,, \\
\| \Drhog(\rho,\theta) \| \le{}& \tfrac{1}{k} (1 + 2 \mu \|D_1S\|_0) \,,
\end{align*}
where $D_i F = \frac{\partial}{\partial x_i} F(x_1, x_2, x_3, x_4)$ for $i = 1, 2, 3$, and $F = R$, or $S$, and $\|\cdot\|_0$ denotes the uniform norm of functions of $(x_1,x_2,x_3,x_4)$ on the domain $(E_1-1, E_2+1) \times \RR^2 \times (0, 1]$. Thus, for any triple $(\mu, k, \delta) \in (0, 1] \times (0, \infty) \times (0, 1]$ such that 
\begin{equation} \label{ineq:ex2:C1}
\begin{split}
& (\Delta(k\theta_1) - 2 \mu \| D_1 R \|_0) - \mu ( \tfrac{1}{k^2} \|S\|_0^2 + 4 \|D_2 S\|_0^2 + \tfrac{4}{k^2} \|D_3 S\|_0^2 )^{\frac{1}{2}} \\
& - \Big( k^2 \big( \|D\Delta(k\theta_1)\| (E+\delta) + \|D\Lambda(k\theta_1)\| + 2 \mu \|D_2R\|_0 \big)^2 \!+ 4 \mu^2 \|D_3R\|_0^2 + \mu^2 \|R\|_0^2 \Big)^{\frac{1}{2}} \\
& - \tfrac{1}{k} (1 + 2 \mu \|D_1S\|_0) > 0
\end{split}
\end{equation}
holds for all $k\theta_1 \in [0, 2\pi]$, (\ref{ineq:C1}) (with a sufficiently small $c_1$) and hence Hypothesis \ref{hypo:C1} are satisfied on the domain $\Ud$. Therefore, by Theorem \ref{thm1}, (\ref{ineq:ex2:Wa}) and (\ref{ineq:ex2:C1}) together guarantee that $\Gamma$ is a $C^1$ positively invariant manifold contained inside $\Ud$. Furthermore, if there also exists an integer $r \ge 2$ such that
\begin{equation} \label{ineq:ex2:Cr}
\begin{split}
& (\Delta(k\theta_1) - 2 \mu \| D_1 R \|_0) - r \mu ( \tfrac{1}{k^2} \|S\|_0^2 + 4 \|D_2 S\|_0^2 + \tfrac{4}{k^2} \|D_3 S\|_0^2 )^{\frac{1}{2}} \\
& - (r + 1) \tfrac{1}{k} (1 + 2 \mu \|D_1S\|_0) > 0
\end{split}
\end{equation}
holds for all $k\theta_1 \in [0, 2\pi]$, then (\refineqCr) (with a sufficiently small $c_r$) and hence Hypothesis \refHypoCr are satisfied on the domain $\Ud$. Consequently, if the triple $(\mu, k, \delta)$ satisfies (\ref{ineq:ex2:Wa})--(\ref{ineq:ex2:Cr}) for some $r \ge 2$, then $\Gamma$ is in fact a $C^r$ manifold by Theorem \ref{thm3}. 

Note that with (\ref{ineq:ex2:Wa})--(\ref{ineq:ex2:Cr}), we can determine for each positive integer $r$ a subset of $(0, 1]$ such that the system (\ref{sys:ChLi:ext}) with any $\mu$ chosen from this subset (if nonempty) possesses a $C^r$ positively invariant manifold $\Gamma$, which is contained inside $\Ud$. In principle, these subsets can be identified in the same way as we obtain Table \ref{tab} in Subsection \ref{ex1}. However, the analysis would be difficult without explicit expressions of $\Delta$, $\Lambda$, $R$, and $S$. On the other hand, if our only concern is what happens given that $\mu$ is sufficiently small, we can derive a relatively simple condition on $\Delta$ and $\Lambda$ to guarantee the existence of the $C^r$ ($r \ge 1$) positively invariant manifold $\Gamma$ for (\ref{sys:ChLi:ext}).

First, we make some important observations: if for some $k > 0$, 
\begin{equation} \label{precond:ex2:C1}
\Delta(\sigma) > k \big( \|D\Delta(\sigma)\| E + \|D\Lambda(\sigma)\| \big) + \tfrac{1}{k} \text{ for all $\sigma \in [0, 2\pi]$,}
\end{equation}
then we can always find an appropriate $\delta \in (0, 1]$ and correspondingly a sufficiently small $\mu_0 > 0$ such that for any $\mu \in (0, \mu_0]$ the triple $(\mu, k, \delta)$ satisfies both (\ref{ineq:ex2:Wa}) and (\ref{ineq:ex2:C1}); and if in addition to (\ref{precond:ex2:C1}) there is an integer $r \ge 2$ such that
\begin{equation*} 
\Delta(\sigma) > (r + 1) \tfrac{1}{k} \text{ for all $\sigma \in [0, 2\pi]$,}
\end{equation*}
then we can further reduce $\mu_0$ if necessary so that for any $\mu \in (0, \mu_0]$ the triple $(\mu, k, \delta)$ satisfies all (\ref{ineq:ex2:Wa})--(\ref{ineq:ex2:Cr}). Next, we notice the following facts: 
\begin{gather*}
\inf_{k > 0} \Big\{ \max_{\sigma \in [0, 2\pi]} \big\{ k \big( \|D\Delta(\sigma)\| E + \|D\Lambda(\sigma)\| \big) + \tfrac{1}{k} \big\} \Big\} = 2 \sqrt{L} \,, \\
\inf_{k > 0} \Big\{ \max_{\sigma \in [0, 2\pi]} \big\{ k \big( \|D\Delta(\sigma)\| E + \|D\Lambda(\sigma)\| \big) + \tfrac{1}{k}, (r + 1) \tfrac{1}{k} \big\} \Big\} = \tfrac{r + 1}{\sqrt{r}} \sqrt{L} \,,
\end{gather*}
where $L := \max_{\sigma \in [0, 2\pi]} \big\{ \|D\Delta(\sigma)\| E + \|D\Lambda(\sigma)\| \big\}$. It follows that if 
\begin{equation} \label{cond:ex2:Cr}
\Delta(\sigma) > \tfrac{r + 1}{\sqrt{r}} \sqrt{L} \text{ for all $\sigma \in [0, 2\pi]$,}
\end{equation}
then there always exist a $k > 0$, an appropriate $\delta \in (0, 1]$, and correspondingly a sufficiently small $\mu_0 > 0$ such that for any $\mu \in (0, \mu_0]$ the triple $(\mu, k, \delta)$ satisfies both (\ref{ineq:ex2:Wa}) and (\ref{ineq:ex2:C1}) when $r = 1$ in (\ref{cond:ex2:Cr}) and satisfies all (\ref{ineq:ex2:Wa})--(\ref{ineq:ex2:Cr}) when $r \ge 2$ in (\ref{cond:ex2:Cr}). Consequently, (\ref{sys:ChLi:ext}) with the corresponding $\mu$ and $k$ possesses a $C^r$ positively invariant manifold 
\begin{equation*}
\Gamma = \big\{ (h(\theta_1,\theta_2, \theta_3; \mu, k), \theta_1, \theta_2, \theta_3) : (\theta_1, \theta_2, \theta_3) \in \RR^2 \times (-2,2) \big\} \,,
\end{equation*}
which is contained in $\Ud$ with the corresponding $\delta$. The function $h( \cdot , \cdot , \cdot ; \mu, k) : \RR^2 \times (-2,2) \rightarrow (E_1-\delta, E_2+\delta)$ is $C^r$. 

Although $\Gamma$ depends on the choice of $(\mu, k, \delta)$ for (\ref{sys:ChLi:ext}) and $\Ud$, we can show by the same arguments as given in Subsection \ref{ex1} that, up to a rescaling in $\theta_1$, $\Gamma$ coincides with any $\Gamma'$ that is associated with any other admissible $(\mu, k', \delta')$ as long as $\mu$ remains the same. In particular, similar to (\ref{hkg-vs-hkpgp}), we have that
\begin{equation*}
h( \tfrac{\sigma}{k}, \theta_2, \theta_3; \mu, k) = h( \tfrac{\sigma}{k'}, \theta_2, \theta_3; \mu, k') \text{ for all } (\sigma, \tau) \in \RR^2 \,.
\end{equation*}
Furthermore, by the periodicity of (\ref{sys:ChLi:ext}) with respect to $\theta_1$ and $\theta_2$, we have that, similar to (\ref{h:periodicity}), $h(\theta_1,\theta_2, \theta_3; \mu, k)$ is $\tfrac{2\pi}{k}$-periodic in $\theta_1$ and $\tfrac{2\pi}{\omega}$-periodic in $\theta_2$. Then, for any $\beta \in [0,1]$, the $C^r$ submanifold $\Gamma |_{\theta_3 = \beta}$ (i.e., the section of $\Gamma$ at $\theta_3 = \beta$) corresponds to a $C^r$ torus 
\begin{equation*}
T_{\mu, \beta} := \big\{ ( \rho(\sigma, \tau; \beta, \mu), \sigma, \tau) : \sigma \in \RR (\Mod{2\pi}),\, \tau \in \RR (\Mod{ \tfrac{2\pi\mu^2}{\omega} }) \big\} \,,
\end{equation*}
where $\rho(\sigma, \tau; \beta, \mu) := h(\tfrac{\sigma}{k}, \tfrac{\tau}{\mu^2}, \beta; \mu, k)$. Recall that $\dot\rho > 0$ for all $(\rho, \theta)$ with $\rho \ge E_2+\delta$, $\dot\rho < 0$ for all $(\rho, \theta)$ with $\rho \le E_1-\delta$, and $\Gamma$ is the largest positively invariant subset of $\Ud$. Thus, $T_{\mu, \beta}$ is the unique invariant torus for the following system with the corresponding $\mu$ and $\beta$:
\begin{equation*} 
\begin{split}
\dot\rho ={}& \Delta(\sigma) \rho + \Lambda(\sigma) + \beta \mu R(\rho,\sigma,\tau/\mu^2,\mu) \,, \\
\dot\sigma ={}& \rho + \beta \mu S(\rho,\sigma,\tau/\mu^2,\mu) \,, \\
\dot\tau ={}& 1 \,.
\end{split}
\end{equation*}
Finally, since $h$ is $C^r$ with respect to $\theta_3$, we have that for every fixed $\mu \in (0, \mu_0]$, $\big\{ T_{\mu, \beta} : \beta \in [0, 1] \big\}$ forms a $C^r$ family of tori. Thus we have established simultaneously the existence of the invariant torus $T_{\mu, 1}$ for (\ref{sys:ChLi}), the existence of the invariant torus $\CM_0 := T_{\mu, 0}$ for (\ref{sys:ChLi0}), and the fact that $T_{\mu, 1}$ is $C^r$ diffeomorphic to $\CM_0$. Notice that the existence of $\CM_0$ for (\ref{sys:ChLi0}) is not among our assumptions. We summarize these results in the following proposition as an answer to the open question posed by Chicone and Liu in \cite{ChLi00}.

\begin{prop}
Suppose (\ref{sys:ChLi}) is class\/ $C^\infty$ and\/ $\Delta(\sigma) > 0$ for all\/ $\sigma \in [0, 2\pi]$. Let\/ $E$ be defined as (\ref{def:const:E}), and define\/ $L := \max_{\sigma \in [0, 2\pi]} \big\{ \|D\Delta(\sigma)\| E + \|D\Lambda(\sigma)\| \big\}$. If there is an integer\/ $r \ge 1$ such that\/ $\Delta(\sigma) > \tfrac{r + 1}{\sqrt{r}} \sqrt{L}$ for all\/ $\sigma \in [0, 2\pi]$, then
\begin{enumerate}
\item (\ref{sys:ChLi0}) has a unique\/ $C^r$ invariant torus\/ $\CM_0$; and
\item there exists a sufficiently small\/ $\mu_0 > 0$ such that for any\/ $\mu \in (0, \mu_0]$, 
(\ref{sys:ChLi}) has a unique\/ $C^r$ invariant torus\/ $T_{\mu, 1}$, which is\/ $C^r$ diffeomorphic to\/ $\CM_0$.
\end{enumerate}
\end{prop}

\subsection{Persistence of a Weakly Normally Hyperbolic Invariant Torus} \label{ex3}

In this subsection, we consider the persistence of a weakly normally hyperbolic invariant torus in the following system:
\begin{equation} \label{GWH}
\begin{split}
\dot{w} ={}& \epsilon F_1(w) + \epsilon^{1+\mu} F_2(w,\theta,\epsilon) \,, \\
\dot{\theta} ={}& \Theta_0 + \epsilon^{\nu} G_1(w,\epsilon) + \epsilon^{1+\gamma} G_2(w,\theta,\epsilon) \,,
\end{split} 
\end{equation}
where $w \in \RR^n$, $\theta \in \Torus^m$, $0 < \epsilon \ll 1$, and the power indices $\mu$, $\nu$, and $\gamma$ satisfy
\begin{align} 
\begin{split} \label{indices}
\mu >{}& 0\,, \\ 
\gamma >{}& 0\,, \\ 
1 \ge \nu \ge{}& 0\,, \\
\mu + \nu >{}& 1\,. 
\end{split}
\end{align}

\begin{asp} \label{asp:GWH:sm}
$F_1$ is $C^{r+1}$ ($r \ge 1$). $F_2$, $G_1$, and $G_2$ are $C^r$ with respect to $(w, \theta)$ and continuous with respect to $\epsilon$ for $\epsilon \in [0, \epsilon_0]$.
\end{asp}

Suppose the averaged equation $\dot{w} = \epsilon F_1(w)$ has a periodic orbit $S$. Note that the existence and the geometry of this periodic orbit is completely determined by the vector field $F_1(w)$ and independent of $\epsilon$. Furthermore, the $C^{r+1}$ smoothness of $F_1$ implies that $S$ is $C^{r+2}$. Thus $\To : = S \times \Torus^m$ is a $C^{r+2}$ invariant $(m+1)$-torus for the following truncated system with any $\epsilon > 0$:
\begin{equation} \label{GWH:tr}
\begin{split}
\dot{w} ={}& \epsilon F_1(w) \,, \\
\dot{\theta} ={}& \Theta_0 + \epsilon^{\nu} G_1(w,\epsilon) \,,
\end{split} 
\end{equation}
where $\epsilon^{1+\mu} F_2(w,\theta,\epsilon)$ and $\epsilon^{1+\gamma} G_2(w,\theta,\epsilon)$, the terms that are higher order in $\epsilon$, are excluded. We assume that $\To$ is normally hyperbolic with respect to the flow of (\ref{GWH:tr}) for any $\epsilon > 0$. This is true if and only if the periodic orbit $S$ is hyperbolic with respect to the flow of the averaged equation $\dot{w} = \epsilon F_1(w)$ for any $\epsilon > 0$. Thus we formulate the hyperbolicity assumption as follows:

\begin{asp} \label{asp:GWH:hyper}
Let\/ $\chi(\zeta)$ with\/ $\chi(0) \in S$ be a periodic solution of\/ $\frac{dw}{d\zeta} = F_1(w)$, and let $\zeta_0$ be the period of $\chi(\zeta)$. The linear variational equation\/ $\frac{dw}{d\zeta} = DF_1( \chi(\zeta) )\, w$ has\/ $n$ linearly independent solutions whose Lyapunov exponents\/ $r_1, ..., r_n$ satisfy 
\begin{gather*}
\begin{split}
\sigma_u :={}& \min \big\{ r_i : 1 \le i \le n_u \big\} > 0 \,, \\
-\sigma_s :={}& \max \big\{ r_i : n_u + 1 \le i \le n_u + n_s \big\} < 0 \,, 
\end{split}
\qquad \text{and} \quad
r_n = 0 \,,
\end{gather*}
where the integers\/ $n_u$ and\/ $n_s$ satisfy\/ $n_u \ge 1$, $n_s \ge 1$, and\/ $n_u + n_s + 1 = n$.
\end{asp}

We will prove the following theorem about the persistence of $\To$ in the full system (\ref{GWH}) for small $\epsilon$.

\begin{thm} \label{thm4}
Consider (\ref{GWH}) with\/ $\mu$, $\nu$, and\/ $\gamma$ satisfying (\ref{indices}). Suppose Assumption \ref{asp:GWH:sm} and Assumption \ref{asp:GWH:hyper} hold. Then there exists an\/ $\epsilon^* \in (0, \epsilon_0]$ such that for (\ref{GWH}) with any\/ $\epsilon \in (0, \epsilon^*]$, there is a unique invariant torus\/ $\Te$ inside an\/ $\CO(1)$-neighborhood (i.e., independent of\/ $\epsilon$) of\/ $\To$. In addition, $\Te$ is\/ $C^r$ diffeomorphic and\/ $\CO( \epsilon^{\mu} )$-close to\/ $\To$. In particular, $\Te$ has the parameterization 
\begin{equation} \label{Torus:para}
\Te = \big\{ (w, \theta) : w = \chi(\zeta) + \hat{n}(\zeta)\, \varrho(\zeta, \theta ; \epsilon),\, \zeta \in \RR (\Mod{\zeta_0}),\, \theta \in \Torus^m \big\} \,, 
\end{equation}
where\/ $\varrho(\cdot, \cdot ; \epsilon) : \RR (\Mod{\zeta_0}) \times \Torus^m \rightarrow \RR^{n-1}$ is an\/ $\CO( \epsilon^{\mu} )$, $C^r$ function and\/ $\hat{n} : \RR (\Mod{\zeta_0}) \rightarrow \RR^{n\times(n-1)}$ is a\/ $C^{r+1}$ matrix function such that for each\/ $\zeta \in \RR (\Mod{\zeta_0})$ the columns of\/ $\hat{n}(\zeta)$ form a basis of the normal space of the periodic orbit\/ $S$ in\/ $\RR^n$ at\/ $\chi(\zeta)$.
\end{thm}

Let $\eta: \RR \rightarrow \RR^{n\times(n-1)}$ be a $C^{r+1}$, bounded matrix function such that the $n \times n$ matrix $\big( \eta(\zeta) \; \frac{d}{d\zeta}\chi(\zeta) \big)$ is nonsingular and its inverse is bounded for all $\zeta \in \RR$. For a sufficiently small $\Delta > 0$, we make a $C^{r+1}$ change of variables in a small neighborhood of the periodic orbit $S$ as follows:
\begin{equation} \label{w:a-b-zeta}
w = \chi(\zeta) + \eta(\zeta) \big( \begin{smallmatrix} a \\ b \end{smallmatrix} \big) \,, 
\end{equation}
where $\zeta \in \RR$, $a \in \RR^{n_u}$ with $\|a\| < \Delta$, and $b \in \RR^{n_s}$ with $\|b\| < \Delta$. By Floquet's theorem, Assumption \ref{asp:GWH:hyper} implies that we can choose a $2 \zeta_0$-periodic $\eta(\zeta)$ such that in the $a$-$b$-$\zeta$ coordinates, the averaged equation $\dot{w} = \epsilon F_1(w)$ is transformed into the following normal form:
\begin{align} \label{ave:NF0}
\begin{split}
\dot a ={}& \epsilon ( A_0 a + A_1(a,b,\zeta) ) \,, \\
\dot b ={}& \epsilon ( B_0 b + B_1(a,b,\zeta) ) \,, \\
\dot \zeta ={}& \epsilon ( 1 + V_1(a,b,\zeta) ) \,,
\end{split}
\end{align}
where both the $n_u \times n_u$ constant matrix $A_0$ and the $n_s \times n_s$ constant matrix $B_0$ are in real Jordan form, and the functions $A_1(a,b,\zeta)$, $B_1(a,b,\zeta)$, and $V_1(a,b,\zeta)$ are all $\CO( \|a\|^2 + \|b\|^2 )$ and $2 \zeta_0$-periodic in $\zeta$. Furthermore, the real parts of the eigenvalues of $A_0$ coincide with the Lyapunov exponents $r_1, ..., r_{n_u}$, and the real parts of the eigenvalues of $B_0$ coincide with the Lyapunov exponents $r_{n_u+1}, ..., r_{n_u+n_s}$.

We will make a further change of coordinates by rescaling individual components of $a$ and $b$ so that we can verify (\ref{ineq:a-a}) for the transformed system. Since $A_0$ and $B_0$ are in real Jordan form, it suffices to illustrate how to rescale the components of $a$ that are associated with the same Jordan block $A_{0,j}$ of $A_0$. Suppose we have
\begin{align*}
A_{0,j} = 
\begin{pmatrix}
r_{\kappa} & \xi_{\kappa} & 1 & 0 & 0 & 0 \\
-\xi_{\kappa} & r_{\kappa} & 0 & 1 & 0 & 0 \\
0 & 0 & r_{\kappa} & \xi_{\kappa} & 1 & 0 \\
0 & 0 & -\xi_{\kappa} & r_{\kappa} & 0 & 1 \\
0 & 0 & 0 & 0 & r_{\kappa} & \xi_{\kappa} \\
0 & 0 & 0 & 0 & -\xi_{\kappa} & r_{\kappa} 
\end{pmatrix} ,
\end{align*}
where $r_{\kappa}$ is one of the positive Lyapunov exponents specified in Assumption \ref{asp:GWH:hyper}, and $r_{\kappa} \pm \xi_{\kappa}\, i$ are the pair of complex eigenvalues of $A_{0,j}$. Let $a_{j,1}, ..., a_{j,6}$ be the corresponding components of $a$ that are associated with $A_{0,j}$. Define
\begin{align*}
\begin{pmatrix}
p_{j,1} \\
\vdots \\
p_{j,6} 
\end{pmatrix} 
:= \Lambda_{j}
\begin{pmatrix}
a_{j,1} \\
\vdots \\
a_{j,6} 
\end{pmatrix} ,
\quad \text{where } 
\Lambda_{j}^{-1} = 
\begin{pmatrix}
1 & 0 & 0 & 0 & 0 & 0 \\
0 & 1 & 0 & 0 & 0 & 0 \\
0 & 0 & \lambda_{j,1} & 0 & 0 & 0 \\
0 & 0 & 0 & \lambda_{j,1} & 0 & 0 \\
0 & 0 & 0 & 0 & \lambda_{j,2} & 0 \\
0 & 0 & 0 & 0 & 0 & \lambda_{j,2}
\end{pmatrix} .
\end{align*}
It follows that 
\begin{align*}
P_{0,j} := \Lambda_j A_{0,j} \Lambda_j^{-1} =
\begin{pmatrix}
r_{\kappa} & \xi_{\kappa} & \lambda_{j,1} & 0 & 0 & 0 \\
-\xi_{\kappa} & r_{\kappa} & 0 & \lambda_{j,1} & 0 & 0 \\
0 & 0 & r_{\kappa} & \xi_{\kappa} & \tfrac{\lambda_{j,2}}{\lambda_{j,1}} & 0 \\
0 & 0 & -\xi_{\kappa} & r_{\kappa} & 0 & \tfrac{\lambda_{j,2}}{\lambda_{j,1}} \\
0 & 0 & 0 & 0 & r_{\kappa} & \xi_{\kappa} \\
0 & 0 & 0 & 0 & -\xi_{\kappa} & r_{\kappa} 
\end{pmatrix} .
\end{align*}
Notice that for any $c > 0$, there exist $\lambda_{j,1}, \lambda_{j,2} > 0$ such that 
\begin{align*}
\langle y, P_{0,j}\, y \rangle \ge ( r_{\kappa} - c ) \|y\|^2 \quad \text{for any $y \in \RR^6$.}
\end{align*}
For each Jordan block of $A_0$ and $B_0$, we apply similar rescaling if necessary to the corresponding components of $a$ and $b$ so that under the change of coordinates $(a, b) \mapsto (p, q)$, (\ref{ave:NF0}) becomes 
\begin{align*} 
\dot p ={}& \epsilon ( P_0 p + P_1(p,q,\zeta) ) \,, \\
\dot q ={}& \epsilon ( Q_0 q + Q_1(p,q,\zeta) ) \,, \\
\dot \zeta ={}& \epsilon ( 1 + Z_1(p,q,\zeta) ) \,,
\end{align*}
with $P_0$ and $Q_0$ now satisfying
\begin{subequations} \label{innerprod:WH:p-q}
\begin{align}
\langle p, P_0 p \rangle \ge{}& \sigma \|p\|^2 \quad \text{for any $p \in \RR^{n_u}$,} \label{innerprod:WH:p} \\
\langle q, Q_0 q \rangle \le{}& -\sigma \|q\|^2 \quad \text{for any $q \in \RR^{n_s}$,} \label{innerprod:WH:q}
\end{align}
\end{subequations}
where $\sigma$ is a positive constant satisfying $0 < \sigma < \min \big\{ \sigma_s, \sigma_u \big\}$. We denote this change of coordinates by $(a, b) = v(p, q)$.

We now obtain a normal form of the full system (\ref{GWH}) near the torus $\To$ in terms of $(p,q,\zeta,\theta)$ as follows:
\begin{align}\label{WHNF}
\begin{split}
\dot p ={}& \epsilon ( P_0 p + P_1(p,q,\zeta) ) + \epsilon^{1+\mu} P_2 (p,q,\zeta,\theta,\epsilon) \,, \\
\dot q ={}& \epsilon ( Q_0 q + Q_1(p,q,\zeta) ) + \epsilon^{1+\mu} Q_2 (p,q,\zeta,\theta,\epsilon) \,, \\
\dot \zeta ={}& \epsilon ( 1 + Z_1(p,q,\zeta) ) + \epsilon^{1+\mu} Z_2 (p,q,\zeta,\theta,\epsilon) \,, \\
\dot \theta ={}& \Theta_0 + \epsilon^{\nu} \Theta_1(p,q,\zeta,\epsilon) + \epsilon^{1+\gamma} \Theta_2 (p,q,\zeta,\theta,\epsilon) \,,
\end{split}
\end{align}
where $\|p\| < \delta_0$ and $\|q\| < \delta_0$ for a certain $\delta_0 > 0$ such that $v(p, q) \in \big\{ (a, b) : \|a\| < \Delta, \|b\| < \Delta \big\}$. In addition to taking $\zeta \in \RR$, we lift $\theta$ to $\RR^m$ in the subsequent analysis. Based on Assumptions \ref{asp:GWH:sm} and \ref{asp:GWH:hyper} and the preceding changes of variables, we can easily verify a set of properties regarding the smoothness, boundedness, and periodicity of the functions on the right-hand side of (\ref{WHNF}). We state these properties in the following lemma while omitting their straightforward verifications.

\begin{lem} \label{lem:WHNF:est}
Let\/ $D_0 := \big\{ (p,q) : \|p\| < \delta_0, \, \|q\| < \delta_0 \big\}$.  
\begin{enumerate}
   \item The constant matrices\/ $P_0$ and\/ $Q_0$ satisfy (\ref{innerprod:WH:p-q}).
   \item The functions\/ $P_1$, $Q_1$, and\/ $Z_1$ are\/ $C^r$ on\/ $D_0 \times \RR$.
   \item The function\/ $\Theta_1$ is\/ $C^r$ with respect to\/ $(p,q,\zeta)$ and continuous with respect to\/ $\epsilon$ on\/ $D_0 \times \RR \times [0, \epsilon_0]$.
   \item The functions\/ $P_2$, $Q_2$, $Z_2$, and\/ $\Theta_2$ are\/ $C^r$ with respect to\/ $(p,q,\zeta,\theta)$ and continuous with respect to\/ $\epsilon$ on\/ $D_0 \times \RR \times \RR^m \times [0, \epsilon_0]$.
   \item There exist positive constants\/ $C_1$ and\/ $C_2$ such that for any\/ $\delta \in [0, \delta_0]$ and any\/ $\|p\| \le \delta$, $\|q\| \le \delta$, and\/ $\zeta \in \RR$,
\begin{align*}
\|P_1\| \le{}& C_1 \delta^2 \,, & 
\|Q_1\| \le{}& C_1 \delta^2 \,, &
\|Z_1\| \le{}& C_1 \delta^2 \,, \\
\|D_{\zeta} P_1\| \le{}& C_1 \delta^2 \,, & 
\|D_{\zeta} Q_1\| \le{}& C_1 \delta^2 \,, &
\|D_{\zeta} Z_1\| \le{}& C_1 \delta^2 \,, \\
\|D_{\xi} P_1\| \le{}& C_2 \delta \,, & 
\|D_{\xi} Q_1\| \le{}& C_2 \delta \,, &
\|D_{\xi} Z_1\| \le{}& C_2 \delta \,, 
\end{align*}
where\/ $\xi = p$ or\/ $q$ in the last three inequalities.
   \item There exist positive constants\/ $C_3$ and\/ $C_4$ such that for any\/ $(p,q) \in D_0$, $\zeta \in \RR$, $\theta \in \RR^m$, and\/ $\epsilon \in [0, \epsilon_0]$,
\begin{align*}
\|P_2\| \le{}& C_3 \,, & 
\|Q_2\| \le{}& C_3 \,, &
\|Z_2\| \le{}& C_3 \,, \\
\|\Theta_1\| \le{}& C_3 \,, &
\|\Theta_2\| \le{}& C_3 \,, \\
\|D_{\xi} P_2\| \le{}& C_4 \,, & 
\|D_{\xi} Q_2\| \le{}& C_4 \,, &
\|D_{\xi} Z_2\| \le{}& C_4 \,, \\
\|D_{\xi} \Theta_1\| \le{}& C_4 \,, &
\|D_{\xi} \Theta_2\| \le{}& C_4 \,,
\end{align*}
where\/ $\xi = p$, $q$, $\zeta$, or\/ $\theta$.
   \item The second to\/ $r$-th derivatives of\/ $P_1$, $P_2$, $Q_1$, $Q_2$, $Z_1$, $Z_2$, $\Theta_1$, and\/ $\Theta_2$ with respect to\/ $(p,q,\zeta)$ or\/ $(p,q,\zeta,\theta)$ are all bounded on the corresponding domains specified in (2), (3), and (4).
   \item The functions\/ $P_1$, $P_2$, $Q_1$, $Q_2$, $Z_1$, $Z_2$, $\Theta_1$, and\/ $\Theta_2$ are\/ $2 \zeta_0$-periodic in\/ $\zeta$ and\/ $2 \pi$-periodic in each component of\/ $\theta$ on the corresponding domains specified in (2), (3), and (4).
\end{enumerate}
\end{lem}

We now show that for every fixed, sufficiently small $\epsilon > 0$, (\ref{WHNF}) has a $C^r$ positively invariant manifold. We introduce another two rescaled variables: $\bq := q/2$ and $\btheta := k^{-1} \theta$ with the scaling factor $k > 0$ to be determined later. By taking $z = (\bq,\zeta,\btheta)$ and treating $p$ as ``$a$'', we organize the rescaled system into the form of (\ref{sys}) as follows:
\begin{align} \label{WHNF:res}
\begin{split} 
\dot p ={}& f(p,z) := 
\epsilon ( P_0 p + P_1(p,2\bq,\zeta) ) + \epsilon^{1+\mu} P_2 (p,2\bq,\zeta,k\btheta,\epsilon) \,, \\
\dot z ={}& g(p,z) := \begin{pmatrix} 
\epsilon ( Q_0 \bq + \tfrac{1}{2} Q_1(p,2\bq,\zeta) ) + \tfrac{1}{2} \epsilon^{1+\mu} Q_2 (p,2\bq,\zeta,k\btheta,\epsilon) \\
\epsilon ( 1 + Z_1(p,2\bq,\zeta) ) + \epsilon^{1+\mu} Z_2 (p,2\bq,\zeta,k\btheta,\epsilon) \\
\tfrac{1}{k} \Theta_0 + \tfrac{1}{k} \epsilon^{\nu} \Theta_1(p,2\bq,\zeta,\epsilon) + \tfrac{1}{k} \epsilon^{1+\gamma} \Theta_2 (p,2\bq,\zeta,k\btheta,\epsilon)
\end{pmatrix} .
\end{split}
\end{align}

We shall restrict $(p,z) = (p,\bq,\zeta,\btheta)$ to the domain 
\begin{equation*}
\Ud := \big\{ \|p\| < \delta \big\} \times \big\{ \|\bq\| < \tfrac{\delta}{2} \big\} \times \RR \times \RR^m \,,
\end{equation*}
where $\delta \in (0, \delta_0]$ is to be determined. Obviously, $\Ud$ satisfies Hypothesis \ref{hypo:U}. As we have done for the previous two examples, we define $\Gaekd \subset \Ud$ to be the set of points that stay in $\Ud$ forever in forward time along solution trajectories of (\ref{WHNF:res}), that is, 
\begin{align}
\Gaekd : = \big\{ (p,z) : \Phi_{\epsilon, k}(t,(p,z)) \in \Ud \text{ for all } t \ge 0 \big\} \,, \label{def:Gaekd}
\end{align}
where $\Phi_{\epsilon, k}(t,(p,z))$ is the flow of (\ref{WHNF:res}). By this definition, $\Gaekd$ is the largest positively invariant subset of $\Ud$ under the flow $\Phi_{\epsilon, k}$. Then we need to prove that $\Gaekd$ satisfies 
\begin{equation} \label{proj:Gamma:ex3}
\Pibot(\Gaekd) = \Pibot(\Ud) = \big\{ \|\bq\| < \tfrac{\delta}{2} \big\} \times \RR \times \RR^m \,,
\end{equation} 
where $\Pibot$ is the projection onto the $z$-coordinate. Again, we will achieve this using the {\Wa} theorem. However, since $p$ sits in $\RR^{n_u}$, some arguments are slightly different. 

\begin{lem} \label{lem:WH:Gamma:Wazi}
There exist\/ $\epsilon_1 \in (0, \epsilon_0]$ and\/ $\delta_1 \in (0, \delta_0]$ satisfying\/ $2 \frac{C_3}{\sigma} \epsilon_1^{\mu} < \delta_1$ such that for (\ref{WHNF:res}) with any\/ $\epsilon \in (0, \epsilon_1]$ and\/ $\Ud$ with any\/ $\delta \in [2 \frac{C_3}{\sigma} \epsilon^{\mu}, \delta_1]$, the positively invariant set\/ $\Gaekd \subset \Ud$ satisfies (\ref{proj:Gamma:ex3}). 
\end{lem}

\begin{proof}
We partition the boundary of $\Ud$ into two subsets $S_1$ and $S_2$ as follows:
\begin{align*}
S_1 :={}& \big\{ \|p\| = \delta \big\} \times \big\{ \|\bq\| \le \tfrac{\delta}{2} \big\} \times \RR \times \RR^m \,, \\
S_2 :={}& \big\{ \|p\| < \delta \big\} \times \big\{ \|\bq\| = \tfrac{\delta}{2} \big\} \times \RR \times \RR^m \,.
\end{align*}
Using (\ref{innerprod:WH:p}) and the estimates in part (5) and part (6) of Lemma \ref{lem:WHNF:est}, we obtain that for any point in $S_1$, 
\begin{align*}
\langle p , \dot p \rangle 
={}& \langle p , \epsilon ( P_0 p + P_1(p,2\bq,\zeta) ) \rangle + \langle p , \epsilon^{1+\mu} P_2 (p,2\bq,\zeta,k\btheta,\epsilon) \rangle \\
\ge{}& \epsilon \langle p , P_0 p \rangle - \epsilon C_1 \delta^2 \|p\| - \epsilon^{1+\mu} C_3 \|p\| \\
\ge{}& \epsilon \sigma \|p\|^2 - \epsilon C_1 \delta^2 \|p\| - \epsilon^{1+\mu} C_3 \|p\| \\
={}& \epsilon \delta ( \sigma \delta - C_1 \delta^2 - \epsilon^{\mu} C_3 ) \,. 
\end{align*}
Similarly, for points in $S_2$, we have 
\begin{align*}
\langle \bq , \dot{\bq} \rangle 
={}& \langle \bq , \epsilon ( Q_0 \bq + \tfrac{1}{2} Q_1(p,2\bq,\zeta) ) \rangle + \langle \bq , \tfrac{1}{2} \epsilon^{1+\mu} Q_2 (p,2\bq,\zeta,k\btheta,\epsilon) \rangle \\
\le{}& \epsilon \langle \bq , Q_0 \bq \rangle + \tfrac{1}{2} \epsilon C_1 \delta^2 \|\bq\| + \tfrac{1}{2} \epsilon^{1+\mu} C_3 \|\bq\| \\
\le{}& - \epsilon \sigma \|\bq\|^2 + \tfrac{1}{2} \epsilon C_1 \delta^2 \|\bq\| + \tfrac{1}{2} \epsilon^{1+\mu} C_3 \|\bq\| \\
={}& - \epsilon \tfrac{\delta}{4} ( \sigma \delta - C_1 \delta^2 - \epsilon^{\mu} C_3 ) \,. 
\end{align*}
Since for $\mu > 0$, $\epsilon^{\mu} \rightarrow 0$ as $\epsilon \rightarrow 0^+$, we can take $\epsilon_1 \in (0, \epsilon_0]$ and $\delta_1 \in (0, \delta_0]$ such that $2 \frac{C_3}{\sigma} \epsilon_1^{\mu} < \delta_1 \le \frac{\sigma}{2C_1}$. Then for any $\epsilon \in (0, \epsilon_1]$ and any $\delta \in [2 \frac{C_3}{\sigma} \epsilon^{\mu}, \delta_1]$,
\begin{align*}
\sigma \delta - C_1 \delta^2 - \epsilon^{\mu} C_3 \ge{}& \sigma (2 \tfrac{C_3}{\sigma} \epsilon^{\mu}) - C_1 (2 \tfrac{C_3}{\sigma} \epsilon^{\mu})^2 - \epsilon^{\mu} C_3 \\
={}& C_3 \epsilon^{\mu} - 4 C_1 C_3^2 \tfrac{1}{\sigma^2} \epsilon^{2\mu} > 0 \,,
\end{align*}
where the last inequality can be guaranteed by choosing $\epsilon_1$ sufficiently small. It follows that $\frac{d}{dt} \|p\|^2 > 0$ for any point in $S_1$ and $\frac{d}{dt} \|\bq\|^2 < 0$ for any point in $S_2$. Then solution trajectories of (\ref{WHNF:res}) leave and enter $\closure(\Ud)$ through points in $S_1$ and $S_2$, respectively. In particular, $S_1$ is the set of points through which trajectories leave $\closure(\Ud)$ immediately in forward time. Then it can be easily verified that $\closure(\Ud)$ is a {\Wa} set. Let $\Ud^0 \subseteq \closure(\Ud)$ be the set of points that do not stay in $\closure(\Ud)$ forever in forward time. By the {\Wa} theorem, there exists a continuous function $\CR : \Ud^0 \times [0,1] \rightarrow \Ud^0$ such that $\CR$ is a strong deformation retraction of $\Ud^0$ onto $S_1$. 

Suppose that for a certain $z^* = (\bq^*,\zeta^*,\btheta^*) \in \big\{ \|\bq\| < \tfrac{\delta}{2} \big\} \times \RR \times \RR$,
\begin{equation} \label{Gamma:zstr}
\Gaekd \textcap \big\{ (p, z^*): \|p\| < \delta \big\} = \emptyset \,,
\end{equation}
where the second set on the left is just the section of $\Ud$ at $z^*$. Since $\big\{ (p, z^*): \|p\| = \delta \big\} \subset S_1 \subset \Ud^0$, (\ref{Gamma:zstr}) implies that the set $\big\{ (p, z^*): \|p\| \le \delta \big\}$ is contained inside the domain $\Ud^0$ of the continuous function $\CR(\,\cdot\,,1) : \Ud^0 \rightarrow S_1$. This allows us to construct a continuous function that maps the closed $n_u$-ball $\big\{ \|p\| \le \delta \big\}$ to its boundary $\big\{ \|p\| = \delta \big\}$:
\begin{equation*}
p \mapsto \Pi \circ \CR((p, z^*), 1) \,,
\end{equation*}
where $\Pi$ is the projection onto the $p$-coordinate. Since $\CR$ is a strong deformation retraction of $\Ud^0$ onto $S_1$, $\CR ((p, z^*), 1) = (p, z^*)$ for any $(p, z^*) \in S_1$. Thus, the above continuous function has the property that for any $p$ with $\|p\| = \delta$,
\begin{equation*}
\Pi \circ \CR((p, z^*), 1) = \Pi (p, z^*) = p \,.
\end{equation*}
The existence of such a continuous function contradicts the fact that there is no retraction that maps a closed $n$-ball onto its boundary (i.e., an $(n\!-\!1)$-sphere). Thus, $\Gaekd \textcap \big\{ (p, z): \|p\| < \delta \big\} \neq \emptyset$ for any $z \in \big\{ \|\bq\| < \tfrac{\delta}{2} \big\} \times \RR \times \RR^m$. This proves (\ref{proj:Gamma:ex3}).
\end{proof}

Next, using (1), (5), and (6) of Lemma \ref{lem:WHNF:est}, we obtain the following estimates regarding $D_p f(p,z)$, $D_z g(p,z)$, $D_z f(p,z)$, and $D_p g(p,z)$ for any $(p,z) \in \Ud$.

\begin{enumerate}

\item For any $p' \in \RR^{n_u}$,
\begin{align*}
\langle p' , D_p f\, p' \rangle 
={}& \langle p', ( \epsilon P_0 + \epsilon D_p P_1 + \epsilon^{1+\mu} D_p P_2 )\, p' \rangle \\
\ge{}& \epsilon \sigma \|p'\|^2 - \epsilon C_2 \delta \|p'\|^2 - \epsilon^{1+\mu} C_4 \|p'\|^2 \\
={}& \au \|p'\|^2 \,,
\end{align*}
where 
\begin{equation*}
\au := \epsilon ( \sigma - C_2 \delta ) - \epsilon^{1+\mu} C_4 \,.
\end{equation*}

\item For any $z' = (\bq',\zeta',\btheta') \in \RR^{n_s} \times \RR \times \RR^m$,
\begin{align*}
\langle z' , D_z g\, z' \rangle 
={}& \langle \bq' , ( \epsilon Q_0 + \epsilon D_q Q_1 + \epsilon^{1+\mu} D_q Q_2 )\, \bq' \rangle 
+ \langle \bq' , \tfrac{1}{2} k \epsilon^{1+\mu} D_{\theta} Q_2\, \btheta' \rangle \\
& \qquad + \langle \bq' , ( \tfrac{1}{2}\epsilon D_{\zeta} Q_1 + \tfrac{1}{2} \epsilon^{1+\mu} D_{\zeta} Q_2 )\, \zeta' \rangle \\
& + \langle \zeta' , ( 2\epsilon D_q Z_1 + 2 \epsilon^{1+\mu} D_q Z_2 )\, \bq' \rangle 
+ \langle \zeta' , k \epsilon^{1+\mu} D_{\theta} Z_2\, \btheta' \rangle \\
& \qquad + \langle \zeta' , ( \epsilon D_{\zeta} Z_1 + \epsilon^{1+\mu} D_{\zeta} Z_2 )\, \zeta' \rangle \\
& + \langle \btheta' , ( 2\tfrac{1}{k} \epsilon^{\nu} D_q \Theta_1 + 2\tfrac{1}{k} \epsilon^{1+\gamma} D_q \Theta_2 )\, \bq' \rangle 
+ \langle \btheta' , \epsilon^{1+\gamma} D_{\theta} \Theta_2\, \btheta' \rangle \\
& \qquad + \langle \btheta' , ( \tfrac{1}{k} \epsilon^{\nu} D_{\zeta} \Theta_1 + \tfrac{1}{k} \epsilon^{1+\gamma} D_{\zeta} \Theta_2 )\, \zeta' \rangle \\
\le{}& - \epsilon \sigma \|\bq'\|^2 + (\epsilon C_2 \delta + \epsilon^{1+\mu} C_4) \|\bq'\|^2
+ \tfrac{1}{2} k \epsilon^{1+\mu} C_4 \|\bq'\|\|\btheta'\| \\
& \qquad + ( \tfrac{1}{2}\epsilon C_1 \delta^2 + \tfrac{1}{2}\epsilon^{1+\mu} C_4 ) \|\bq'\| |\zeta'| \\
& + ( 2\epsilon C_2 \delta + 2\epsilon^{1+\mu} C_4) |\zeta'| \|\bq'\|
+ k \epsilon^{1+\mu} C_4 |\zeta'| \|\btheta'\| \\
& \qquad + (\epsilon C_1 \delta^2 + \epsilon^{1+\mu} C_4 ) |\zeta'|^2 \\
& + ( 2\tfrac{1}{k} \epsilon^{\nu} C_4 + 2\tfrac{1}{k} \epsilon^{1+\gamma} C_4) \|\btheta'\|\|\bq'\|
+ \epsilon^{1+\gamma} C_4 \|\btheta'\|^2 \\
& \qquad + (\tfrac{1}{k} \epsilon^{\nu} C_4 + \tfrac{1}{k} \epsilon^{1+\gamma} C_4) \|\btheta'\| |\zeta'| \\
\le{}& \ell_u \|z'\|^2 \,,
\end{align*}
where 
\begin{align*}
\ell_u :={}&
\epsilon ( 3C_2\delta + \tfrac{3}{2}C_1\delta^2 ) + \epsilon^{1+\mu} \tfrac{9}{2} C_4 + \epsilon^{1+\gamma} C_4 \\
& + k \epsilon^{1+\mu} \tfrac{3}{2} C_4 + \tfrac{1}{k} ( \epsilon^{\nu} 3C_4 + \epsilon^{1+\gamma} 3C_4 ) \,.
\end{align*}
To obtain the last inequality, we have discarded the negative term $- \epsilon \sigma \|\bq'\|^2$ and then used $\|\bq'\| \le \|z'\|$, $|\zeta'| \le \|z'\|$, and $\|\btheta'\| \le \|z'\|$.

\item 
\begin{align*}
\|D_z f\| 
\le{}& ( \| 2\epsilon D_q P_1 \| + \| 2\epsilon^{1+\mu} D_q P_2 \| )
+ ( \| \epsilon D_{\zeta} P_1 \| + \| \epsilon^{1+\mu} D_{\zeta} P_2 \| ) \\
& \qquad + \| k \epsilon^{1+\mu} D_{\theta} P_2 \| \\
\le{}& ( 2\epsilon C_2\delta + 2\epsilon^{1+\mu} C_4 ) + ( \epsilon C_1\delta^2 + \epsilon^{1+\mu} C_4 ) + k \epsilon^{1+\mu} C_4 \\
={}& \Lfz \,,
\end{align*}
where
\begin{align*}
\Lfz := \epsilon ( 2C_2\delta + C_1\delta^2 ) + \epsilon^{1+\mu} 3C_4 + k \epsilon^{1+\mu} C_4 \,.
\end{align*}

\item 
\begin{align*}
\|D_p g\| 
\le{}& ( \| \tfrac{1}{2}\epsilon D_p Q_1 \| + \| \tfrac{1}{2}\epsilon^{1+\mu} D_p Q_2 \| )
+ ( \| \epsilon D_p Z_1 \| + \| \epsilon^{1+\mu} D_p Z_2 \| ) \\
& \qquad + ( \| \tfrac{1}{k} \epsilon^{\nu} D_p \Theta_1 \| + \| \tfrac{1}{k}\epsilon^{1+\gamma} D_p \Theta_2 \| ) \\
\le{}& ( \tfrac{1}{2}\epsilon C_2\delta + \tfrac{1}{2}\epsilon^{1+\mu} C_4 ) + ( \epsilon C_2\delta  + \epsilon^{1+\mu} C_4 ) + ( \tfrac{1}{k} \epsilon^{\nu} C_4 + \tfrac{1}{k}\epsilon^{1+\gamma} C_4 ) \\
={}& \Lgp \,,
\end{align*}
where
\begin{align*}
\Lgp := \epsilon \tfrac{3}{2} C_2\delta + \epsilon^{1+\mu} \tfrac{3}{2}C_4 + \tfrac{1}{k} ( \epsilon^{\nu} C_4 + \epsilon^{1+\gamma} C_4 ) \,.
\end{align*}

\end{enumerate}

In order to verify Hypothesis \ref{hypo:C1} and Hypothesis \refHypoCr for (\ref{WHNF:res}) on $\Ud$, we define an auxiliary function $\CK : (0, \epsilon_0] \times (0, \infty) \times (0, \delta_0] \rightarrow \RR$ as follows:
\begin{align*}
\CK(\epsilon,k,\delta) :={}& \au - ( r \ell_u + (r+1)\Lgp + \Lfz ) \\
={}& \epsilon \sigma - \epsilon ( K_1 \delta + K_2 \delta^2 + \epsilon^{\mu} K_3 + \epsilon^{\gamma} K_4 + k \epsilon^{\mu} K_5 + \tfrac{1}{k} ( \epsilon^{\nu-1} K_6 + \epsilon^{\gamma} K_7 ) ) \,,
\end{align*}
where $r$ is the degree of smoothness referred to in (2) and (3) of Lemma \ref{lem:WHNF:est} and the constants $K_i$ are defined as follows: 
\begin{align*}
K_1 :={}& \tfrac{9}{2}(r+1)C_2 \,, \\
K_2 :={}& (\tfrac{3}{2}r+1)C_1 \,, \\
K_3 :={}& (6r+\tfrac{11}{2})C_4 \,, \\
K_4 :={}& r C_4 \,, \\
K_5 :={}& (\tfrac{3}{2}r+1) C_4 \,, \\
K_6 :={}& (4r+1)C_4 \,, \\
K_7 :={}& (4r+1)C_4 \,.
\end{align*}

\begin{lem} \label{lem:WH:CK}
There exist\/ $\epsilon_2 \in (0, \epsilon_0]$ and\/ $\delta_2 \in (0, \delta_0]$ satisfying\/ $2 \frac{C_3}{\sigma} \epsilon_2^{\mu} < \delta_2$ such that\/ $\CK(\epsilon,k_{\epsilon},\delta) > 0$ for any\/ $\epsilon \in (0, \epsilon_2]$, any\/ $\delta \in [2 \frac{C_3}{\sigma} \epsilon^{\mu}, \delta_2]$, and an appropriately chosen\/ $k_{\epsilon}$ that depends on\/ $\epsilon$.
\end{lem}

\begin{proof}
Note that for any $\epsilon > 0$,
\begin{align*}
\inf_{k \in (0, \infty)} \big\{ k \epsilon^{\mu} K_5 + \tfrac{1}{k} ( \epsilon^{\nu-1} K_6 + \epsilon^{\gamma} K_7 ) \big\} 
= \big( \epsilon^{\mu+\nu-1} K_5 K_6 + \epsilon^{\mu+\gamma} K_5 K_7 \big)^{\frac{1}{2}} \,.
\end{align*}
In particular, for every $\epsilon > 0$, there exists a $k_{\epsilon} > 0$ such that
\begin{align} \label{ke}
k_{\epsilon} \epsilon^{\mu} K_5 + \tfrac{1}{k_{\epsilon}} ( \epsilon^{\nu-1} K_6 + \epsilon^{\gamma} K_7 ) 
\le \big( \epsilon^{\mu+\nu-1} K_5 K_6 + \epsilon^{\mu+\gamma} K_5 K_7 \big)^{\frac{1}{2}} + \tfrac{1}{3} \sigma \,.
\end{align}

By (\ref{indices}), we have $\epsilon^{\mu} \rightarrow 0$, $\epsilon^{\gamma} \rightarrow 0$, $\epsilon^{\mu+\gamma} \rightarrow 0$, and $\epsilon^{\mu+\nu-1} \rightarrow 0$ for $\epsilon \rightarrow 0^+$. Thus, there exist $\epsilon_2 \in (0, \epsilon_0]$ and $\delta_2 \in (0, \delta_0]$ satisfying $2 \frac{C_3}{\sigma} \epsilon_2^{\mu} < \delta_2$ such that for any $\epsilon \in (0, \epsilon_2]$ and any $\delta \in [2 \frac{C_3}{\sigma} \epsilon^{\mu}, \delta_2]$, 
\begin{align*}
\tfrac{2}{3} \sigma \ge{}&
K_1 \delta + K_2 \delta^2 + \epsilon^{\mu} K_3 + \epsilon^{\gamma} K_4 + \big( \epsilon^{\mu+\nu-1} K_5 K_6 + \epsilon^{\mu+\gamma} K_5 K_7 \big)^{\frac{1}{2}} + \tfrac{1}{3} \sigma \\
\ge{}& K_1 \delta + K_2 \delta^2 + \epsilon^{\mu} K_3 + \epsilon^{\gamma} K_4 + k_{\epsilon} \epsilon^{\mu} K_5 + \tfrac{1}{k_{\epsilon}} ( \epsilon^{\nu-1} K_6 + \epsilon^{\gamma} K_7 ) \,.
\end{align*}
It follows that $\CK(\epsilon,k_{\epsilon},\delta) \ge \tfrac{1}{3} \epsilon \sigma > 0$.
\end{proof}

Lemma \ref{lem:WH:CK} guarantees that both (\ref{ineq:C1}) and (\refineqCr) hold for (\ref{WHNF:res}) on $\Ud$. Thus, by combining Lemma \ref{lem:WH:Gamma:Wazi} and Lemma \ref{lem:WH:CK}, we establish the existence and the $C^r$ smoothness of a positively invariant manifold for (\ref{WHNF:res}) as stated in the following lemma. 

\begin{lem} \label{lem:WH:Gamma}
There exist\/ $\epsilon^* \in (0, \epsilon_0]$ and\/ $\delta^* \in (0, \delta_0]$ satisfying\/ $2 \frac{C_3}{\sigma} {\epsilon^*}^{\mu} < \delta^*$ such that for (\ref{WHNF:res}) with any\/ $\epsilon \in (0, \epsilon^*]$ and any\/ $k = k_{\epsilon}$ that satisfies (\ref{ke}), the positively invariant set\/ $\Gaekedst \subset \Udst$ is the graph of a\/ $C^r$ function\/ $h( \cdot, \cdot, \cdot; \epsilon, k_{\epsilon} ) : \Pibot(\Udst) \rightarrow \big\{ \|p\| < \delta^* \big\}$
that satisfies
\begin{equation} \label{WH:h:Lip}
\| h( \bq_2, \zeta_2, \btheta_2 ; \epsilon, k_{\epsilon} ) - h( \bq_1, \zeta_1, \btheta_1 ; \epsilon, k_{\epsilon} ) \| 
\le \| \bq_2 - \bq_1 \| + | \zeta_2 - \zeta_1 | + \| \btheta_2 - \btheta_1 \|  
\end{equation}
for any\/ $(\bq_1, \zeta_1, \btheta_1)$ and\/ $(\bq_2, \zeta_2, \btheta_2) \in \Pibot(\Udst)$. Furthermore, for each\/ $\delta \in [2 \frac{C_3}{\sigma} \epsilon^{\mu}, \delta^*]$, the positively invariant set\/ $\Gaeked = \Gaekedst \textcap \Ud$. 
\end{lem}

\begin{proof}
Take $\epsilon^* = \min \big\{ \epsilon_1, \epsilon_2 \big\}$ and $\delta^* = \min \big\{ \delta_1, \delta_2 \big\}$. Clearly, $2 \frac{C_3}{\sigma} {\epsilon^*}^{\mu} < \delta^*$. Thus, for (\ref{WHNF:res}) with any $\epsilon \in (0, \epsilon^*]$ and any $k = k_{\epsilon}$ that satisfies (\ref{ke}), Lemma \ref{lem:WH:Gamma:Wazi} and Lemma \ref{lem:WH:CK} guarantee that both Theorem \ref{thm1} and Theorem \ref{thm3} are applicable on the domain $\Ud = \big\{ \|p\| < \delta \big\} \times \big\{ \|\bq\| < \tfrac{\delta}{2} \big\} \times \RR \times \RR^m$ with any $\delta \in [2 \frac{C_3}{\sigma} \epsilon^{\mu}, \delta^*]$. In particular, for $\delta = \delta^*$, the positively invariant set $\Gaekedst \subset \Udst$ is the graph of a $C^r$ function 
\begin{align*}
h( \cdot, \cdot, \cdot ; \epsilon, k_{\epsilon} ) : 
\Pibot(\Udst) = \big\{ \|\bq\| < \tfrac{\delta^*}{2} \big\} \times \RR \times \RR^m \rightarrow \big\{ \|p\| < \delta^* \big\}
\end{align*}
that satisfies (\ref{WH:h:Lip}) according to Theorem \ref{thm1} or Theorem \ref{thm3}. Furthermore, for each $\delta \in [2 \frac{C_3}{\sigma} \epsilon^{\mu}, \delta^*]$, $\Pibot(\Gaeked) = \Pibot(\Ud)$ by Lemma \ref{lem:WH:Gamma:Wazi}, and $\Gaeked \subseteq \Gaekedst \textcap \Ud$ by the definition (\ref{def:Gaekd}). Thus $\Gaeked = \Gaekedst \textcap \Ud$.
\end{proof}

We also have the following facts regarding the function $h$.

\begin{enumerate}

\item For any $\epsilon \in (0, \epsilon^*]$ and any $k_{\epsilon}$ and $k'_{\epsilon}$ that both satisfy (\ref{ke}), 
\begin{align} \label{ex3:res:equiv}
h( \bq, \zeta, \tfrac{\theta}{k_{\epsilon}} ; \epsilon, k_{\epsilon} ) = 
h( \bq, \zeta, \tfrac{\theta}{k'_{\epsilon}} ; \epsilon, k'_{\epsilon} )
\end{align}
for all $(\bq, \zeta, \btheta) \in \big\{ \|\bq\| < \tfrac{\delta^*}{2} \big\} \times \RR \times \RR^m$.

\item For any $\epsilon \in (0, \epsilon^*]$ and any $k_{\epsilon}$ that satisfies (\ref{ke}), $h( \bq, \zeta, \btheta ; \epsilon, k_{\epsilon} )$ is $2\zeta_0$-periodic in $\zeta$ and $\frac{2\pi}{k_{\epsilon}}$-periodic in each component of $\btheta$.

\item For any $\epsilon \in (0, \epsilon^*]$ and any $k_{\epsilon}$ that satisfies (\ref{ke}), 
\begin{align*} 
\| h( \bq, \zeta, \btheta ; \epsilon, k_{\epsilon} ) \| \le 2 \tfrac{C_3}{\sigma} \epsilon^{\mu}
\end{align*}
for all $(\bq, \zeta, \btheta) \in \big\{ \|\bq\| < \tfrac{C_3}{\sigma} \epsilon^{\mu} \big\} \times \RR \times \RR^m$.

\end{enumerate}

The first two properties should be familiar by now, and they follow from the definition (\ref{def:Gaekd}) and the same arguments as given in Subsections \ref{ex1} and \ref{ex2} for the first two examples. The third property is a simple consequence of the fact that $\Gaeked = \Gaekedst \textcap \Ud$ for $\delta = 2 \frac{C_3}{\sigma} \epsilon^{\mu}$.

By (\ref{ex3:res:equiv}), we can define a $C^r$ function $\CP(\cdot, \cdot, \cdot ; \epsilon) : \big\{ \|q\| < \delta^* \big\} \times \RR \times \RR^m \rightarrow \big\{ \|p\| < \delta^* \big\}$ for each $\epsilon \in (0, \epsilon^*]$ as follows:
\begin{align*}
\CP( q, \zeta, \theta ; \epsilon) := h( \tfrac{q}{2}, \zeta, \tfrac{\theta}{k_{\epsilon}} ; \epsilon, k_{\epsilon} ) \,,
\end{align*}
where $k_{\epsilon}$ can be taken to be any value that satisfies (\ref{ke}). Then it follows immediately from (\ref{WH:h:Lip}) that 
\begin{equation} \label{WH:CP:Lip}
\| \CP( q_2, \zeta_2, \theta_2 ; \epsilon) - \CP( q_1, \zeta_1, \theta_1 ; \epsilon) \|
\le \tfrac{1}{2} \| q_2 - q_1 \| + | \zeta_2 - \zeta_1 | + \tfrac{1}{k_{\epsilon}} \| \theta_2 - \theta_1 \|
\end{equation}
for any $(q_1, \zeta_1, \theta_1)$ and $(q_2, \zeta_2, \theta_2) \in \big\{ \|q\| < \delta^* \big\} \times \RR \times \RR^m$. In addition, by the properties (2) and (3) above, $\CP(q, \zeta, \theta ; \epsilon)$ is $2\zeta_0$-periodic in $\zeta$ and $2\pi$-periodic in each component of $\theta$, and 
\begin{equation} \label{ex3:CP:e-bound}
\| \CP( q, \zeta, \theta ; \epsilon ) \| \le 2 \tfrac{C_3}{\sigma} \epsilon^{\mu}
\end{equation}
for all $(q, \zeta, \theta) \in \big\{ \|q\| < 2 \tfrac{C_3}{\sigma} \epsilon^{\mu} \big\} \times \RR \times \RR^m$.

We now return to the normal form (\ref{WHNF}). The results stated in the next lemma are now obvious.

\begin{lem} \label{lem:WH:W+}
For any\/ $\epsilon \in (0, \epsilon^*]$, (\ref{WHNF}) has a\/ $C^r$ positively invariant manifold\/ $W^+_{\epsilon}$, which is the graph of the\/ $C^r$ function\/ $\CP(\cdot, \cdot, \cdot ; \epsilon) : \big\{ \|q\| < \delta^* \big\} \times \RR \times \RR^m \rightarrow \big\{ \|p\| < \delta^* \big\}$. In addition, $W^+_{\epsilon}$ is the largest positively invariant subset of\/ $\big\{ \|p\| < \delta^* \big\} \times \big\{ \|q\| < \delta^* \big\} \times \RR \times \RR^m$. 
\end{lem}

Notice that all the analysis starting from the formulation of (\ref{WHNF:res}) up to Lemma \ref{lem:WH:W+} can be adapted (with only minor modifications) for the time reversal of (\ref{WHNF}). In particular, we let $\bp := p/2$ while keeping $\btheta = k^{-1} \theta$, and then form a system in the spirit of (\ref{WHNF:res}) but using the time reversal of (\ref{WHNF}) and treating $q$ and $(\bp,\zeta,\btheta)$ as ``$a$'' and ``$z$'' of (\ref{sys}), respectively. Then, by following all the previous steps, we construct a $C^r$ function $\CQ(\cdot, \cdot, \cdot ; \epsilon) : \big\{ \|p\| < \delta^* \big\} \times \RR \times \RR^m \rightarrow \big\{ \|q\| < \delta^* \big\}$ whose properties are completely analogous to those of the function $\CP(\cdot, \cdot, \cdot ; \epsilon)$, i.e.: 
\begin{equation} \label{WH:CQ:Lip}
\| \CQ( p_2, \zeta_2, \theta_2 ; \epsilon) - \CQ( p_1, \zeta_1, \theta_1 ; \epsilon) \|
\le \tfrac{1}{2} \| p_2 - p_1 \| + | \zeta_2 - \zeta_1 | + \tfrac{1}{k_{\epsilon}} \| \theta_2 - \theta_1 \| 
\end{equation}
for any $(p_1, \zeta_1, \theta_1)$ and $(p_2, \zeta_2, \theta_2) \in \big\{ \|p\| < \delta^* \big\} \times \RR \times \RR^m$; $\CQ(p, \zeta, \theta ; \epsilon)$ is $2\zeta_0$-periodic in $\zeta$ and $2\pi$-periodic in each component of $\theta$; and 
\begin{align} \label{ex3:CQ:e-bound}
\| \CQ( p, \zeta, \theta ; \epsilon ) \| \le 2 \tfrac{C_3}{\sigma} \epsilon^{\mu}
\end{align}
for all $(p, \zeta, \theta) \in \big\{ \|p\| < 2 \tfrac{C_3}{\sigma} \epsilon^{\mu} \big\} \times \RR \times \RR^m$. Then we establish the existence and the $C^r$ smoothness of a negatively invariant manifold for (\ref{WHNF}) as stated in the next lemma.

\begin{lem} \label{lem:WH:W-}
For any\/ $\epsilon \in (0, \epsilon^*]$, (\ref{WHNF}) has a\/ $C^r$ negatively invariant manifold\/ $W^-_{\epsilon}$, which is the graph of the\/ $C^r$ function\/ $\CQ(\cdot, \cdot, \cdot ; \epsilon) : \big\{ \|p\| < \delta^* \big\} \times \RR \times \RR^m \rightarrow \big\{ \|q\| < \delta^* \big\}$. In addition, $W^-_{\epsilon}$ is the largest negatively invariant subset of\/ $\big\{ \|p\| < \delta^* \big\} \times \big\{ \|q\| < \delta^* \big\} \times \RR \times \RR^m$.
\end{lem}

Define $\Me := W^+_{\epsilon} \textcap W^-_{\epsilon}$ for each $\epsilon \in (0, \epsilon^*]$. Clearly, $\Me$ is invariant (i.e., in both forward time and backward time) under the flow of (\ref{WHNF}) with the corresponding $\epsilon$. Below we show that $\Me$ is in fact the graph of a $C^r$ function that maps $\RR \times \RR^m$ into $\big\{ \|p\| < 2 \tfrac{C_3}{\sigma} \epsilon^{\mu} \big\} \times \big\{ \|q\| < 2 \tfrac{C_3}{\sigma} \epsilon^{\mu} \big\}$ and thus a $C^r$ manifold.

\begin{lem} \label{lem:WH:Me}
For any\/ $\epsilon \in (0, \epsilon^*]$, (\ref{WHNF}) has a\/ $C^r$ invariant manifold\/ $\Me = W^+_{\epsilon} \textcap W^-_{\epsilon}$, which is the graph of a\/ $C^r$ function\/ $\rho(\cdot, \cdot ; \epsilon) : \RR \times \RR^m \rightarrow \big\{ \|p\| \le 2 \tfrac{C_3}{\sigma} \epsilon^{\mu} \big\} \times \big\{ \|q\| \le 2 \tfrac{C_3}{\sigma} \epsilon^{\mu} \big\}$. In addition, $\Me$ is the largest invariant subset of\/ $\big\{ \|p\| < \delta^* \big\} \times \big\{ \|q\| < \delta^* \big\} \times \RR \times \RR^m$.
\end{lem}

\begin{proof} 
We consider an arbitrary, fixed $\epsilon \in (0, \epsilon^*]$ throughout this proof. Define the map $\CF_{\epsilon} : \big\{ \|p\| \le \delta^* \big\} \times \big\{ \|q\| \le \delta^* \big\} \times \RR \times \RR^m \rightarrow \big\{ \|p\| \le \delta^* \big\} \times \big\{ \|q\| \le \delta^* \big\}$ as follows:
\begin{align*} 
\CF_{\epsilon}(p, q, \zeta, \theta) := 
( \bar\CP(q, \zeta, \theta ; \epsilon), \bar\CQ(p, \zeta, \theta ; \epsilon) ) \,,
\end{align*}
where $\bar\CP(\cdot, \cdot, \cdot ; \epsilon)$ and $\bar\CQ(\cdot, \cdot, \cdot ; \epsilon)$ are the continuous extensions of $\CP(\cdot, \cdot, \cdot ; \epsilon)$ and $\CQ(\cdot, \cdot, \cdot ; \epsilon)$ onto $\{ \|q\| \le \delta^* \big\} \times \RR \times \RR^m$ and $\{ \|p\| \le \delta^* \big\} \times \RR \times \RR^m$, respectively. 

Consider the norm $\|(p,q)\|_1 := \|p\| + \|q\|$. By (\ref{WH:CP:Lip}) and (\ref{WH:CQ:Lip}), we have that for any $(\zeta, \theta) \in \RR \times \RR^m$ and any $(p_1, q_1)$ and $(p_2, q_2) \in \big\{ \|p\| \le \delta^* \big\} \times \big\{ \|q\| \le \delta^* \big\}$,
\begin{align*}
& \| \CF_{\epsilon}(p_2, q_2, \zeta, \theta) - \CF_{\epsilon}(p_1, q_1, \zeta, \theta) \|_1 \\
& \hspace{0.3in} = \| \bar\CP( q_2, \zeta, \theta ; \epsilon) - \bar\CP( q_1, \zeta, \theta ; \epsilon) \| + \| \bar\CQ( p_2, \zeta, \theta ; \epsilon) - \bar\CQ( p_1, \zeta, \theta ; \epsilon) \| \\
& \hspace{0.3in} \le \tfrac{1}{2} \| q_2 - q_1 \| + \tfrac{1}{2} \| p_2 - p_1 \| \\
& \hspace{0.3in} = \tfrac{1}{2} \| (p_2 - p_1, q_2 - q_1) \|_1 \,.
\end{align*}
Thus, for each $(\zeta, \theta) \in \RR \times \RR^m$, $\CF_{\epsilon}(\cdot, \cdot, \zeta, \theta)$ is a contraction on $\big\{ \|p\| \le \delta^* \big\} \times \big\{ \|q\| \le \delta^* \big\}$ under the norm $\|\cdot\|_1$. Then we can define a function $\rho(\cdot, \cdot ; \epsilon) : \RR \times \RR^m \rightarrow \big\{ \|p\| \le \delta^* \big\} \times \big\{ \|q\| \le \delta^* \big\}$ such that for each $(\zeta, \theta) \in \RR \times \RR^m$, $(p, q) = \rho(\zeta, \theta ; \epsilon )$ is the unique solution to the equation
\begin{equation} \label{eq:CFe}
(p, q) = \CF_{\epsilon}(p, q, \zeta, \theta) \,.
\end{equation}
Clearly, we have $\big\{ ( \rho(\zeta, \theta ; \epsilon ), \zeta, \theta ) : (\zeta, \theta) \in \RR \times \RR^m \big\} = W^+_{\epsilon} \textcap W^-_{\epsilon}$, which verifies that $\Me$ is the graph of the function $\rho(\cdot, \cdot ; \epsilon)$.

Next, we consider the restriction of $\CF_{\epsilon}$ on $\big\{ \|p\| \le 2 \tfrac{C_3}{\sigma} \epsilon^{\mu} \big\} \times \big\{ \|q\| \le 2 \tfrac{C_3}{\sigma} \epsilon^{\mu} \big\} \times \RR \times \RR^m$. By (\ref{ex3:CP:e-bound}) and (\ref{ex3:CQ:e-bound}),
we have that 
\begin{align*} 
( \bar\CP(q, \zeta, \theta ; \epsilon), \bar\CQ(p, \zeta, \theta ; \epsilon) )
\in \big\{ \|p\| \le 2 \tfrac{C_3}{\sigma} \epsilon^{\mu} \big\} \times \big\{ \|q\| \le 2 \tfrac{C_3}{\sigma} \epsilon^{\mu} \big\}
\end{align*}
for all $(p, q, \zeta, \theta) \in \big\{ \|p\| \le 2 \tfrac{C_3}{\sigma} \epsilon^{\mu} \big\} \times \big\{ \|q\| \le 2 \tfrac{C_3}{\sigma} \epsilon^{\mu} \big\} \times \RR \times \RR^m$. Thus, for each $(\zeta, \theta) \in \RR \times \RR^m$, $\CF_{\epsilon}(\cdot, \cdot, \zeta, \theta)$ is also a contraction on $\big\{ \|p\| \le 2 \tfrac{C_3}{\sigma} \epsilon^{\mu} \big\} \times \big\{ \|q\| \le 2 \tfrac{C_3}{\sigma} \epsilon^{\mu} \big\}$ under the norm $\|\cdot\|_1$, and the unique solution $(p, q)$ to (\ref{eq:CFe}) is in fact contained in the set $\big\{ \|p\| \le 2 \tfrac{C_3}{\sigma} \epsilon^{\mu} \big\} \times \big\{ \|q\| \le 2 \tfrac{C_3}{\sigma} \epsilon^{\mu} \big\}$. Therefore, $\rho(\cdot, \cdot ; \epsilon)$ maps $\RR \times \RR^m$ into $\big\{ \|p\| \le 2 \tfrac{C_3}{\sigma} \epsilon^{\mu} \big\} \times \big\{ \|q\| \le 2 \tfrac{C_3}{\sigma} \epsilon^{\mu} \big\}$.

Finally, we show that $\rho(\zeta, \theta ; \epsilon)$ is $C^r$ with respect to $(\zeta, \theta)$ on $\RR \times \RR^m$ using the implicit function theorem. Define a $C^r$ function $\CG_{\epsilon} : \big\{ \|p\| < \delta^* \big\} \times \big\{ \|q\| < \delta^* \big\} \times \RR \times \RR^m \rightarrow \RR^{n_u} \times \RR^{n_s}$ as follows:
\begin{align*}
\CG_{\epsilon} (p, q, \zeta, \theta) := 
\begin{pmatrix}
p - \CP(q, \zeta, \theta ; \epsilon) \\
q - \CQ(p, \zeta, \theta ; \epsilon)
\end{pmatrix} .
\end{align*}
Note that at any $(\rho(\zeta, \theta ; \epsilon ), \zeta, \theta) \in \big\{ \|p\| \le 2 \tfrac{C_3}{\sigma} \epsilon^{\mu} \big\} \times \big\{ \|q\| \le 2 \tfrac{C_3}{\sigma} \epsilon^{\mu} \big\} \times \RR \times \RR^m$,
\begin{align*}
D_{(p,q)} \CG_{\epsilon} = 
\begin{pmatrix}
I_{n_u \times n_u} & - D_q \CP \\
- D_p \CQ & I_{n_s \times n_s}
\end{pmatrix}
\end{align*}
is always invertible since $\| D_q \CP \| \le \tfrac{1}{2}$ by (\ref{WH:CP:Lip}) and $\| D_p \CQ \| \le \tfrac{1}{2}$ by (\ref{WH:CQ:Lip}). Then, by the implicit function theorem, the function $\rho(\cdot, \cdot ; \epsilon)$, which satisfies $\CG_{\epsilon}(\rho(\zeta, \theta ; \epsilon), \zeta, \theta) \equiv 0$ for all  $(\zeta, \theta) \in \RR \times \RR^m$, is $C^r$ on $\RR \times \RR^m$.
\end{proof}

Finally, we return to the original system (\ref{GWH}). By the periodicity of $\CP(q, \zeta, \theta ; \epsilon)$ and $\CQ(p, \zeta, \theta ; \epsilon)$ with respect to $\zeta$ and $\theta$, we have that $\rho(\zeta, \theta ; \epsilon)$ is $2\zeta_0$-periodic in $\zeta$ and $2\pi$-periodic in each component of $\theta$. Recall the changes of variables (\ref{w:a-b-zeta}) and $(a,b) = v(p,q)$. Then Lemma \ref{lem:WH:Me} implies that for any $\epsilon \in (0, \epsilon^*]$, (\ref{GWH}) has an invariant set 
\begin{equation*}
\Te := \big\{ (w, \theta) : w = \chi(\zeta) + \eta(\zeta)\, v( \rho(\zeta, \theta ; \epsilon) ),\, \zeta \in \RR (\Mod{ 2\zeta_0 }),\, \theta \in \Torus^m \big\} \,, 
\end{equation*}
which is contained in an $\CO( \epsilon^{\mu} )$-neighborhood of $\To$ since the function $\rho(\cdot, \cdot ; \epsilon)$ is $\CO( \epsilon^{\mu} )$. To show that $\Te$ can be parameterized by (\ref{Torus:para}), which uses the $\zeta_0$-periodic basis $\hat{n}(\zeta)$, we take a section $\Sigma(\zeta_1)$ at an arbitrary $\zeta_1 \in [0, \zeta_0)$ as follows:
\begin{equation*}
\Sigma(\zeta_1) := \big\{ (w, \theta) : w = \chi(\zeta_1) + \eta(\zeta_1)\, v( p, q ),\, \|p\| < \delta^*,\, \|q\| < \delta^* \big\} \,. 
\end{equation*}
Clearly, we have 
\begin{equation*}
R(\zeta_1) := \big\{ (w, \theta) : w = \chi(\zeta_1) + \eta(\zeta_1)\, v( \rho(\zeta_1, \theta ; \epsilon) ),\, \theta \in \Torus^m \big\} \subseteq \Te \textcap \Sigma(\zeta_1) \,. 
\end{equation*}
Suppose there exist a $\zeta_2 \in [0, 2\zeta_0)$ different from $\zeta_1$ and a $\theta_2 \in \Torus^m$ such that $(w_2, \theta_2) \in \Te \textcap \Sigma(\zeta_1)$ for $w_2 = \chi(\zeta_2) + \eta(\zeta_2)\, v( \rho(\zeta_2, \theta_2 ; \epsilon) )$. Since $(w_2, \theta_2) \in \Sigma(\zeta_1)$ and $\rho(\cdot, \cdot ; \epsilon)$ is $\CO( \epsilon^{\mu} )$, there exists $(\tilde{p}, \tilde{q})$ with both $\|\tilde{p}\|$ and $\|\tilde{q}\|$ being $\CO( \epsilon^{\mu} )$ such that $w_2 = \chi(\zeta_1) + \eta(\zeta_1)\, v( \tilde{p}, \tilde{q} )$. Furthermore, since $(w_2, \theta_2) \in \Te$ and $\Te$ is contained in an $\CO( \epsilon^{\mu} )$-neighborhood of $\To$, the solution trajectory of (\ref{GWH}) that passes through $(w_2, \theta_2)$ is contained inside an $\CO( \epsilon^{\mu} )$-neighborhood of $\To$ forever in both forward time and backward time due to the invariance of $\Te$ under the flow of (\ref{GWH}). Then for the normal form (\ref{WHNF}), the solution trajectory that passes through $(\tilde{p}, \tilde{q}, \zeta_1, \theta_2)$ stays inside $\big\{ \|p\| < \delta^* \big\} \times \big\{ \|q\| < \delta^* \big\} \times \RR \times \RR^m$ forever in both forward time and backward time. Recall that $\Me$ is the largest invariant subset of $\big\{ \|p\| < \delta^* \big\} \times \big\{ \|q\| < \delta^* \big\} \times \RR \times \RR^m$. Thus $(\tilde{p}, \tilde{q}, \zeta_1, \theta_2) \in \Me$, which implies that $(\tilde{p}, \tilde{q}) = \rho(\zeta_1, \theta_2 ; \epsilon)$. It follows that $(w_2, \theta_2) \in R(\zeta_1)$. Then we have $R(\zeta_1) = \Te \textcap \Sigma(\zeta_1)$. Note that $\Te \textcap \Sigma(\zeta_1)$ is $C^r$ diffeomorphic to the intersection of $\Te$ and the section $\big\{ (w, \theta) : w = \chi(\zeta) + \hat{n}(\zeta)\,\xi,\, \xi \in \RR^{n-1},\, \|\xi\| < \Delta_0 \big\}$ for a certain $\Delta_0 > 0$ and any $\zeta \in [0, \zeta_0)$. Therefore, $\Te$ can be parameterized by (\ref{Torus:para}), and it is the unique invariant torus for (\ref{GWH}) inside an $\CO(1)$-neighborhood (i.e., independent of $\epsilon$) of $\To$.

\appendix

\section{The {\Wa} Principle} \label{appendix:Wa}

We follow the presentation of Conley \cite{Co78}. Let $\CX$ be a topological space and $\varphi : \RR \times \CX \rightarrow \CX$ be a flow. For a set $W \subset \CX$, we define the following sets:
\begin{align*}
\Wo :={}& \big\{ x \in W : \exists\, t > 0 \mbox{ such that } \varphi (t,x) \not\in W \big\} \,, \\ 
\Wm :={}& \big\{ x \in W : \varphi ([0,t),x) \not\subseteq W \mbox{ for all } t > 0 \big\} \,, 
\end{align*}
where $\Wo$ is the set of points that do not stay in $W$ forever under the flow $\varphi$ in forward time, and $\Wm$ is the set of points that immediately leave $W$ in forward time. Clearly, $\Wm \subseteq \Wo \subseteq W$.

\begin{defnWa} The set $W$ is called a {\Wa} set if the following conditions are satisfied:
\begin{enumerate}
    \item[(W1)] If $x \in W$ and $\varphi([0,t],x) \subset {\rm{cl}}(W)$ then $\varphi([0,t],x) \subset W$.
    \item[(W2)] $\Wm$ is closed relative to $\Wo$.
\end{enumerate}
\end{defnWa}

\begin{Waze} If\/ $W$ is a {\Wa} set, then\/ $\Wm$ is a strong deformation retract of\/ $\Wo$ and\/ $\Wo$ is open relative to\/ $W$.
\end{Waze}

The proof of the above theorem can be found in \cite{Co78,KaMiMr04}. Here we only recall that when $\Wm$ is a strong deformation retract of $\Wo$, there exists a continuous function $r : \Wo \times [0,1] \rightarrow \Wo$ such that: (1) for all $x \in \Wo$, $r(x,0) = x$ and $r(x,1) \in \Wm$; and (2) for all $x \in \Wm$ and all $\sigma \in [0,1]$, $r(x,\sigma) = x$. The function $r$ is called a {\it strong deformation retraction}.

\bibliographystyle{siam}
\bibliography{IMrefs}

\end{document}